%% file: marginal_2026.tex
\documentclass[11pt,twoside]{article}

\newcommand{\ifims}[2]{#1} 
\newcommand{\ifAMS}[2]{#1}   
\newcommand{\ifau}[4]{#1}  
\newcommand{\ifbook}[2]{#1}   
\newcommand{\ifapp}[2]{#2}  
\newcommand{\iffourG}[2]{#2}  

\input myfrontmatter

\input mydef

\input statdef

\def\AFN{\mathbbmsl{U}}
\def\IFN{\IF}

\def\thetitle{Convergence of alternating minimization and sup-norm accuracy bounds in perturbed optimization}
\def\theruntitle{Alternating minimization and sup-norm bounds}

\def\theabstract{
Let the objective function \( \fs \) depends on the target variable \( \targv \) along with a 
nuisance variable \( \nuiv \): \( \fs(\upsv) = \fs(\targv,\nuiv) \).
Consider the task of identifying the marginal solution \( \targvs = \arginf_{\targv} \inf_{\nuiv} \fs(\targv,\nuiv) \).
This paper discusses three related problems.
The \emph{plug-in} approach, widely used, e.g., in inverse problems, suggests using a preliminary guess (pilot)
\( \hat{\nuiv} \) and apply the solution of the partial optimization 
\( \hat{\targv} = \arginf_{\targv} \fs(\targv,\hat{\nuiv}) \).
The main question to address within this approach is the required quality of the pilot, ensuring the prescribed 
accuracy of \( \hat{\targv} \).
The popular \emph{alternating optimization} approach suggests the following procedure:
given a starting guess \( \targv_{0} \), for \( t \geq 1 \), define
\(
	\nuiv_{t} 
	=
	\arginf_{\nuiv} \fs(\targv_{t-1},\nuiv)
\), and then 
\(
	\targv_{t} 
	=
	\arginf_{\targv} \fs(\targv,\nuiv_{t})
\).
The main question here is the set of conditions ensuring a convergence of \( \targv_{t} \) to \( \targvs \).
Finally, the paper discusses an interesting connection between marginal optimization and \emph{sup-norm estimation}.
The basic idea is to consider one component of the variable \( \upsv \) as a target and the rest as nuisance.
In all cases, we provide accurate closed-form results under realistic assumptions.
The results are illustrated by an example for Bradley--Terry--Luce model
of ranking from pairwise comparisons. 
}

\def\kwdp{62F10,62E17}
\def\kwds{62J12}

\def\thekeywords{partial optimization; semiparametric bias; expansions using perturbed optimization}

\def\thankstitle{The article was prepared within the framework of the HSE University Basic Research Program}

\aua
{Vladimir Spokoiny}
{Spokoiny, V. }
{
    Weierstrass Institute Berlin,  
    HSE University 
    and IITP RAS, Moscow,
    \\
    Mohrenstr. 39, 10117 Berlin, Germany
    }
{spokoiny@wias-berlin.de}
{Weierstrass-Institute and Humboldt University Berlin}
{
\thankstitle
}

\begin{document}
\thispagestyle{empty}
{
\title{\thetitle}
\theauthors

\maketitle
\begin{abstract}
{\footnotesize \theabstract}
\end{abstract}

\ifAMS
    {\par\noindent\emph{AMS 2010 Subject Classification:} Primary \kwdp. Secondary \kwds}
    {\par\noindent\emph{JEL codes}: \kwdp}

} 

\newpage
\tableofcontents

\Chapter{Introduction}
\label{Spertintr}

For an objective  function \( \fs(\cdot) \), consider an unconstrained optimization problem
\begin{EQA}
	\upsvs
	&=&
	\arginf_{\upsv} \fs(\upsv) \, .
\label{8cicjwe6tdkcodlds}
\end{EQA}
Perturbed optimization addresses the following question: how 
do \( \upsvs \) and \( \fs(\upsvs) \) change if the objective function 
\( \fs \) is perturbed in some special way?
This question originates from statistical and machine learning applications,
with many particular examples and applications
including uncertainty quantification, accuracy guarantees and rates for statistical learning,
bias due to penalizations, etc. 
\cite{Sp2024PO} established several rather precise results in the case of a smooth convex function \( \fs(\upsv) \);
%
see Section~\ref{Squadnquad} for a brief overview.
%
This paper studies a slightly different but related question.
Suppose that the full variable \( \upsv \) is composed of two parts:
\( \upsv = (\targv,\nuiv) \), where \( \targv \) is of primary interest
while \( \nuiv \) is a nuisance variable. 
Of course, a global optimization 
\begin{EQA}[c]
	\upsvs
	=
	(\targv,\nuivs)
	=
	\arginf_{(\targv,\nuiv)} \fs(\targv,\nuiv)
\end{EQA}
delivers the \( \targv \)-component of the solution.
However, a joint optimization can be a hard problem.
Instead, one can consider two partial problems
\begin{EQA}[rcccl]
	\targv_{\nuiv}
	&=&
	\arginf_{\targv} \fs(\targv,\nuiv)
	\, ,
	\qquad
	\nuiv_{\targv}
	&=&
	\arginf_{\nuiv} \fs(\targv,\nuiv)
	\, ,
\label{sugw725jvliw6ydfew}
\end{EQA}
when one component is fixed while the other is optimized. 
Each of these two problems can be much easier than the original 
full-dimensional optimization, especially if the objective function
\( \fs(\targv,\nuiv) \) is convex in \( \targv \) with \( \nuiv \) fixed
and another way around.
We mention three different setups where partial optimization 
appers in a natural way.

\medskip
\par\textbf{Setup 1: Plug-in approach in semiparametric statistical inference.}
A pragmatic plug-in procedure is frequently used in 
inverse problems or statistical semiparametric estimation.
Let the fidelity function \( \LL(\cdot) \) (loss, empirical risk, negative log-likelihood, contrast, etc.)
depend on the parameter of interest \( \targv \) (often low-dimensional)
and on the nuisance parameter \( \nuiv \) (usually high-dimensional or even functional).
A \emph{profile} estimation procedure suggests estimating the full parameter
\( \upsv = (\targv,\nuiv) \) and then using its \( \targv \)-component.
This estimator possesses many nice theoretical properties like asymptotic efficiency and normality; see e.g.
\cite{MR1245941}, \cite{robinson1988root}, \cite{Do1994}, \cite{van2000asymptotic}, \cite{ChCh2018}.
Unfortunately, full-dimensional optimization is often hard to implement and study.
At the same time, partial optimization of \( \LL(\targv,\nuiv) \) w.r.t. \( \targv \) for \( \nuiv \) fixed can be a much simpler problem.
This suggests a \emph{plug-in} approach: given a preliminary/pilot estimate \( \hat{\nuiv} \), define \( \hat{\targv} \) by maximization 
of \( \LL(\targv,\hat{\nuiv}) \):
\begin{EQA}
	\hat{\targv}
	& \eqdef &
	\arginf_{\targv} \LL(\targv,\hat{\nuiv}) 
\label{dhyiu3iwextq4e3esarg}
\end{EQA}
with a starting guess (pilot) \( \hat{\nuiv} \), apply \( \targv_{\hat{\nuiv}} \) 
in place of \( \targvs \).
We refer to 
\cite{NoPuSp2024}, 
\cite{baLoLu2020}, 
\cite{ChMo2022}, 
for an overview and further references 
related to the plugin approach in random design regression.
A more general Error-in-Operator problem has been discussed in \cite{Sp2025EiO}.
The error of this plug-in method can be studied within 
perturbed optimization problem: 
the true objective function \( \fs(\targv,\nuivs) \) of the target variable
\( \targv \) is replaced with \( \fs(\targv,\hat{\nuiv}) \).
We study the error \( \targv_{\hat{\nuiv}} - \targvs \) induced by this perturbation in Section~\ref{Ssemiopt} ahead
and show that this error is essentially linear in \( \hat{\nuiv} - \nuivs \).

\medskip
\par\textbf{Setup 2: Alternating optimization and EM-type procedures.} \,
The other widely used approach is a so-called \emph{alternating optimization} (AO).
Suppose to be given a starting guess \( \targv_{0} \) and consider the following AO procedure: 
for \( t \geq 1 \)
\begin{EQ}[rcccl]
	\nuiv_{t} 
	&=&
	\arginf_{\nuiv} \fs(\targv_{t-1},\nuiv)
	\, ,
	\qquad
	\targv_{t} 
	&=&
	\arginf_{\targv} \fs(\targv,\nuiv_{t})
	\, .
\label{tdydyc6cfw5cgedhrwfxxI}
\end{EQ}
Here, we assume for simplicity of presentation that each partial solution is a singleton.
This is obviously the case when \( \fs \) is strongly convex in \( \targv \) for \( \nuiv \) fixed
and another way around.
A vast body of literature exists on this subject.
A comprehensive survey on coordinate descent and its variants, including alternating optimization can be found e.g. in 
\cite{wright2015coordinate},
%
\cite{boyd2004convex},
and in
\cite{bauschke2017convex}.
For early theoretical analysis of AO convergence for convex problems, see 
\cite{luo1992convergence}.
%
\cite{bezdek2003convergence}
discusses convergence properties of AO in the context of clustering.
\cite{tseng2001convergence}
establishes convergence guarantees for non-smooth problems.
Block coordinate descent methods in the broader context of optimization
is addressed in 
\cite{bertsekas1999nonlinear}.

We also mention some 
applications in machine learning and statistics.
\cite{xu2013block}
extends the AO theory to multi-convex problems.
%
Non-negative matrix factorization by AO has been studied in
\cite{hsieh2011fast}.
%
\cite{mairal2015incremental}
discussed AO in the context of majorization-minimization (MM) algorithms.
%
\cite{beck2013convergence}
provides general convergence guarantee bounds for AO.
%
General convergence theory of AO for non-smooth and non-convex problems 
has been studied in
\cite{razaviyayn2013unified}
%
and in 
\cite{xu2017globally}.

A common point for all the obtained results is local convergence:
an AO procedure converges at a linear rate provided a good starting guess \( \targv_{0} \).
This seems to be an intrinsic feature of the approach.
However, for practical application, one would need more specific conditions
on this local vicinity, its shape, and dependence 
on the curvature of the objective function \( \fs \) and its smoothness properties.
Section~\ref{SoptAO} addresses these questions and 
states an accurate result about linear convergence 
of the AO procedure to the solution \( (\targvs,\nuivs) \).

\medskip
\par\textbf{Setup 3: Optimization with sup-norm errors.} \,
In many applied areas, one would be interested in the componentwise accuracy of the solution \( \hat{\upsv} \).
This naturally leads to sup-norm error \( \| \hat{\upsv} - \upsvs \|_{\infty} \).
In this paper, we consider one typical example from statistical inference:
\emph{ranking from pairwise comparison} problems, when 
players (goods, assets) are ranked using binary pairwise comparisons.
The Bradley-Terry-Luce (BTL) model (\cite{BT1952}, \cite{Lu1959}) 
assumes that the outcome of each comparison between \( j \)th and \( m \)th palyers is a Bernoulli 
experiment with the canonical parameter depending on the difference \( \ups_{j} - \ups_{m} \)
of individual skills (merits) of each player.  
A detailed historical overview and the state-of-the-art results for the BTL model can be found in \cite{GSZ2023}.
Existing results confirm that the ranking problem leads to sup-norm estimation, and a careful treatment of this problem is rather involved.
The methods of analysis, proof strategy, and the final results differ significantly from the case of a more standard \( \ell_{2} \)-norm.
\medskip

\par\textbf{Prior works.} \,
A recent paper \cite{Sp2024PO} studied the problem of perturbed optimization 
for a smooth convex function \( \fs \) and its linear/quadratic/smooth additive perturbation.
The obtained results describe the change of the solution and the value of optimization problem \eqref{8cicjwe6tdkcodlds} 
if the objective function \( \fs(\upsv) \) is corrupted by some perturbation. 
Applications to statistical inference can be found e.g. in \cite{Sp2024,Sp2025EiO}.
Consider for simplicity a linear perturbation \( \langle \Av, \upsv \rangle \), that is,
\( \fp(\upsv) = \fs(\upsv) + \langle \Av, \upsv \rangle \).
For \( \IFN \eqdef \nabla^{2} \fs(\upsvs) \) and \( \upsvp \eqdef \arginf_{\upsv} \fp(\upsv) \),
the main results  describe the differences \( \upsvp - \upsv \) and \( \fp(\upsvp) - \fp(\upsvs) \)
in a self-consistent (simplified) form:
\begin{EQA}[rcl]
\label{ufufdj47ujkwfjfee}
	\| \IFN^{1/2} (\upsvp - \upsvs) \|
	& \leq &
	\frac{3}{2} \| \IFN^{-1/2} \Av \|
	\, ,
	\\
    \| \IFN^{1/2} (\upsvp - \upsvs + \IFN^{-1} \Av) \|
    & \leq &
    \frac{3\dltwu_{3}}{4} \| \IFN^{-1/2} \Av \|^{2} 
	\, ,
\label{DGttGtsGDGm13rGa2i}
    \\
    \Bigl| 2 \fp(\upsvp) - 2 \fp(\upsvs) + \| \IFN^{-1/2} \Av \|^{2} \Bigr|
    & \leq &
    \frac{\dltwu_{3}}{2} \, \| \IFN^{-1/2} \Av \|^{3} \, .
    \qquad
\label{3d3Af12DGttGa2i}
\end{EQA}
Here \( \dltwu_{3} \) is a small constant entering the third-order smoothness condition \nameref{LLsT3ref}.
Section~\ref{Squadnquad} ahead provides a more detailed overview and the related conditions.
Apart convexity and smoothness of \( \fs(\cdot) \), 
it is required that \( \dltwu_{3} \, \| \IFN^{-1/2} \Av \| \ll 1 \). 
This condition ensures the concentration result \eqref{ufufdj47ujkwfjfee} and 
the remainders in the right-hand side of \eqref{DGttGtsGDGm13rGa2i} and \eqref{3d3Af12DGttGa2i}
are smaller than the corresponding leading terms of the presented expansions.
For statistical inference, this condition \( \dltwu_{3} \, \| \IFN^{-1/2} \Av \| \) can be viewed as
so-called \emph{critical dimension} condition: 
``effective parameter dimension is smaller than effective sample size'',
where the \emph{effective dimension} \( \dimA \) can be defined as an upper bound on the perturbation energy 
\( \| \IFN^{-1/2} \Av \|^{2} \),
while the effective sample size \( \neff \) corresponds to the smallest eigenvalue of the Fisher information matrix \( \IFN \).
Statistical inference on the parameter of interest requires even stronger condition \( \dimA^{2} \ll \neff \).
General results on Laplace approximation from \cite{katsevich2024laplaceapproximationaccuracyhigh} can be viewed as a lower bound in this problem
and indicate that this condition is really unavoidable. 
Surprisingly, \cite{GSZ2023} established very strong inference results 
for the Bradley-Terry-Luce (BTL) problem of ranking of \( \dimp \) players 
from pairwise comparisons under a much weaker condition \( \neff \gg \log(\dimp) \).
The distinguishing feature of the study in \cite{GSZ2023} is that 
each component of the parameter vector is analyzed separately, leading to sup-norm losses. 
One of the aims of this paper is to revisit the estimation problem with sup-norm losses; 
see Section~\ref{Ssupnormco} ahead.

\medskip
\par\textbf{This paper contribution.} \,
This paper provides a general setup of partial optimization; the obtained results
are used for two specific problems: 
alternating optimization and sup-norm estimation.
Below, we briefly address each of these contributions.

\par\textbf{(1) Partial optimization:} \,
It is well known that for a smooth function \( \fs(\targv,\nuiv) \), the partial solution 
\( \targvs_{\nuiv} = \arginf_{\targv} \fs(\targv,\nuiv) \) is nearly linear in \( \nuiv \) in the form
\begin{EQA}[c]
	\targv_{\nuiv} - \targvs
	\approx
	- \IFT_{\targv\targv}^{-1} \, \IFT_{\targv\nuiv} (\nuiv - \nuivs)
	\, .
\end{EQA}
Here \( \IFT_{\targv\targv} , \, \IFT_{\targv\nuiv}, \, \IFT_{\nuiv\targv}, \, \IFT_{\nuiv\nuiv} \) are the blocks of the full-dimensional
information matrix \( \IFT \eqdef \nabla^{2} \fs(\targv,\nuiv) \).
However, accurate guarantees for this approximation are important.
One of the main results of this paper 
can be written (in a simplified form) as 
\begin{EQA}[c]
	\bigl\| 
		\IFT_{\targv\targv}^{1/2} 
		\bigl\{ \targv_{\nuiv} - \targvs 
		+ \IFT_{\targv\targv}^{-1} \IFT_{\targv\nuiv} (\nuiv - \nuivs) \bigr\} 
	\bigr\|
	\leq 
	\dltwuns \, \| \IFT_{\nuiv\nuiv}^{1/2} (\nuiv - \nuivs) \|_{\nano}^{2}
	\, ,
\label{sdc7eujkkkkd5w3ecddw}
\end{EQA}
where \( \dltwuns \) is a small constant related to the smoothness 
condition \nameref{LLpT3ref} on the function \( \fs \) and \( \| \cdot \|_{\nano} \)
is a norm in \( \R^{\dimq} \) in which smoothness of \( \fs \) w.r.t. \( \nuiv \) is measured;
see Section~\ref{Ssemiopt} ahead for more details.
The closed self-consistent form of the remainder is essential.
In particular, it helps for studying many iterative algorithms.
Later, we apply this bound to two special cases.

\par\textbf{(2) Alternating optimization:} \,
Suppose to be given a starting guess \( \targv_{0} \) and consider
the following alternating optimization (AO) procedure: for \( t \geq 1 \)
\begin{EQ}[rcccl]
	\nuiv_{t} 
	&=&
	\arginf_{\nuiv} \fs(\targv_{t-1},\nuiv)
	\, ,
	\qquad
	\targv_{t} 
	&=&
	\arginf_{\targv} \fs(\targv,\nuiv_{t})
	\, .
\label{tdydyc6cfw5cgedhrwfxxi}
\end{EQ}
The question under study is a set of conditions ensuring a convergence 
of this procedure to the solution \( (\targvs,\nuivs) \) as \( t \) grows.
Let \( \IFT(\upsv) = \nabla^{2} \fs(\upsv) \) and \( \IFT = \IFT(\upsvs) \)
with the blocks \( \IFT_{\targv\targv} , \, \IFT_{\targv\nuiv}, \, \IFT_{\nuiv\targv}, \, \IFT_{\nuiv\nuiv} \).
Assume that each diagonal block \( \IFT_{\targv\targv} \) and \( \IFT_{\nuiv\nuiv} \) is strongly positive.
Define
\begin{EQA}[c]
	\IFC \eqdef 
	\IFT_{\targv\targv}^{-1/2} \, \IFT_{\targv\nuiv} \, \IFT_{\nuiv\nuiv}^{-1/2} .
\label{jccybr9vgoey4rdy3d}
\end{EQA}
Applicability of the AO method and its convergence at a linear rate requires the condition \( \| \IFC \, \IFC^{\T} \| < 1 \);
see Section~\ref{SoptAO}.

\par\textbf{(3) Sup-norm perturbed optimization:} \,
Let \( \fs(\upsv) \) be a strongly convex function and \( \upsvs = \arginf_{\upsv} \fs(\upsv) \).
Define \( \IFT = \nabla^{2} \fs(\upsvs) = (\IFT_{jm}) \).
Let also \( \fp(\upsv) = \fs(\upsv) + \langle \AAv,\upsv \rangle \)
and \( \upsvp = \arginf_{\upsv} \fp(\upsv) \).
We intend to bound the corresponding change \( \upsvp - \upsvs \) in a sup-norm.
The basic idea is to consider the sup-norm estimation 
as a kind of linewise optimization with one fixed component \( \nui_{j} \) as a target
and the remaining components \( \nui_{m} \) as nuisance. 
Fix a diagonal matrix \( \DF \eqdef \diag(\DPN_{1}, \ldots, \DPN_{\dimp}) \) with \( \DPN_{j}^{2} = \IFT_{jj} \)
and define \( \BBF \eqdef \DF^{-1} \IFT \, \DF^{-1} \).
Introduce also the value \( \crosssup \) and \( \crossB \):
\begin{EQA}[rcccl]
	\crosssup
	& \eqdef &
	\sup_{\| \uv \|_{\infty} \leq 1} \| (\BBF - \Id_{\dimp}) \uv \|_{\infty}
	&=&
	\| \BBF - \Id_{\dimp} \|_{1}
	=
	\max_{j = 1,\ldots,\dimp} \frac{1}{\DPN_{j}}
	\sum_{m \neq j} \frac{|\IFT_{jm}|}{\DPN_{m}}
	\, ,
	\qquad
\label{7tdsyf8iuwopkrtg4576i}
	\\
	\crossB
	& \eqdef &
	\sup_{\| \uv \|_{\infty} \leq 1} \| \BBF^{-1} \uv \|_{\infty} 
	&=&
	\| \BBF^{-1} \|_{1}
	\, .
\label{hdcukekei76erddfreei}
\end{EQA}
The results presented in Section~\ref{Ssupnormco} ahead required 
\( \crosssup < 1 \) and \( \crossB \) separated away from zero.
This particularly means that the information matrix \( \IFT \) has 
a stochastically dominant structure, and each diagonal element 
\( \BBF_{jj} = 1 \) exceeds the sum \( \sum_{m \neq j} |\BBF_{jm}| \).
This enables us 
to derive self-consistent sup-norm bounds in perturbed optimization:
\begin{EQA}[ccl]
\label{usdhyw6hikhurnetrspp}
	\bigl\| \DF^{-1} \bigl\{ \IFT (\upsvd - \upsvs) + \AAv \bigr\} \bigr\|_{\infty}
	& \leq &
	\dltwu_{\CONSTzz} \, \| \DF^{-1} \AAv \|_{\infty}^{2}
	\, ,
	\\
\label{usdhyw6hikhurnetrsppDi}
	\| \DF (\upsvd - \upsvs) + \DF^{-1} \AAv \|_{\infty} 
	& \leq &
	\crossB^{-1} \bigl( \| (\BBF - \Id_{\dimp}) \, \DF^{-1} \AAv \|_{\infty} + \dltwu_{\CONSTzz} \, \| \DF^{-1} \AAv \|_{\infty}^{2} \bigr)
	\, ,
	\qquad
\end{EQA}
where \( \dltwu_{\CONSTzz} \) is an explicit small constant depending on local smoothness of \( \fs \);
see Section~\ref{Ssupnormco} for a precise statement.

\medskip
\par\textbf{Organization of the paper.} \,
Section~\ref{Ssemiopt} studies the problem of partial and marginal optimization
with a particular focus on the semiparametric bias caused by using 
a wrong value of the nuisance parameter.
The obtained results will be used in Section~\ref{SoptAO} for 
proving the convergence of the AO procedure.
Section~\ref{Ssupnormco} applies the bound from Section~\ref{Ssemiopt}
to the problem of perturbed optimization in the sup-norm. 
Section~\ref{ScoBTL} presents a numerical example for the Bradley-Terry-Luce (BTL) model.
Section~\ref{Squadnquad} overviews one useful result from \cite{Sp2024PO} 
on linearly perturbed optimization, which will be used in all our derivations.
The proofs are collected at the end of each section.
%


{
\renewcommand{\Section}[1]{\section{#1}}
\renewcommand{\Subsection}[1]{\subsection{#1}}

\input marg_opt

	\subsection{Partial optimization. Proofs}
	Here we collect the proofs of Proposition~\ref{PbiassemiN}
	and Theorem~\ref{PsemiAvex} and Theorem~\ref{Pconcsupp}.
	\input marg_opt_proofs

\input marg_opt_alt

\input alternate_proof.tex

\input marg_opt_sup

	\subsection{Sup-norm optimization. Proofs}
	Here we collect the proofs of the results in this section.
\input marg_opt_sup_proofs

}
\section*{Funding and competing interests}
The article was prepared within the framework of the HSE University Basic Research Program.

\bibliographystyle{apalike}

\bibliography{exp_ts,listpubm-with-url}

\newpage
\appendix
{   
\input BTL_short

\input BTL_short_proofs

\input localbounds_short

}



\end{document}

%% file: myfrontmatter.tex
\usepackage{ifthen}
\usepackage[ampersand]{easylist}
\usepackage{amsmath,amssymb,amsthm}
\usepackage[bb=stixtwo]{mathalpha}
\usepackage{mathtools}
\usepackage{natbib}
\usepackage{epsfig,graphicx}
\usepackage{comment}
\usepackage{color}
\usepackage{srcltx}
\usepackage[mathscr]{eucal}
\usepackage[math]{easyeqn}
\usepackage{etoolbox}
\usepackage{tocloft}
\usepackage{nccmath}
\usepackage{hyperref}
\hypersetup{
            colorlinks,
            linkcolor=hookersgreen,
            linktoc=blue,
            citecolor=ultramarine,
            urlcolor=black,
            filecolor=black
            }
\usepackage{nameref}
\usepackage{multirow}

\usepackage{accents}
\usepackage{mathrsfs}
\usepackage{bbm}
\usepackage{algorithm}
\usepackage{algpseudocode}

\ifims{
\textheight=22cm
\textwidth=14.8cm
\topmargin=0pt
\oddsidemargin=1.0cm
\evensidemargin=1.0cm
\linespread{1.3}

}{ 
}

\numberwithin{equation}{section}
\numberwithin{figure}{section}
\newcounter{example}[section]
\numberwithin{example}{section}
\newcounter{remark}[section]
\numberwithin{remark}{section}
\newtheorem{theorem}{Theorem}[section]
\newtheorem{proposition}[theorem]{Proposition}
\newtheorem{lemma}[theorem]{Lemma}

\newtheorem{exmp}[example]{Example}
\newtheorem{rmrk}[remark]{Remark}
\newenvironment{example}{\begin{exmp}\rm}{\end{exmp}}
\newenvironment{remark}{\begin{rmrk}\rm}{\end{rmrk}}

\ifbook{
    \newcommand{\Chapter}[1]{\section{#1}}
    \newcommand{\Section}[1]{\subsection{#1}}
    \newcommand{\Subsection}[1]{\subsubsection{#1}}

  }
  {
    \newcommand{\Chapter}[1]{\chapter{#1}}
    \newcommand{\Section}[1]{\section{#1}}
    \newcommand{\Subsection}[1]{\subsection{#1}}

  }

%% file: mydef.tex
\renewcommand{\(}{$\,}
\renewcommand{\)}{\,$}

\def\nquad{\hspace{-1cm}}

\def\eqdef{\stackrel{\operatorname{def}}{=}}

\DeclareMathAlphabet{\mathbbmsl}{U}{bbm}{bx}{sl}

\DeclareMathSymbol{\Alpha}{\mathalpha}{operators}{"41}
\DeclareMathSymbol{\Beta}{\mathalpha}{operators}{"42}
\DeclareMathSymbol{\Epsilon}{\mathalpha}{operators}{"45}
\DeclareMathSymbol{\Zeta}{\mathalpha}{operators}{"5A}
\DeclareMathSymbol{\Eta}{\mathalpha}{operators}{"48}
\DeclareMathSymbol{\Iota}{\mathalpha}{operators}{"49}
\DeclareMathSymbol{\Kappa}{\mathalpha}{operators}{"4B}
\DeclareMathSymbol{\Mu}{\mathalpha}{operators}{"4D}
\DeclareMathSymbol{\Nu}{\mathalpha}{operators}{"4E}
\DeclareMathSymbol{\Omicron}{\mathalpha}{operators}{"4F}
\DeclareMathSymbol{\Rho}{\mathalpha}{operators}{"50}
\DeclareMathSymbol{\Tau}{\mathalpha}{operators}{"54}
\DeclareMathSymbol{\Chi}{\mathalpha}{operators}{"58}
\DeclareMathSymbol{\omicron}{\mathord}{letters}{"6F}

\newcommand{\cc}[1]{\mathscr{#1}}
\newcommand{\bb}[1]{\boldsymbol{#1}}

\DeclareFontFamily{U}{mathx}{\hyphenchar\font45}
\DeclareFontShape{U}{mathx}{m}{n}{
<5><6><7><8><9><10>
<10.95><12><14.4><17.28><20.74><24.88>
mathx10
}{}
\DeclareSymbolFont{mathx}{U}{mathx}{m}{n}
\DeclareFontSubstitution{U}{mathx}{m}{n}
\DeclareMathAccent{\widebar}{0}{mathx}{"73}
\DeclareMathAccent{\barlong}{0}{mathx}{"74}

\renewcommand{\bar}[1]{\widebar{#1}}
\renewcommand{\hat}[1]{\widehat{#1}}
\renewcommand{\tilde}[1]{\widetilde{#1}}

\makeatletter
\def\mathcenterto#1#2{\mathclap{\phantom{#1}\mathclap{#2}}\phantom{#1}}
\let\old@widetilde\widetilde
\def\widetildeto#1#2{\mathcenterto{#2}{\old@widetilde{\mathcenterto{#1}{#2\,}}}}
\let\old@widehat\widehat
\def\widehatto#1#2{\mathcenterto{#2}{\old@widehat{\mathcenterto{#1}{#2\,}}}}
\makeatother

\newcommand{\thankstitle}[1]{\ifthenelse{\equal{#1}{}}{}{\thanks{#1}}}
\newcommand{\thanksau}[1]{\ifthenelse{\equal{#1}{}}{}{\thanks{#1}}}

\newcommand{\aua}[6]
{\def\authora{#1}
\def\runauthora{#2}
\def\addressa{#3}
\def\emaila{#4}
\def\affiliationa{#5}
\def\thanksa{#6}}

\def\theauthors{
\ifau{ 
  \author{
    \authora
    \thanksau{\thanksa}
    \\[5.pt]
    \addressa \\
    \texttt{ \emaila}
  }
}
{  
  \author{
    \authora
    \thanksau{\thanksa}
    \\[5.pt]
    \addressa \\
    \texttt{ \emaila}
    \and
    \authorb
    \thanksau{\thanksb}
    \\[5.pt]
    \addressb \\
    \texttt{ \emailb}
  }
}
{   
  \author{
    \authora
    \thanksau{\thanksa}
    \\[5.pt]
    \addressa \\
    \texttt{ \emaila}
    \and
    \authorb
    \thanksau{\thanksb}
    \\[5.pt]
    \addressb \\
    \texttt{ \emailb}
    \and
    \authorc
    \thanksau{\thanksc}
    \\[5.pt]
    \addressc \\
    \texttt{ \emailc}
  }
} {   
  \author{
    \authora
    \thanksau{\thanksa}
    \\[5.pt]
    \addressa \\
    \texttt{ \emaila}
    \and
    \authorb
    \thanksau{\thanksb}
    \\[5.pt]
    \addressb \\
    \texttt{ \emailb}
    \and
    \authorc
    \thanksau{\thanksc}
    \\[5.pt]
    \addressc \\
    \texttt{ \emailc}
    \and
    \authord
    \thanksau{\thanksd}
    \\[5.pt]
    \addressd \\
    \texttt{ \emaild}
  }
}
}

\renewcommand{\Gamma}{\varGamma}
\renewcommand{\Pi}{\varPi}
\renewcommand{\Sigma}{\varSigma}
\renewcommand{\Delta}{\varDelta}
\renewcommand{\Lambda}{\varLambda}
\renewcommand{\Psi}{\varPsi}
\renewcommand{\Phi}{\varPhi}
\renewcommand{\Theta}{\varTheta}
\renewcommand{\Omega}{\varOmega}
\renewcommand{\Xi}{\varXi}
\renewcommand{\Upsilon}{\varUpsilon}

\def\arginf{\operatornamewithlimits{arginf}}

\def\argmin{\operatornamewithlimits{argmin}}

\def\av{\bb{a}}

\def\ev{\bb{e}}

\def\uv{\bb{u}}

\def\wv{\bb{w}}

\def\zv{\bb{z}}

\def\Av{\bb{A}}

\def\Uv{\bb{U}}

\def\epsv{\bb{\varepsilon}}

\def\Psiv{\bb{\Psi}}

\usepackage{color}

\definecolor{blue(pigment)}{rgb}{0.2, 0.2, 0.6}
\definecolor{ultramarine}{rgb}{0.07, 0.04, 0.56}
\definecolor{darkspringgreen}{rgb}{0.09, 0.45, 0.27}
\definecolor{hookersgreen}{rgb}{0.0, 0.44, 0.0}
\definecolor{hgreen}{rgb}{0.0, 0.44, 0.0}
\definecolor{plum(traditional)}{rgb}{0.56, 0.27, 0.52}
\definecolor{purple(html/css)}{rgb}{0.5, 0.0, 0.5}
\definecolor{magenta(dye)}{rgb}{0.79, 0.08, 0.48}

%% file: statdef.tex
\def\bulletc{\par\noindent\( \bullet \) \,}

\def\fp{\fs^{\circ}}
\def\bvp{\bvn_{\circ}}
\def\upsvp{\upsv^{\circ}}
\def\Tensp{\Tens^{\circ}}

\def\AUv{\Delta}

\def\pert{t}
\def\Targ{\mathcal{X}}

\def\neff{\mathbbmsl{N}}
\def\hspm{\hspace{1pt}}

\def\inftyi{\scalebox{0.66}{$\infty$}}
\def\inftyii{\scalebox{0.44}{$\infty$}}

\def\nano{\circ}
\def\dual{*}
\def\duali{**}

\def\crossnorm{\rho_{\dual}}

\def\crossgrad{\rho_{\nano}}
\def\crosssup{\rho_{\dual}}
\def\crosssupp{\rho_{\duali}}
\def\crosssupt{\kappa_{2}}
\def\crossD{\rho_{\!_{\Delta}}}
\def\crossB{\nu}

\def\AFN{\mathbb{Z}}

\def\dmax{\kappa}

\def\AAv{{\mathcal{A}}}
\def\CCv{{\mathcal{C}}}

\def\bvn{\bb{\mu}}

\def\IFN{\IF}

\def\DFN{\DPN}

\def\AFN{\mathbb{U}}
\def\HFN{\HPN}
\def\HFN{\mathbb{H}}

\def\Tens{\mathcal{T}}

\def\prmt{\ups}

\def\prmtv{\bb{\prmt}}

\def\prmtvs{\prmtv^{*}}

\def\targ{x}
\def\targv{\bb{\targ}}
\def\targvs{\targv^{*}}

\def\targvn{\targv^{\circ}}

\def\hL{h}

\def\Matr{\mathfrak{M}}
\def\Graph{\mathcal{G}}
\def\GrV{\mathcal{V}}
\def\GrE{\mathcal{E}}

\def\DeltaB{\mathcal{E}}

\def\Eta{\mathcal{H}}

\def\nbin{N}

\def\HPN{H}

\def\dltwb{\omega}

\def\dltwu{\tau}

\def\dltwun{\dltwu_{12}}
\def\dltwunn{\dltwu_{21}}
\def\dltwuns{\dltwu_{\nano}}

\def\dltwbun{\dltwb}

\def\DF{\mathcal{\dblk}}

\def\R{\mathbbmsl{R}}
\def\E{\mathbbmsl{E}}

\def\P{\mathbbmsl{P}}

\def\kappa{\varkappa}

\def\diag{\operatorname{diag}}

\def\Var{\operatorname{Var}}

\def\T{\top}

\def\ex{\mathrm{e}}

\def\Id{\mathbbmsl{I}}
\def\Ind{\operatorname{1}\hspace{-4.3pt}\operatorname{I}}

\def\alp{\alpha}

\def\alpG{\alp}

\def\alpv{\bb{\alpha}}

\def\AFN{\mathbbmsl{U}}
\def\Avm{\bb{M}}

\def\bias{\mathsf{b}}

\def\biasp{b}

\def\biasp{\bias_{\circ}}

\def\BBF{\cc{B}}

\def\cdens{\phi}

\def\CONSTi{\mathtt{C}}

\def\CONSTzz{\mathbb{c}}

\def\DPN{D}

\def\DF{\cc{D}}

\def\dimA{\mathbb{p}}

\def\dimttl{\bar{\dimp}}

\def\dimp{p}

\def\dimGrV{p}

\def\dimq{q}

\def\Eta{\cc{H}}

\def\fs{f}

\def\fd{f^{\circ}}

\def\fn{g}

\def\fG{f_{\GP}}

\def\gp{g}

\def\GP{G}

\def\IF{\mathbbmsl{F}}

\def\IFN{\IFL}

\def\IFtotal{F}
\def\IFT{\mathscr{\IFtotal}}

\def\IFC{\mathcal{P}}

\def\IFP{\IFN_{\circ}}

\def\IFmax{\IF_{\max}}
\def\IFmin{\IF_{\min}}
\def\IFrange{\kappa}

\def\Kappa{\cc{K}}

\def\LT{L}
\def\LGP{\LT_{\GP}}

\def\LL{\cc{L}}

\def\nup{\eta}
\def\nupv{\bb{\nup}}

\def\nui{s}
\def\Nui{\mathcal{S}}

\def\nuiv{\bb{\nui}}

\def\nuivs{\nuiv^{*}}

\def\QP{Q}

\def\rhos{\rho^{*}}

\def\rr{\mathtt{r}}

\def\rru{\rr_{\circ}}

\def\rrta{\rr_{\targv}}
\def\rrnu{\rr_{\nuiv}}

\def\rrinf{\rr_{\inftyi}}
\def\rrinfi{\rr_{\inftyii}}

\def\rrnano{\rr_{\nano}}

\def\rrn{\rr}

\def\Tau{T}

\def\ups{\upsilon}
\def\upsv{\bb{\ups}}

\def\upsvd{\upsv^{\circ}}
\def\upsvs{\upsv^{*}}

\def\upsvn{\upsvd}

\def\ups{\upsilon}
\def\upsv{\bb{\ups}}

\def\upsd{\ups^{\circ}}

\def\upss{\ups^{*}}

\def\Ups{\varUpsilon}
\def\Upsd{\Ups^{\circ}}

\def\wv{\bb{w}}

\def\WV{\mathcal{W}}

\def\xx{\mathtt{x}}

\def\zz{\mathfrak{z}}

\def\zzinf{\zz_{\infty}}

%% file: marg_opt.tex

\def\IFN{\IF}

\Section{Partial and marginal optimization}
\label{Ssemiopt}
This section discusses the problem of conditional/partial and marginal optimization. 
Consider a function \( \fs(\prmtv) \) of a parameter \( \prmtv \in \R^{\dimttl} \)
which can be represented as \( \prmtv = (\targv,\nuiv) \), where \( \targv \in \R^{\dimp} \)
is the target subvector while \( \nuiv \in \Nui \subseteq \R^{\dimq} \) is a nuisance variable,
\( \dimttl = \dimp+\dimq \).
Our goal is to study the target component \( \targvs \) of the solution to the optimization problem 
\( \prmtvs = (\targvs,\nuivs) = \arginf_{\prmtv \in \Ups} \fs(\prmtv) \)
with a local set \( \Ups = \{ \upsv = (\targv,\nuiv) \colon \nuiv \in \Nui \} \).

\Subsection{Partial optimization}
\label{Spartopt}
For a fixed value of the nuisance variable \( \nuiv \in \Nui \), consider 
\( \fs_{\nuiv}(\targv) = \fs(\targv,\nuiv) \) as a function of \( \targv \) only.
Below we assume that \( \fs_{\nuiv}(\targv) \) is strongly convex in \( \targv \) for any \( \nuiv \in \Nui \).
Define 
\begin{EQA}[rcl]
	\targv_{\nuiv}
	& \eqdef &
	\arginf_{\targv \colon (\targv,\nuiv) \in \Ups} \fs_{\nuiv}(\targv) \, .
\label{jfoiuy2wedfv7tr2qsdzxdfso}
\end{EQA}
Our goal is to describe variability of the partial solution \( \targv_{\nuiv} \) in \( \nuiv \) in terms of \( \targv_{\nuiv} - \targvs \)
and \( \fs(\prmtvs) - \fs_{\nuiv}(\targv_{\nuiv}) \).
Introduce 
\begin{EQ}[rcccc]
	\Av_{\nuiv}
	& \eqdef &
	\nabla \fs_{\nuiv}(\targvs)
	&=&
	\nabla_{\targv} \fs(\targvs,\nuiv) 
	\, ,
	\\
	\IFN_{\nuiv}
	& \eqdef &
	\nabla^{2} \fs_{\nuiv}(\targvs) 
	&=&
	\nabla_{\targv\targv}^{2} \fs(\targvs,\nuiv) 
	\, .
\label{ge8qwefygw3qytfyju8qfhbdso}
\end{EQ}
Now we introduce local smoothness condition of each function \( \fs_{\nuiv}(\cdot) \) around \( \targv_{\nuiv} \).
Let for any \( \nuiv \in \Nui \), some radius \( \rrn_{\nuiv} \) and 
a positive definite symmetric \( \dimp \)-matrix \( \DPN_{\nuiv} \) be fixed
to descibe the local vicinity \( \{ \targv \colon \| \DPN_{\nuiv} (\targv - \targv_{\nuiv}) \| \leq \rrn_{\nuiv} \} \).
Here \( \| \cdot \| \) is the Euclidean norm.
We assume \( \DPN_{\nuiv}^{2} \leq \dmax^{2} \, \IFN_{\nuiv} \) for some constant \( \dmax \) meaning 
\( \dmax^{2} \, \IFN_{\nuiv} - \DPN_{\nuiv}^{2} \) positive definite.

\begin{description}

    \item[\label{LLsoT3ref} \( \bb{(\mathcal{T}_{3|\nuiv}^{*})} \)]
      \emph{For \( \nuiv \in \Nui \), 
      it holds
\begin{EQA}
	\sup_{\uv \in \R^{\dimp} \colon \| \DPN_{\nuiv} \uv \| \leq \rrn_{\nuiv}}  \,\, \sup_{\zv \in \R^{\dimp}} \,\,
	\frac{\bigl| \langle \nabla^{3} \fs_{\nuiv}(\targv_{\nuiv} + \uv), \zv^{\otimes 3} \rangle \bigr|}
		 {\| \DPN_{\nuiv} \zv \|^{3}}
	& \leq &
	\dltwu_{3} \, .
\label{d6f53ye5vry4fddfgeyd}
\end{EQA}
}

\end{description}

Our first result describes the \emph{semiparametric bias} \( \targvs - \targv_{\nuiv} \) caused by using
the value \( \nuiv \) of the nuisance variable in place of \( \nuivs \).
The idea is to treat the function \( \fs_{\nuiv}(\targv) \) as a perturbation of the function \( \fs_{\nuivs}(\targv) \).

\begin{proposition}
\label{PbiassemiN}
Assume the following conditions:

\bulletc
The solution 
\( (\targvs,\nuivs) = \arginf_{\targv,\nuiv} \fs(\targv,\nuiv) \) 
exists and unique.

\bulletc
\( \forall\nuiv \in \Nui \), the function 
\( \fs_{\nuiv}(\targv) \) is strongly convex, 
\( \targv_{\nuiv} = \arginf_{\targv} \fs_{\nuiv}(\targv) \), 
\( \IFN_{\nuiv} = \nabla^{2} \fs_{\nuiv}(\targv_{\nuiv}) \),
and
\( \Av_{\nuiv} = \nabla \fs_{\nuiv}(\targvs) \).

\bulletc
\( \forall \nuiv \in \Nui \), 
condition \nameref{LLsoT3ref} holds at \( \targv_{\nuiv} \) with \( \DPN_{\nuiv}^{2} \), \( \rrn_{\nuiv} \), and \( \dltwu_{3} \) such that 
\begin{EQA}[c]
	\DPN_{\nuiv}^{2} \leq \dmax^{2} \, \IFN_{\nuiv} 
	\, ,
	\quad
	\rrn_{\nuiv} \geq \frac{3}{2} \| \DPN_{\nuiv} \, \IFN_{\nuiv}^{-1} \Av_{\nuiv} \| 
	\, ,
	\quad
	\dmax^{2} \dltwu_{3} \| \DPN_{\nuiv} \, \IFN_{\nuiv}^{-1} \Av_{\nuiv} \| < \frac{4}{9} 
	\, .
\label{8difiyfc54wrboer7bjfrO}
\end{EQA}
Then 
\begin{EQ}[rcccl]
    \| \DPN_{\nuiv} (\targv_{\nuiv} - \targvs) \| 
    & \leq &
    \frac{3}{2} \| \DPN_{\nuiv} \, \IFN_{\nuiv}^{-1} \Av_{\nuiv} \|
    \, ,
    &&
    \\
    \| \DPN_{\nuiv}^{-1} \IFN_{\nuiv} (\targv_{\nuiv} - \targvs + \IFN_{\nuiv}^{-1} \Av_{\nuiv}) \|
    & \leq &
    \frac{3\dltwu_{3}}{4} \| \DPN_{\nuiv} \, \IFN_{\nuiv}^{-1} \Av_{\nuiv} \|^{2} 
    & \leq &
    \frac{\dltwu_{3} \, \rrn_{\nuiv}}{2} \| \DPN_{\nuiv} \, \IFN_{\nuiv}^{-1} \Av_{\nuiv} \|
	\, .
	\qquad
\label{jhcvu7ejdytur39e9frtfw}
\end{EQ}
Also,
\begin{EQA}
	\bigl| 
		2 \fs_{\nuiv}(\targv_{\nuiv}) - 2 \fs_{\nuiv}(\targvs) + \| \IFN_{\nuiv}^{-1/2} \Av_{\nuiv} \|^{2} 
	\bigr|
	& \leq &
	\frac{5\dltwu_{3}}{2} \, \| \DPN_{\nuiv} \, \IFN_{\nuiv}^{-1} \Av_{\nuiv} \|^{3} \, .
	\qquad
\label{gtxddfujhyfdytet6ywerfd}
\end{EQA}
\end{proposition}

\begin{remark}
Conditions in \eqref{8difiyfc54wrboer7bjfrO} describe the range of applicability of the proposed approach.
Suppose for a moment that \( \DPN_{\nuiv}^{2} = \IFN_{\nuiv} \).
By \eqref{jhcvu7ejdytur39e9frtfw}, 
\( \| \DPN_{\nuiv} (\targv_{\nuiv} - \targvs) \| \approx \| \DPN_{\nuiv}^{-1} \Av_{\nuiv} \| \).
Therefore, the radius \( \rrn_{\nuiv} \) in the smoothness condition \nameref{LLsoT3ref} should be larger than
\( \| \DPN_{\nuiv}^{-1} \Av_{\nuiv} \| \), while the self-concordant constant \( \dltwu_{3} \) should 
be sufficiently small to ensure \( \dltwu_{3} \, \rrn_{\nuiv} \ll 1 \).
This is required in \eqref{8difiyfc54wrboer7bjfrO} in a precise manner.
In statistical applications, condition \( \dltwu_{3} \rrn_{\nuiv} \ll 1 \) is translated in the so-called
``critical dimension'' condition: effective number of parameters is smaller than effective sample size;
see \cite{Sp2024} for more discussion.
\end{remark}

\Subsection{Marginal optimization under (semi)orthogonality}
\label{SortLaplsemi}
Here we study variability of the value \( \targv_{\nuiv} = \arginf_{\targv} \fs(\targv,\nuiv) \)
w.r.t. the nuisance para\-meter \( \nuiv \).
It appears that local quadratic approximation of the function \( \fs \) in a vicinity of the extreme point
\( \prmtvs \) yields a nearly linear dependence of \( \targv_{\nuiv} \) on \( \nuiv \).
We illustrate this fact on the case of a quadratic function \( \fs(\cdot) \).
Consider the Hessian \( \IFT = \nabla^{2} \fs(\prmtvs) \) in the block form:
\begin{EQA}
	\IFT
	& \eqdef &
	\nabla^{2} \fs(\prmtvs)
	=
	\begin{pmatrix}
	\IFT_{\targv\targv} & \IFT_{\targv\nuiv}
	\\
	\IFT_{\nuiv\targv} & \IFT_{\nuiv\nuiv}
	\end{pmatrix} 
\label{hwe78yf2diwe76tfw67etfwtbso}
\end{EQA}
with \( \IFT_{\nuiv\targv} = \IFT_{\targv\nuiv}^{\T} \).
If \( \fs(\prmtv) \) is quadratic then \( \IFT \) and its blocks do not depend on \( \prmtv \).

\begin{lemma}
\label{Lpartmaxq}
Let \( \fs(\prmtv) \) be {quadratic}, {strongly convex}, and \( \nabla \fs(\prmtvs) = 0 \).
Then
\begin{EQA}
	\targv_{\nuiv} - \targvs 
	&=& 
	- \IFT_{\targv\targv}^{-1} \, \IFT_{\targv\nuiv} \, (\nuiv - \nuivs) .
\label{0fje7fhihy84efiewkw}
\end{EQA}
\end{lemma}

\begin{proof}
The condition \( \nabla \fs(\prmtvs) = 0 \) implies 
\( \fs(\prmtv) = \fs(\prmtvs) + (\prmtv - \prmtvs)^{\T} \IFT \, (\prmtv - \prmtvs)/2 \) 
with \( \IFT = \nabla^{2} \fs(\prmtvs) \).
For \( \nuiv \) fixed, the point \( \targv_{\nuiv} \) minimizes 
\( (\targv - \targvs)^{\T} \IFT_{\targv\targv} \, (\targv - \targvs) /2 + (\targv - \targvs)^{\T} \IFT_{\targv\nuiv} \, (\nuiv - \nuivs) \)
and thus, \( \targv_{\nuiv} - \targvs = - \IFT_{\targv\targv}^{-1} \, \IFT_{\targv\nuiv} \, (\nuiv - \nuivs) \).
\end{proof} 

Observation \eqref{0fje7fhihy84efiewkw} is a bit discouraging 
because the bias \( \targv_{\nuiv} - \targvs \) has the same magnitude as the nuisance parameter \( \nuiv - \nuivs \).
However, the condition \( \IFT_{\targv\nuiv} = 0 \) yields \( \targv_{\nuiv} \equiv \targvs \) 
and the bias vanishes.
If \( \fs(\prmtv) \) is not quadratic, the \emph{orthogonality} condition \( \nabla_{\nuiv} \nabla_{\targv} \, \fs(\targv,\nuiv) \equiv 0 \)
for all \( (\targv,\nuiv) \in \WV \) still ensures a vanishing bias.

\begin{lemma}
\label{Lpartmax}
Let \( \fs(\targv,\nuiv) \) be continuously differentiable and 
\( \nabla_{\nuiv} \nabla_{\targv} \, \fs(\targv,\nuiv) \equiv 0 \).
Then the point 
\( \targv_{\nuiv} = \arginf_{\targv} f(\targv,\nuiv) \) does not depend on \( \nuiv \).
\end{lemma}

\begin{proof}
The condition \( \nabla_{\nuiv} \nabla_{\targv} \, \fs(\targv,\nuiv) \equiv 0 \) implies the decomposition 
\( \fs(\targv,\nuiv) = \fs_{1}(\targv) + \fs_{2}(\nuiv) \) for some functions 
\( \fs_{1} \) and \( \fs_{2} \).
This in turn yields \( \targv_{\nuiv} \equiv \targvs \).
\end{proof}

In some cases, one can check \emph{semi-orthogonality} condition 
\begin{EQA}
	\nabla_{\nuiv} \nabla_{\targv} \, \fs(\targvs,\nuiv)  
	&=&
	0,
	\qquad
	\forall \nuiv \in \Nui \, .
\label{jdgyefe74erjscgygfydhwse}
\end{EQA}
\ifapp{A typical example is given by nonlinear regression; see Section~\ref{SnonlinLA}.}{}

\begin{lemma}
\label{Lsemiorto}
Assume \eqref{jdgyefe74erjscgygfydhwse}.
Then 
\begin{EQA}[rcccl]
	\nabla_{\targv} \, \fs(\targvs,\nuiv)
	& \equiv &
	0,
	\qquad
	\nabla_{\targv\targv}^{2} \fs(\targvs,\nuiv) 
	& \equiv &
	\nabla_{\targv\targv}^{2} \fs(\targvs,\nuivs),
	\qquad
	\nuiv \in \Nui .
\label{0xyc5ftr4j43vefvuruendt}
\end{EQA}
Moreover, if \( \fs(\targv,\nuiv) \) is convex in \( \targv \) given \( \nuiv \) then
\begin{EQA}
	\targv_{\nuiv} 
	& \eqdef & 
	\arginf_{\targv} \fs(\targv,\nuiv) \equiv \targvs ,
	\qquad
	\forall \nuiv \in \Nui \,.
\label{irdtyrnjjutu45r7gjhr}
\end{EQA}
\end{lemma}

\begin{proof}
Consider \( \Av_{\nuiv} = \nabla_{\targv} \, \fs(\targvs,\nuiv) \).
Obviously \( \Av_{\nuivs} = 0 \).
Moreover, \eqref{jdgyefe74erjscgygfydhwse} implies that \( \Av_{\nuiv} \) 
does not depend on \( \nuiv \) and thus, vanishes everywhere.
As \( \fs \) is convex in \( \targv \), 
this implies \( \fs(\targvs,\nuiv) = \inf_{\targv} \fs(\targv,\nuiv) \)
and \( \targv_{\nuiv} = \arginf_{\targv} \fs(\targv,\nuiv) \equiv \targvs \).
Similarly, by \eqref{jdgyefe74erjscgygfydhwse}, it holds
\( \nabla_{\nuiv} \nabla_{\targv\targv} \, \fs(\targvs,\nuiv) \equiv 0 \) 
and \eqref{0xyc5ftr4j43vefvuruendt} follows.
Convexity of \( \fs(\targv,\nuiv) \) in \( \targv \) for \( \nuiv \) fixed and \eqref{0xyc5ftr4j43vefvuruendt} 
imply \eqref{irdtyrnjjutu45r7gjhr}.
\end{proof}

Orthogonality or semi-orthogonality \eqref{jdgyefe74erjscgygfydhwse} conditions are rather 
restrictive and fulfilled only in special situations.
Later we explore a general situation and present some condition allowing to control the semiparametric bias.

\Subsection{Semiparametric bias. A linear approximation}
\label{Sonepoint}


This section presents a first-order expansion for the semiparametric bias \( \targv_{\nuiv} - \targvs \) under  
conditions on the cross-derivatives of \( \fs(\targv,\nuiv) \) for \( \nuiv = \nuivs \)
or \( \targv = \targvs \).
Suppose that the nuisance variable \( \nuiv \) is already localized to a small vicinity \( \Nui \) of \( \nuivs \).
A typical example is given by the level set \( \Nui \) of the form
\begin{EQA}
	\Nui
	&=&
	\bigl\{ \nuiv \colon \| \HPN (\nuiv - \nuivs) \|_{\nano} \leq \rrnano \bigr\} ,
\label{6chbdtw2hydfuguye4223w23}
\end{EQA}
where \( \HPN \) is a metric tensor in \( \R^{\dimq} \), 
\( \| \cdot \|_{\nano} \) is a norm in \( \R^{\dimq} \), and \( \rrnano > 0 \).
Often \( \| \cdot \|_{\nano} \) is the usual \( \ell_{2} \)-norm.
However, in some situations, it is beneficial to use different topologies for the target parameter \( \targv \)
and the nuisance parameter \( \nuiv \).
One example is given by estimation in \( \sup \)-norm with \( \| \cdot \|_{\nano} = \| \cdot \|_{\infty} \).

Later, we assume that for any \( \nuiv \in \Nui \), 
the function \( \fs_{\nuiv} \) satisfies \nameref{LLsoT3ref} with
\( \DPN_{\nuiv}^{2} \equiv \DPN^{2} \), \( \rrn_{\nuiv} \equiv \rrn \), and a small constant \( \dltwu_{3} \).
In addition, introduce the following condition.

\begin{description}


    \item[\label{LLpT3ref} \( \bb{(\mathcal{T}_{3,\Nui}^{*})} \)]
      \emph{
      It holds with some \( \dltwun, \dltwunn \)
}
\begin{EQA}[rccl]
	\sup_{\nuiv \in \Nui} \,\, 
	\sup_{\zv \in \R^{\dimp}, \wv \in \R^{\dimq}} \,
	&
	\frac{| \langle \nabla_{\targv\nuiv\nuiv}^{3} \fs(\targvs,\nuiv), \zv \otimes \wv^{\otimes 2} \rangle |}
		 {\| \DPN \zv \| \, \| \HPN \wv \|_{\nano}^{2}}
	& \leq &
	\dltwun \, ,
	\qquad
\label{c6ceyecc5e5etcT2}
	\\
	\sup_{\nuiv \in \Nui} \,\, 
	\sup_{\zv \in \R^{\dimp}, \wv \in \R^{\dimq}} \,
	&
	\frac{| \langle \nabla_{\targv\targv\nuiv}^{3} \fs(\targvs,\nuiv), \zv^{\otimes 2} \otimes \wv \rangle |}
		 {\| \DPN \zv \|^{2} \, \| \HPN \wv \|_{\nano}}
	& \leq &
	\dltwunn \, .
\label{c6ceyecc5e5etctwhcyegwc}
\end{EQA}
\end{description}


\begin{remark}
Condition \nameref{LLpT3ref} only involves mixed derivatives 
\( \nabla_{\targv\nuiv\nuiv}^{3} \fs(\targvs,\nuiv) \) and \( \nabla_{\targv\targv\nuiv}^{3} \fs(\targvs,\nuiv) \) 
of \( \fs(\targvs,\nuiv) \) for the fixed value \( \targv = \targvs \),
while condition \nameref{LLsoT3ref} only concerns smoothness of \( \fs(\targv,\nuiv) \) 
in \( \targv \) for \( \nuiv \) fixed.
Therefore, the combination of \nameref{LLsoT3ref} and \nameref{LLpT3ref} is much weaker than the full dimensional 
condition \nameref{LLsT3ref}.
\end{remark}

%
%
%
%
Let \( \IFT_{\targv\targv} \, , \IFT_{\targv\nuiv} \, , \IFT_{\nuiv\targv} \, , \IFT_{\nuiv\nuiv} \) be the blocks of 
\( \IFT = \nabla^{2} \fs(\upsvs) \); see \eqref{hwe78yf2diwe76tfw67etfwtbso}.
Introduce the operator \( \IFC \eqdef \DPN^{-1} \IFT_{\targv\nuiv} \, \HPN^{-1} \colon \R^{\dimq} \to \R^{\dimp} \) and
its dual norm
\begin{EQA}
\label{f8vuehery6gv65ftehwee}
	\crossnorm
	& \eqdef &
	\| \IFC \|_{\dual} 
	=
	\sup_{\zv \colon \| \zv \|_{\nano} \leq 1} \| \IFC \zv \| 
	\, .
\end{EQA}
Also, define technical values \( \dltwbun \eqdef \dltwunn \, \rrnano < 1 \),
\begin{EQA}[rcccl]
	\crossgrad 
	& \eqdef &
	\frac{3}{2(1 - \dltwbun)}\Bigl( \crossnorm + \frac{\dltwun \, \rrnano}{2} \Bigr) 
	\, ,
	\qquad
	\dltwuns
	& \eqdef &
	\frac{1}{1 - \dltwbun} \, 
	\biggl( \crossnorm \, \dltwunn + \frac{\dltwun}{2}
		+ \frac{\crossgrad^{2} \, \dltwu_{3}}{3} 
	\biggr) \, .
	\qquad
\label{36gfijh94ejdvtwekoised}
\end{EQA}
In a typical situation with \( \crossnorm \leq 1 \),
\( \dltwun \leq \dltwu_{3} \), 
\( \dltwunn \leq \dltwu_{3} \), and \( \dltwu_{3} \, \rrnano \leq 1/4 \), it holds
\( \dltwbun \leq 1/4 \),
\( \crossgrad \leq 2 \crossnorm + \dltwu_{3} \, \rrnano \), and
\( \dltwuns \leq 4 \dltwu_{3} \).


%
The next result provides an expansion of the semiparametric bias \( \targv_{\nuiv} - \targvs \) under \nameref{LLpT3ref} and \nameref{LLsoT3ref}
with a leading linear term as in quadratic case of Lemma~\ref{Lpartmaxq}. 
%

\begin{theorem}
\label{PsemiAvex}
Let \( \Nui \) be given by \eqref{6chbdtw2hydfuguye4223w23}.
Assume the following conditions:
\bulletc
The solution 
\( \upsvs = (\targvs,\nuivs) = \arginf_{\upsv = (\targv,\nuiv) \colon \nuiv \in \Nui} \fs(\upsv) \) 
exists and is unique. 

\bulletc
For any \( \nuiv \in \Nui \), the function \( \fs_{\nuiv}(\targv) \) is strongly convex and  satisfies \nameref{LLsoT3ref}  
with \( \rrn_{\nuiv} \equiv \rrn \) and \( \DPN_{\nuiv}^{2} \equiv \DPN^{2} \leq \IFN \),
where \( \IFN = \IFT_{\targv\targv} \).

\bulletc
The function \( \fs(\targv,\nuiv) \) satisfies \nameref{LLpT3ref}.

\bulletc
With \( \crossgrad \) from \eqref{36gfijh94ejdvtwekoised}, it holds
\begin{EQA}[c]
	\dltwunn \, \rrnano \leq \frac{1}{4} \, ,
	\quad
	\rr 
	\geq 
	\crossgrad \, 
	\rrnano
	\, ,
	\quad
	\crossgrad \, 
	\dltwu_{3} \, \rrnano
	\leq 
	\frac{2}{3} 
	\, .
	\qquad
\label{6cuydjd5eg3jfggu8eyt4h}
\end{EQA}
Then the partial solution \( \targv_{\nuiv} \) satisfies 
\begin{EQA}
	\| \DPN (\targv_{\nuiv} - \targvs) \|
	& \leq &
	\crossgrad \, \| \HPN (\nuiv - \nuivs) \|_{\nano} \,. 
\label{rhDGtuGmusGU0a2spt}
\end{EQA}
Moreover, for any linear mapping \( \QP \) on \( \R^{\dimp} \), it holds with \( \dltwuns \) from \eqref{36gfijh94ejdvtwekoised}
\begin{EQA}[c]
\label{36gfijh94ejdvtwekoise}
	\| \QP \{ \targv_{\nuiv} - \targvs + \IFN^{-1} \IFT_{\targv\nuiv} (\nuiv - \nuivs) \} \|
	\leq 
	\| \QP \, \IFN^{-1} \DPN \| \, \dltwuns
	\, \| \HPN (\nuiv - \nuivs) \|_{\nano}^{2}
	\; .
	\qquad
\end{EQA}
\end{theorem}

\begin{remark}
\label{Rdmax2}
An extension to the case \( \DPN^{2} \leq \dmax^{2} \IFN \) can be done by replacing everywhere
\( \dltwu_{3} \),
\( \dltwunn \),
\( \dltwun \),
\( \dltwuns \),
\( \crossnorm \)
with \( \dmax^{3} \dltwu_{3} \), \( \dmax^{2} \dltwunn \), 
\( \dmax \, \dltwun \), 
\( \dmax \, \dltwuns \), 
\( \dmax^{-1} \crossnorm \)
respectively.
\end{remark}

\Subsection{A linear perturbation}
Let \( \fd(\targv,\nuiv) \) be obtained by a linear perturbation of \( \fs(\targv,\nuiv) \):
\begin{EQA}
	\fd(\targv,\nuiv) - \fd(\targvs,\nuivs)
	&=&
	\fs(\targv,\nuiv) - \fs(\targvs,\nuivs) + \langle \AAv,\targv - \targvs \rangle + \langle \CCv,\nuiv - \nuivs \rangle 
	\, ;
	\qquad
\label{s6cdt76ws6656ysjwsdftc6}
\end{EQA}
cf. \eqref{4hbh8njoelvt6jwgf09}.
Let also \( \Nui \) be given by \eqref{6chbdtw2hydfuguye4223w23} and \( \nuiv \in \Nui \).
We are interested in quantifying the distance between \( \targvs \) and \( \targvn_{\nuiv} \), where 
\begin{EQA}
	\targvn_{\nuiv}
	& \eqdef &
	\arginf_{\targv} \fd(\targv,\nuiv) .
\label{8dcudxjd5t6t5mkwsyddf}
\end{EQA}
The linear perturbation \( \langle \CCv,\nuiv - \nuivs \rangle \) does not depend on \( \targv \) and it can be ignored.

\begin{theorem}
\label{Pconcsupp}
Assume the conditions of Theorem~\ref{PsemiAvex}. Let also \( \AAv \) in \eqref{s6cdt76ws6656ysjwsdftc6} satisfy
\begin{EQA}[c]
	\rrn \geq 2 \| \DPN \, \IFN^{-1} \AAv \| \, ,
	\quad
	\dltwu_{3} \| \DPN \, \IFN^{-1} \AAv \| < \frac{1}{3} \, .
\end{EQA}
Then with \( \dltwbun = \dltwunn \, \rrnano \leq 1/4 \) and \( \crossgrad \) from \eqref{36gfijh94ejdvtwekoised},
it holds for any \( \nuiv \in \Nui \)
\begin{EQA}
	\| \DPN (\targvn_{\nuiv} - \targvs) \|
	& \leq &
	\crossgrad \, \| \HPN (\nuiv - \nuivs) \|_{\nano} 
	+ \frac{3 }{2(1 - \dltwbun)} \, \| \DPN \, \IFN^{-1} \AAv \| 
	\, .
\label{ysd7euwejfg653wthf873k}
\end{EQA}
Moreover, for any linear mapping \( \QP \), it holds with 
\( \dltwuns \) from \eqref{36gfijh94ejdvtwekoised} 
\begin{EQA}[l]
	\!\!\!
	\| \QP \{ \targvn_{\nuiv} - \targvs + \IFN^{-1} \IFT_{\targv\nuiv} (\nuiv - \nuivs) + \IFN^{-1} \AAv \} \|
	\\
	\leq 
	\| \QP \, \IFN^{-1} \DPN \| \, \Bigl\{ (\dltwuns + \dltwunn) \, \| \HPN (\nuiv - \nuivs) \|_{\nano}^{2}
		+ (2 \dltwu_{3} + \dltwunn/2) \, \| \DPN \, \IFN^{-1} \AAv \|^{2}
	\Bigr\}	
	\, .
	\qquad 
\label{usdhyw6hikhurnetr256}
\end{EQA}
\end{theorem}

\begin{remark}
If \( \dltwun \leq \dltwu_{3} \), 
\( \dltwunn \leq \dltwu_{3} \), and \( \dltwu_{3} \, \rrnano \leq 1/4 \), then
\begin{EQA}[c]
	(\dltwuns + \dltwunn) \, \| \HPN (\nuiv - \nuivs) \|_{\nano}^{2}
		+ \bigl( 2 \dltwu_{3} + \frac{\dltwunn}{2} \bigr) \, \| \DPN \, \IFN^{-1} \AAv \|^{2}
	\leq 
	\frac{5 \dltwu_{3}}{2} \, \bigl\{ 2 \| \HPN (\nuiv - \nuivs) \|_{\nano}^{2}
		+ \| \DPN \, \IFN^{-1} \AAv \|^{2} \bigr\}
	\, .
\end{EQA}
The results \eqref{ysd7euwejfg653wthf873k} and \eqref{usdhyw6hikhurnetr256} describe the difference
\( \targvn_{\nuiv} - \targvs \) caused by two sources:
a linear perturbation \( \langle \AAv,\targv - \targvs \rangle \) of the objective function
and the use of a wrong nuisance value \( \nuiv \) in place of \( \nuivs \).
The required conditions are \( \dltwu_{3} \| \DPN \, \IFN^{-1} \AAv \| \ll 1 \) and 
\( \dltwunn \, \rrnano \ll 1 \).
\end{remark}

%% file: marg_opt_proofs.tex

%
\def\fn{\fs^{\circ}}
\def\IFN{\IF}


\paragraph{Proof of Proposition~\ref{PbiassemiN}}
The function \( \fn_{\nuiv}(\targv) = \fs_{\nuiv}(\targv) - \langle \targv,\Av_{\nuiv} \rangle \) is convex 
and \( \nabla \fn_{\nuiv}(\targvs) = 0 \) yielding 
\( \targvs = \argmin_{\targv} \fn_{\nuiv}(\targv) \).
Now Theorem~\ref{PFiWigeneric2} yields \eqref{jhcvu7ejdytur39e9frtfw} and by \eqref{3d3Af12DGttGa2}
\begin{EQA}
	\Bigl| 2 \fn_{\nuiv}(\targvs) - 2 \fn_{\nuiv}(\targv_{\nuiv}) + \| \IFN_{\nuiv}^{-1/2} \Av_{\nuiv} \|^{2} \Bigr|
	& \leq &
	\frac{\dltwu_{3}}{2} \, \| \DPN_{\nuiv} \, \IFN_{\nuiv}^{-1} \Av_{\nuiv} \|^{3} 
	\, .
\label{ytduw237832jfg656335hdi}
\end{EQA}
Equivalently,
\begin{EQA}
	\bigl| 
		2 \fs_{\nuiv}(\targvs) - 2 \fs_{\nuiv}(\targv_{\nuiv}) - 2 \langle \targvs - \targv_{\nuiv} , \Av_{\nuiv} \rangle 
		+ \| \IFN_{\nuiv}^{-1/2} \Av_{\nuiv} \|^{2} 
	\bigr|
	& \leq &
	\frac{\dltwu_{3}}{2} \, \| \DPN_{\nuiv} \, \IFN_{\nuiv}^{-1} \Av_{\nuiv} \|^{3} \, .
	\qquad
\label{gtxdbcbcgdfswerfdsdew}
\end{EQA}
By \eqref{jhcvu7ejdytur39e9frtfw}
\begin{EQA}
	&& \nquad
	\bigl| \langle \targvs - \targv_{\nuiv} , \Av_{\nuiv} \rangle + \| \IFN_{\nuiv}^{-1/2} \Av_{\nuiv} \|^{2} \bigr|
	=
	\bigl| \langle \IFN_{\nuiv}^{1/2} (\targvs - \targv_{\nuiv}) , \IFN_{\nuiv}^{-1/2} \Av_{\nuiv} \rangle 
	+ \langle \IFN_{\nuiv}^{-1/2} \Av_{\nuiv},\IFN_{\nuiv}^{-1/2} \Av_{\nuiv} \rangle \bigr|
	\\
	&=&
	\bigl| \langle \IFN_{\nuiv}^{1/2} (\targvs - \targv_{\nuiv} + \IFN_{\nuiv}^{-1} \Av_{\nuiv}) , 
		\IFN_{\nuiv}^{-1/2} \Av_{\nuiv} \rangle 
	\bigr|
	\\
	& \leq &
	\| \DPN_{\nuiv}^{-1} \, \IFN_{\nuiv} (\targvs - \targv_{\nuiv} + \IFN_{\nuiv}^{-1} \Av_{\nuiv}) \| \, \,
		\| \DPN_{\nuiv} \, \IFN_{\nuiv}^{-1} \Av_{\nuiv} \|
	\leq 
	\frac{3\dltwu_{3}}{4} \, \| \DPN_{\nuiv} \, \IFN_{\nuiv}^{-1} \Av_{\nuiv} \|^{3}
\label{duy7cyf6re6eytedhge2}
\end{EQA}
This and \eqref{gtxdbcbcgdfswerfdsdew} imply \eqref{gtxddfujhyfdytet6ywerfd}.

\paragraph{Proof of Theorem~\ref{PsemiAvex}}
First, we bound the variability of 
\( \IFN_{\nuiv} = \nabla_{\targv\targv} \fs(\targvs,\nuiv) \) over \( \Nui \). 

\begin{lemma}
\label{L3IFNNui}
Assume \nameref{LLpT3ref}.
Then for any \( \nuiv \in \Nui \), with \( \IFN = \IFN_{\nuivs} \),
\begin{EQA}
	\| \DPN^{-1} \, (\IFN_{\nuiv} - \IFN) \, \DPN^{-1} \|
	& \leq &
	\dltwunn \, \| \HPN (\nuiv - \nuivs) \|_{\nano} 
	\, .
\label{8c7c67cc6c63kdldlvvudw}
\end{EQA}
\end{lemma}

\begin{proof}
Let \( \nuiv \in \Nui \).
By \eqref{c6ceyecc5e5etctwhcyegwc}, for any \( \zv \in \R^{\dimp} \), it holds
\begin{EQA}
	&& \!\!\!\!
	\bigl| \bigl\langle \DPN^{-1} \, (\IFN_{\nuiv} - \IFN) \, \DPN^{-1}, \zv^{\otimes 2} \bigr\rangle \bigr|
	=
	\bigl| \bigl\langle \IFN_{\nuiv} - \IFN, (\DPN^{-1} \zv)^{\otimes 2} \bigr\rangle \bigr|
	\\
	& \leq &
	\sup_{t \in [0,1]}
	\bigl| \langle \nabla_{\targv\targv\nuiv}^{3} \fs(\targvs,\nuivs + t (\nuiv - \nuivs)), (\DPN^{-1} \zv)^{\otimes 2} \otimes (\nuiv - \nuivs) \rangle \bigr|
	\leq 
	\dltwunn \, \| \zv \|^{2}  \, \| \HPN (\nuiv - \nuivs) \|_{\nano} \, .
\label{9487654r5tghasdfgsup}
\end{EQA}
This yields \eqref{8c7c67cc6c63kdldlvvudw}.
\end{proof}

The next result describes some corollaries of \eqref{8c7c67cc6c63kdldlvvudw}. 

\begin{lemma}
\label{LIFvari}
Let \( \DPN^{2} \leq \dmax^{2} \IFN \), and for some other positive symmetric \( \dimp \)-matrix \( \IFN_{1} \) 
\begin{EQA}
	\| \DPN^{-1} \, (\IFN_{1} - \IFN) \, \DPN^{-1} \|
	& \leq &
	\dmax^{-2} \dltwbun 
\label{icfuifcu7dy7e3gw3ft6e}
\end{EQA}
with \( \dltwbun < 1 \).
Then
\begin{EQA}[ccc]
\label{8c7c67cc6c63kdldlvvudw2}
	\| \IFN^{-1/2} \, (\IFN_{1} - \IFN) \, \IFN^{-1/2} \|
	& \leq &
	\dltwbun 
	\, ,
	\\
	\| \IFN^{1/2} \, (\IFN_{1}^{-1} - \IFN^{-1}) \, \IFN^{1/2} \|
	& \leq &
	\frac{\dltwbun}{1 - \dltwbun} 
	\, ,
\label{g98g75re3gf76r4egu4e}
\end{EQA}
and
\begin{EQA}
	\frac{1}{1 + \dltwbun} \, \| \DPN \, \IFN^{-1} \, \DPN \|
	\leq 
	\| \DPN \, \IFN_{1}^{-1} \DPN \|
	& \leq &
	\frac{1}{1 - \dltwbun} \, \| \DPN \, \IFN^{-1} \, \DPN \|
	\, .
\label{g98g75re3gf76r4egu4ee}
\end{EQA}
Furthermore, for any vector \( \uv \)
\begin{EQA}[rcccl]
\label{uciuvcdiu3eoooodffcy}
	(1 - \dltwbun) \| \DPN^{-1} \IFN \uv \|
	& \leq &
	\| \DPN^{-1} \IFN_{1} \uv \|
	& \leq &
	(1 + \dltwbun) \| \DPN^{-1} \IFN \uv \| \, ,
	\qquad
	\\
	\frac{1 - 2 \dltwbun}{1 - \dltwbun} \| \DPN \, \IFN^{-1} \uv \|
	& \leq &
	\| \DPN \, \IFN_{1}^{-1} \uv \|
	& \leq &
	\frac{1}{1 - \dltwbun} \| \DPN \, \IFN^{-1} \uv \| \, .
	\qquad
\label{uciuvcdiu3eoooodffcym}
\end{EQA}
\end{lemma}

\begin{proof}
Rescaling of \( \DPN \) by \( \dmax^{-1} \) reduces the proof to \( \dmax = 1 \).
Now \eqref{icfuifcu7dy7e3gw3ft6e} implies \eqref{8c7c67cc6c63kdldlvvudw2} 
because of \( \IFN^{-1} \leq \DPN^{-2} \).
Define now \( \Uv \eqdef \IFN^{-1/2} \, (\IFN_{1} - \IFN) \, \IFN^{-1/2} \). 
Then \( \| \Uv \| \leq \dltwbun \) and
\begin{EQA}
	\| \IFN^{1/2} \, (\IFN_{1}^{-1} - \IFN^{-1}) \, \IFN^{1/2} \|
	&=&
	\| (\Id + \Uv)^{-1} - \Id \|
	\leq 
	\frac{1}{1 - \dltwbun} \| \Uv \|
\label{uefu83yr83y38frkfbg7e}
\end{EQA}
yielding \eqref{g98g75re3gf76r4egu4e}.
Further, \( \IFN_{1} (\IFN_{1}^{-1} - \IFN^{-1}) \, \IFN = \IFN - \IFN_{1} \), \( \IFN^{-1} \leq \DPN^{-2} \), and \eqref{icfuifcu7dy7e3gw3ft6e} imply
\begin{EQA}
	\| \DPN \, (\IFN_{1}^{-1} - \IFN^{-1}) \, \DPN \|
	&=&
	\| \DPN \, \IFN_{1}^{-1} \IFN_{1} (\IFN_{1}^{-1} - \IFN^{-1}) \, \IFN \, \IFN^{-1} \DPN \|
	\\
	&=&
	\| \DPN \, \IFN_{1}^{-1} \; \DPN \, \DPN^{-1} (\IFN_{1} - \IFN) \, \DPN^{-1} \DPN \; \IFN^{-1} \DPN \|
	\\
	& \leq &
	\| \DPN \, \IFN_{1}^{-1} \DPN \| \, \| \DPN \, \IFN^{-1} \DPN \| \, \| \DPN^{-1} (\IFN_{1} - \IFN) \DPN^{-1} \| 
	\leq 
	\dltwbun \| \DPN \, \IFN_{1}^{-1} \DPN \| 
	\, ,
\label{odfu7eujgvy55r44rf}
\end{EQA}
and \eqref{g98g75re3gf76r4egu4ee} follows.
Also, by \( \DPN^{2} \leq \IFN \)
\begin{EQA}
	\| \DPN^{-1} \IFN_{1} \uv \|
	& \leq &
	\| \DPN^{-1} \IFN \uv \| + \| \DPN^{-1} (\IFN_{1} - \IFN) \DPN^{-1} \DPN \uv \|
	\leq 
	\| \DPN^{-1} \IFN \uv \| + \dltwbun \, | \DPN \uv \|
	\\
	& \leq & 
	\| \DPN^{-1} \IFN \uv \| + \dltwbun \| \DPN^{-1} \IFN \uv \| 
	\leq 
	(1 + \dltwbun) \| \DPN^{-1} \IFN \uv \| \, ,
	\qquad
\label{uciuvcdiu3ery7httrwbgfcy}
\end{EQA}
and \eqref{uciuvcdiu3eoooodffcy} follows.
Similarly, by \( \IFN_{1}^{-1} (\IFN_{1} - \IFN) \IFN^{-1} = \IFN^{-1} - \IFN_{1}^{-1} \), it holds
\begin{EQA}
	\| \DPN \, (\IFN_{1}^{-1} - \IFN^{-1}) \uv \| 
	& = &
	\| \DPN \, \IFN_{1}^{-1} (\IFN_{1} - \IFN) \IFN^{-1} \uv \|
	=
	\| \DPN \, \IFN_{1}^{-1} \; \DPN \DPN^{-1} (\IFN_{1} - \IFN) \, \DPN^{-1} \DPN \; \IFN^{-1} \uv \|
	\\
	& \leq &
	\| \DPN^{-1} \, (\IFN_{1} - \IFN) \, \DPN^{-1} \| \,\,
	\| \DPN \, \IFN_{1}^{-1} \, \DPN \| \,\, \| \DPN \, \IFN^{-1} \uv \|
	\leq 
	\frac{\dltwbun}{1 - \dltwbun} \| \DPN \, \IFN^{-1} \uv \|
	\, ,
\label{rtiuru8bttg3hfy7vrgf}
\end{EQA}
and \eqref{uciuvcdiu3eoooodffcym} follows as well.
\end{proof}

Definition \eqref{f8vuehery6gv65ftehwee} implies the following bound.
\begin{lemma}
\label{Lsemidual}
For any \( \nuiv \in \Nui \) and any linear mapping \( \QP \), it holds
\begin{EQA}
	\| \QP \IFT_{\targv\nuiv} (\nuiv - \nuivs) \|
	& \leq &
	\| \QP \IFT_{\targv\nuiv} \, \HPN^{-1} \|_{\dual} \, \| \HPN (\nuiv - \nuivs) \|_{\nano} \, .
\label{difkjv78frhjrjhfh6wekd}
\end{EQA}
\end{lemma}

The next lemma shows that \( \Av_{\nuiv} = \nabla \fs_{\nuiv}(\targvs) \) is nearly linear in \( \nuiv \).

\begin{lemma}
\label{LvarDVa}
Assume \eqref{c6ceyecc5e5etcT2}.
Then with \( \crossnorm = \| \DPN^{-1} \IFT_{\targv\nuiv} \, \HPN^{-1} \|_{\dual} \),
it holds for any \( \nuiv \) with \( \| \HPN (\nuiv - \nuivs) \|_{\nano} \leq \rru \)
\begin{EQA}[rcl]
	\| \DPN^{-1} \{ \Av_{\nuiv} - \IFT_{\targv\nuiv} (\nuiv - \nuivs) \} \|
	& \leq &
	\frac{\dltwun}{2} \, \| \HPN (\nuiv - \nuivs) \|_{\nano}^{2} 
	\leq 
	\frac{\dltwun \, \rru}{2} \, \| \HPN (\nuiv - \nuivs) \|_{\nano} 
	\, ,
\label{4cbijjm09j9hhrjfifkrs3}
\end{EQA}
and
\begin{EQA}[rcl]
	\| \DPN^{-1} \Av_{\nuiv} \|
	& \leq &
	\crossnorm \, \| \HPN (\nuiv - \nuivs) \|_{\nano}
	+ \frac{\dltwun}{2} \, \| \HPN (\nuiv - \nuivs) \|_{\nano}^{2}
	\leq 
	\Bigl( \crossnorm + \frac{\dltwun \, \rru}{2} \Bigr) \, \| \HPN (\nuiv - \nuivs) \|_{\nano} 
	\, .
	\qquad \quad
\label{4cbijjm09j9hhrjfifkrs}
\end{EQA}
\end{lemma}

\begin{proof}
Fix \( \nuiv \in \Nui \) and define \( \av_{\nuiv}(t) \eqdef \Av_{\nuivs + t (\nuiv - \nuivs)} \).
Then \( \av_{\nuiv}(0) = \Av_{\nuivs} = 0 \), \( \av_{\nuiv}(1) = \Av_{\nuiv} \), 
and
\begin{EQA}
	\av_{\nuiv}(1) - \av_{\nuiv}(0) 
	&=&
	\int_{0}^{1} \av'_{\nuiv}(t) \, dt \, ,
\label{dsf7hgduw3heiuwed53ww3}
\end{EQA}
where \( \av'_{\nuiv}(t) = \frac{d}{dt} \av_{\nuiv}(t) \) for \( t \in [0,1] \).
Similarly, by \( \av'_{\nuiv}(0) = \IFT_{\targv\nuiv} (\nuiv - \nuivs) \), we derive 
\begin{EQA}
	&& \nquad
	\av_{\nuiv}(1) - \av_{\nuiv}(0) - \av'_{\nuiv}(0)
	\\
	&=&
	\int_{0}^{1} \{ \av'_{\nuiv}(t) - \av'_{\nuiv}(0) \} \, dt
	=
	\int_{0}^{1} \int_{0}^{t} \av''_{\nuiv}(u) \, du \, dt
	=
	\int_{0}^{1} (1 - u) \, \av''_{\nuiv}(u) \, du \, ,
\label{g8eh35fgtg7t76868jdy}
\end{EQA}
where 
\( \av''_{\nuiv}(t) = \frac{d^{2}}{dt^{2}} \av_{\nuiv}(t) \).
By condition \nameref{LLpT3ref} 
\begin{EQA}
	\bigl| \langle \av''_{\nuiv}(t), \zv \rangle \bigr|
	&=&
	\bigl| \bigl\langle \nabla_{\targv\nuiv\nuiv}^{3} \fs(\targvs,\nuivs + t (\nuiv - \nuivs)), \zv \otimes (\nuiv - \nuivs)^{\otimes 2} \bigr\rangle \bigr|
	\leq 
	\dltwun \, \| \DPN \zv \| \, \| \HPN (\nuiv - \nuivs) \|_{\nano}^{2} \, 
\label{5cbgcfyt6webhwefuy6}
\end{EQA}
and hence,
\begin{EQA}
	\| \DPN^{-1} \av''_{\nuiv}(t) \|
	&=&
	\sup_{\zv \colon \| \zv \| \leq 1} \bigl| \langle \DPN^{-1} \av''_{\nuiv}(t), \zv \rangle \bigr|
	=
	\sup_{\zv \colon \| \zv \| \leq 1}  
		\bigl| \langle \av''_{\nuiv}(t), \DPN^{-1} \zv \rangle \bigr|
	\\
	& \leq &
	\dltwun \sup_{\zv \colon \| \zv \| \leq 1}
	\| \zv \| \, \| \HPN (\nuiv - \nuivs) \|_{\nano}^{2} \, 
	=
	\dltwun  \| \HPN (\nuiv - \nuivs) \|_{\nano}^{2} \, .
\label{hgvytew3hw36fdfnhrg}
\end{EQA}
This yields
\begin{EQA}
	\| \DPN^{-1} \{ \Av_{\nuiv} - \IFT_{\targv\nuiv} (\nuiv - \nuivs) \} \|
	& \leq &
	\dltwun \, \| \HPN (\nuiv - \nuivs) \|_{\nano}^{2} \, \int_{0}^{1} (1 - t) \, dt
	\leq 
	\frac{\dltwun \, \| \HPN (\nuiv - \nuivs) \|_{\nano}^{2}}{2} \, 
\label{4cbijjm09j9hhrjfifkr}
\end{EQA}
as claimed in \eqref{4cbijjm09j9hhrjfifkrs3}. 
Lemma~\ref{Lsemidual} implies \eqref{4cbijjm09j9hhrjfifkrs}.
\end{proof}

\noindent
Now Proposition~\ref{PbiassemiN} helps to show \eqref{rhDGtuGmusGU0a2spt} and to bound
\( \targv_{\nuiv} - \targvs + \IFN_{\nuiv}^{-1} \Av_{\nuiv} \).

\begin{lemma}
Let \( \| \HPN (\nuiv - \nuivs) \|_{\nano} \leq \rru \) and
\( \dltwunn \, \rru \leq \dltwbun < 1 \).
Then with \( \crossgrad \) from \eqref{36gfijh94ejdvtwekoised}
\begin{EQA}
	\| \DPN (\targv_{\nuiv} - \targvs) \|
	& \leq &
	\frac{3 }{2} \| \DPN \, \IFN_{\nuiv}^{-1} \Av_{\nuiv} \| 
	\leq 
	\crossgrad \, \| \HPN (\nuiv - \nuivs) \|_{\nano}
	\, . 
\label{rhDGtuGmusGU0a2sp}
\end{EQA}
Moreover, 
\begin{EQA}
	\| \DPN^{-1} \IFN_{\nuiv} (\targv_{\nuiv} - \targvs + \IFN_{\nuiv}^{-1} \Av_{\nuiv}) \|
	& \leq &
	\frac{\dltwu_{3} \, \crossgrad^{2}}{3} \, 
	\| \HPN (\nuiv - \nuivs) \|_{\nano}^{2}
	\, .
	\qquad
\label{d73jejg7tre3fdy3fujell}
\end{EQA}
\end{lemma}

\begin{proof}
By \( \DPN^{2} \leq \IFN \), 
\eqref{uciuvcdiu3eoooodffcym} of Lemma~\ref{LIFvari}, and \eqref{4cbijjm09j9hhrjfifkrs} 
\begin{EQA}
	\frac{3 }{2} \| \DPN \, \IFN_{\nuiv}^{-1} \Av_{\nuiv} \|
	& \leq &
	\frac{3 }{2} \| \DPN \, \IFN_{\nuiv}^{-1} \DPN \| \,\, \| \DPN^{-1} \Av_{\nuiv} \|
	\leq 
	\frac{3}{2(1 - \dltwbun)} \, \| \DPN^{-1} \Av_{\nuiv} \|
	\leq 
	\crossgrad \, \| \HPN (\nuiv - \nuivs) \|_{\nano} \, .
\label{6dvgycuw76yeficf7wnfo}
\end{EQA}
This enables us to apply Proposition~\ref{PbiassemiN} which implies \eqref{rhDGtuGmusGU0a2sp} and
\begin{EQA}
	&& \nquad
	\| \DPN^{-1} \IFN_{\nuiv} (\targv_{\nuiv} - \targvs + \IFN_{\nuiv}^{-1} \Av_{\nuiv}) \|
	\leq 
	\frac{3 \dltwu_{3}}{4} \, \| \DPN \, \IFN_{\nuiv}^{-1} \Av_{\nuiv} \|^{2}
	\leq 
	\frac{\dltwu_{3} \, \crossgrad^{2}}{3} \, \, \| \HPN (\nuiv - \nuivs) \|_{\nano}^{2}
	\, .
\label{suhyd78wyedw378yedh2}
\end{EQA}
This completes the proof.
\end{proof}

Now we can finalize the proof of the main result.
With \( \AUv_{\nuiv} \eqdef \IFT_{\targv\nuiv} (\nuiv - \nuivs) \), it holds
\begin{EQA}
	&& \nquad
	\targv_{\nuiv} - \targvs + \IFN^{-1} \AUv_{\nuiv}
	=
	\targv_{\nuiv} - \targvs + \IFN_{\nuiv}^{-1} \Av_{\nuiv}
	- \IFN_{\nuiv}^{-1} \Av_{\nuiv}
	+ \IFN^{-1} \AUv_{\nuiv} 
	\\
	& = &
	\targv_{\nuiv} - \targvs + \IFN_{\nuiv}^{-1} \Av_{\nuiv}
	- \IFN_{\nuiv}^{-1} (\Av_{\nuiv} - \AUv_{\nuiv}) 
	- (\IFN_{\nuiv}^{-1} - \IFN^{-1}) \AUv_{\nuiv}
	\, .
\label{36gfijh94ejdvtwekoDg}
\end{EQA}
The use of \eqref{8c7c67cc6c63kdldlvvudw} and \eqref{f8vuehery6gv65ftehwee} yields
\begin{EQA}
	&& \nquad 
	\| \DPN^{-1} \IFN_{\nuiv} (\IFN^{-1} - \IFN_{\nuiv}^{-1}) \AUv_{\nuiv} \|
	= 
	\| \DPN^{-1} (\IFN - \IFN_{\nuiv}) \IFN^{-1} \AUv_{\nuiv} \|
	\\
	& \leq &
	\| \DPN^{-1} (\IFN - \IFN_{\nuiv}) \DPN^{-1} \| \,\, 
	\| \DPN \, \IFN^{-1} \DPN \| \,\, \| \DPN^{-1} \IFT_{\targv\nuiv} \HPN^{-1} \|_{\dual} \,\, \| \HPN (\nuiv - \nuivs) \|_{\nano}
	\\
	& \leq &
	\dltwunn \, \| \HPN (\nuiv - \nuivs) \|_{\nano} \, \,
	\crossnorm \, \| \HPN (\nuiv - \nuivs) \|_{\nano}.
	\qquad
\label{jufuf7yfy7tg643yrh}
\end{EQA}
This together with \eqref{4cbijjm09j9hhrjfifkrs} and \eqref{d73jejg7tre3fdy3fujell} implies
\begin{EQA}[c]
	\| \DPN^{-1} \IFN_{\nuiv} (\targv_{\nuiv} - \targvs + \IFN^{-1} \AUv_{\nuiv}) \|
	\leq 
	\Bigl( \frac{\dltwu_{3} \, \crossgrad^{2}}{3} \, 
	+ \crossnorm \, \dltwunn + \frac{\dltwun}{2}
	\Bigr) \| \HPN (\nuiv - \nuivs) \|_{\nano}^{2}
	\, .
\label{7tfce7sdfy7w7y73wrytk3f7}
\end{EQA}
The use of \eqref{uciuvcdiu3eoooodffcy} allows to bound
\begin{EQA}[c]
	\| \DPN^{-1} \IFN (\targv_{\nuiv} - \targvs + \IFN^{-1} \AUv_{\nuiv}) \|
	\leq
	\frac{1}{1 - \dltwbun} \, 
	\| \DPN^{-1} \IFN_{\nuiv} (\targv_{\nuiv} - \targvs + \IFN^{-1} \AUv_{\nuiv}) \| 
	\, ,
\label{9jjjjh3d4dd32x32sg3thj}
\end{EQA}
and \eqref{36gfijh94ejdvtwekoise} follows.

\paragraph{Proof of Theorem~\ref{Pconcsupp}}
Theorem~\ref{PFiWigeneric2} applied to \( \fs_{\nuiv}(\targv) \) and 
\( \fn_{\nuiv}(\targv) = \fs_{\nuiv}(\targv) - \langle \AAv, \targv \rangle \) implies  by \eqref{uciuvcdiu3eoooodffcym}
\begin{EQA}
	\| \DPN (\targvn_{\nuiv} - \targv_{\nuiv}) \|
	& \leq &
	\frac{3}{2} \| \DPN \, \IFN_{\nuiv}^{-1} \AAv \|
	\leq 
	\frac{3}{2(1 - \dltwbun)} \| \DPN \, \IFN^{-1} \AAv \| 
	\, .
\label{djvhehfief6yh3jfgiu2}
\end{EQA}
This and \eqref{rhDGtuGmusGU0a2spt} imply \eqref{ysd7euwejfg653wthf873k}.
Next we check \eqref{usdhyw6hikhurnetr256} using the decomposition
\begin{EQA}
	&& \nquad
	\targvn_{\nuiv} - \targvs + \IFN^{-1} \IFT_{\targv\nuiv} (\nuiv - \nuivs) 
	+ \IFN^{-1} \AAv
	\\
	&=&
	(\targvn_{\nuiv} - \targv_{\nuiv} + \IFN_{\nuiv}^{-1} \AAv) 
	- (\IFN_{\nuiv}^{-1} \AAv - \IFN^{-1} \AAv)
	+ \{ \targv_{\nuiv} - \targvs + \IFN^{-1} \IFT_{\targv\nuiv} (\nuiv - \nuivs) \} 
	\, .
\label{yhdt56swhew454erjhfgi}
\end{EQA}
Theorem~\ref{PsemiAvex} evaluates the last term.
Lemma~\ref{L3IFNNui} helps to bound
\begin{EQA}
	\| \DPN^{-1} \IFN (\IFN_{\nuiv}^{-1} - \IFN^{-1}) \IFN_{\nuiv} \, \DPN^{-1} \|
	& = &
	\| \DPN^{-1} (\IFN_{\nuiv} - \IFN) \DPN^{-1} \|
	\leq 
	\dltwunn \| \HPN (\nuiv - \nuivs) \|_{\nano} 
	\, .
\label{cuyedgfcuedgtfbvj}
\end{EQA}
This yields by \( \dltwbun \leq 1/4 \) and \( (4/3)ab \leq a^{2} + b^{2}/2 \) for any \( a,b \)
\begin{EQA}
	\| \QP (\IFN_{\nuiv}^{-1} - \IFN^{-1}) \AAv \|
	& \leq &
	\| \QP \, \IFN^{-1} \DPN \| 
	\, \| \DPN^{-1} \IFN_{\nuiv} (\IFN_{\nuiv}^{-1} - \IFN^{-1}) \IFN \, \DPN^{-1} \|
	\, \| \DPN \, \IFN_{\nuiv}^{-1} \AAv \|
	\\
	& \leq &
	\| \QP \, \IFN^{-1} \DPN \| 
	\, \dltwunn \, \| \HPN (\nuiv - \nuivs) \|_{\nano} \,
	\, \frac{1}{1 - \dltwbun} \| \DPN \, \IFN^{-1} \AAv \|
	\\
	& \leq &
	\| \QP \, \IFN^{-1} \DPN \| 
	\, \dltwunn \, \bigl( \| \HPN (\nuiv - \nuivs) \|_{\nano}^{2} + \| \DPN \, \IFN^{-1} \AAv \|^{2}/2 \bigr)
	\, .
\label{7duewjf87he4fyjwucyfne}
\end{EQA}
Moreover, 
for \( \epsv_{\nuiv} \eqdef \targvn_{\nuiv} - \targv_{\nuiv} + \IFN_{\nuiv}^{-1} \AAv \), it holds
\begin{EQA}
	\| \DPN^{-1} \IFN_{\nuiv} \, \epsv_{\nuiv} \|
	& \leq &
	\frac{3 \dltwu_{3}}{4} \| \DPN \, \IFN_{\nuiv}^{-1} \AAv \|^{2}  
	\leq 
	\frac{3 \dltwu_{3}}{4(1 - \dltwbun)^{2}} \| \DPN \, \IFN^{-1} \AAv \|^{2}
	\, ,
\label{ud82uegkftryyde6ewhde}
\end{EQA} 
and by \eqref{uciuvcdiu3eoooodffcym} and \( \dltwbun \leq 1/4 \)
\begin{EQA}
	\| \QP \, \epsv_{\nuiv} \|
	& \leq &
	\| \QP \, \IFN^{-1} \DPN \| \,\, \| \DPN^{-1} \IFN \, \epsv_{\nuiv} \|
	\leq 
	\| \QP \, \IFN^{-1} \DPN \| \, \frac{1 - \dltwbun}{1 - 2 \dltwbun} \, \| \DPN^{-1} \IFN_{\nuiv} \, \epsv_{\nuiv} \|
	\\
	& \leq &
	\| \QP \, \IFN^{-1} \DPN \| \,\, 
	\frac{3 \dltwu_{3}}{4(1 - \dltwbun) (1 - 2 \dltwbun)} \| \DPN \, \IFN^{-1} \AAv \|^{2} 
	\leq 
	\| \QP \, \IFN^{-1} \DPN \| \, 2 \dltwu_{3} \, \| \DPN \, \IFN^{-1} \AAv \|^{2}
	\, .
\label{ysdudhwuw8ey87wejue4bf}
\end{EQA}
The obtained bounds imply \eqref{usdhyw6hikhurnetr256}.

%% file: marg_opt_alt.tex

\def\DFN{\DPN}
\def\HFN{\HPN}

\Section{Alternating optimization}
\label{SoptAO}
Let \( \fs(\upsv) = \fs(\targv,\nuiv) \) depend on a target variable \( \targv \)
and the nuisance variable \( \nuiv \).
Consider the full-dimensional optimization problem:
\begin{EQA}[c]
	\upsvs
	=
	(\targvs,\nuivs)
	=
	\arginf_{(\targv,\nuiv)} \fs(\targv,\nuiv)
	\, .
\end{EQA}
This problem is not convex and can be hard to solve.
Instead, we assume that a partial optimization
of \( \fs(\targv,\nuiv) \) w.r.t. \( \targv \) with \( \nuiv \) fixed
and similarly w.r.t. \( \nuiv \) for a \( \targv \) fixed is a simple 
convex programming.
Suppose to be given a starting guess \( \targv_{0} \) and consider
the following alternating optimization (AO) procedure: for \( t \geq 1 \)
\begin{EQ}[rcccl]
	\nuiv_{t} 
	&=&
	\arginf_{\nuiv} \fs(\targv_{t-1},\nuiv)
	\, ,
	\qquad
	\targv_{t} 
	&=&
	\arginf_{\targv} \fs(\targv,\nuiv_{t})
	\, .
\label{tdydyc6cfw5cgedhrwfxx}
\end{EQ}
The question under study is a set of conditions ensuring a convergence 
of this procedure to the solution \( (\targvs,\nuivs) \) as \( t \) grows.
Let \( \IFT(\upsv) = \nabla^{2} \fs(\upsv) \) and \( \IFT = \IFT(\upsvs) \).
We will use the block representation \eqref{hwe78yf2diwe76tfw67etfwtbso}
of this matrix and assume that each diagonal 
block \( \IFT_{\targv\targv} \) and \( \IFT_{\nuiv\nuiv} \) is strongly positive.
Define
\begin{EQA}[c]
	\IFC \eqdef 
	\IFT_{\targv\targv}^{-1/2} \, \IFT_{\targv\nuiv} \, \IFT_{\nuiv\nuiv}^{-1/2} .
\label{jccybr9vgoey4rdy3di}
\end{EQA}
The major condition for applicability of the AO method is \( \| \IFC \, \IFC^{\T} \| < 1 \). 
First, we explain the basic idea 
for the case of a quadratic function \( \fs \) with 
\( \IFT \equiv \nabla^{2} \fs \).

\begin{lemma}
\label{LquadAO}
Let \( \fs \) be a quadratic function \( \fs \) with 
\( \IFT \equiv \nabla^{2} \fs \).
For any \( t \geq 1 \), 
the value \( \targv_{t} \) from \eqref{tdydyc6cfw5cgedhrwfxx} satisfies
with \( \IFC \) from \eqref{jccybr9vgoey4rdy3d}
\begin{EQA}[c]
	\IFT_{\targv\targv}^{1/2}(\targv_{t} - \targvs)
	=
	\IFC \, \IFC^{\T} \IFT_{\targv\targv}^{1/2} (\targv_{t-1} - \targvs)
	\, .
\end{EQA}
\end{lemma}

The proof follows from Lemma~\ref{Lpartmaxq} applied twice to 
\( \targv_{t} - \targvs \) and then to \( \nuiv_{t} - \nuivs \).
A successful application of this identity yields a linear convergence 
of \( \IFT_{\targv\targv}^{1/2}(\targv_{t} - \targvs) \) to zero with the rate \( \| \IFC \, \IFC^{\T} \| \).

We now present conditions ensuring a local convergence of the AO method in a general case.
These conditions mainly concern the marginal behavior of \( \fs \) in \( \targv \) 
with \( \nuiv \) fixed or another way around, 
and involve a metric tensor \( \DFN \) for \( \targv \)-component and 
\( \HFN \) for the \( \nuiv \)-component.
Consider a local set in \( \Ups \) of the form \( \Targ \times \Nui \),
where 
\begin{EQA}[rcccl]
	\Targ 
	&=& 
	\bigl\{ \targv \colon \| \DFN (\targv - \targvs) \| \leq \rrta \bigr\} 
	\, ,
	\qquad 
	\Nui 
	&=& 
	\bigl\{ \nuiv \colon \| \HFN (\nuiv - \nuivs) \| \leq \rrnu \bigr\} 
	\, .
\end{EQA}
Later we assume that \( \fs(\targv,\nuiv) \) as a function of \( \targv \in \Targ \)
is smooth and strongly convex in the sense ``\( \nabla^{2}_{\targv\targv} \fs(\targv,\nuiv) \) positive definite'' for \( \nuiv \in \Nui \) fixed,
and similarly in \( \nuiv \in \Nui \) for \( \targv \in \Targ \) fixed.
%

\begin{description}

    \item[\label{LLxoT3ref} \( \bb{(\mathcal{T}_{3,\Targ,\Nui}^{*})} \)]
      \emph{For any \( \nuiv \in \Nui \), 
      the function \( \fs(\targv,\nuiv) \) is strongly convex in 
      \( \targv \in \Targ \) and 
      it holds}
\begin{EQA}
	\sup_{\targv \in \Targ} \,\,
	\sup_{\uv \in \R^{\dimp}}  \,\, 
	\frac{\bigl| \langle \nabla_{\targv\targv\targv}^{3} \fs(\targv,\nuiv), \uv^{\otimes 3} \rangle \bigr|}
		 {\| \DFN \uv \|^{3}}
	& \leq &
	\dltwu_{3} \, .
\label{d6f53ye5vry4fddfgeyd}
\end{EQA}
      \emph{For any \( \targv \in \Targ \), 
      the function \( \fs(\targv,\nuiv) \) is strongly convex in \( \nuiv \) and 
      it holds}
\begin{EQA}
	\sup_{\nuiv \in \Nui} \,\,
	\sup_{\wv \in \R^{\dimq}}  \,\, 
	\frac{\bigl| \langle \nabla_{\nuiv\nuiv\nuiv}^{3} \fs(\targv,\nuiv), \wv^{\otimes 3} \rangle \bigr|}
		 {\| \HFN \wv \|^{3}}
	& \leq &
	\dltwu_{3} \, .
\label{d6f53ye5vry4fddfgeydr}
\end{EQA}
\emph{For any \( \upsv = (\targvs,\nuiv) \) with \( \nuiv \in \Nui \) or \( \upsv = (\targv,\nuivs) \) with \( \targv \in \Targ \),
it holds}
\begin{EQA}[rccl]
	\sup_{\zv \in \R^{\dimp}, \wv \in \R^{\dimq}} \,
	&
	\frac{| \langle \nabla_{\targv\nuiv\nuiv}^{3} \fs(\upsv), \zv \otimes \wv^{\otimes 2} \rangle |}
		 {\| \DFN \zv \| \, \| \HFN \wv \|^{2}}
	& \leq &
	\dltwun \, ,
	\qquad
\label{c6ceyecc5e5etcT2r}
	\\
	\sup_{\zv \in \R^{\dimp}, \wv \in \R^{\dimq}} \,
	&
	\frac{| \langle \nabla_{\targv\targv\nuiv}^{3} \fs(\upsv), \zv^{\otimes 2} \otimes \wv \rangle |}
		 {\| \DFN \zv \|^{2} \, \| \HFN \wv \|}
	& \leq &
	\dltwunn \, .
\label{c6ceyecc5e5etctwhcyegwcr}
\end{EQA}
\end{description}

If \( \dltwu_{3} \), \( \dltwun \), and \( \dltwunn \) are of the same order, then this condition is nearly equivalent 
to the full-dimensional condition \nameref{LLsT3ref}.
We will assume that \( \dltwun, \dltwu_{3} \) are sufficiently small to ensure 
\( (\dltwun + \dltwu_{3}) \, (\rrta + \rrnu) \ll 1 \).
For ease of formulation, later assume \( \dltwun = \dltwunn \).
With \( \IFC = \DFN^{-1} \IFT_{\targv\nuiv} \, \HFN^{-1} \) from \eqref{jccybr9vgoey4rdy3di}, 
define the technical values 
\begin{EQA}[rcl]
	\crossgrad
	& \eqdef &
	\frac{3}{2(1 - \dltwbun)}
	\Bigl( \| \IFC \| 
	+ \frac{\dltwbun}{2} 
	\Bigr) 
	\, ,
	\qquad
	\text{where } \quad
	\dltwbun \eqdef \dltwun \, (\rrnu \vee \rrta)
	\, ,
\label{f8vuehery6gv65ftehwee1r}
	\\
	\dltwuns
	& \eqdef &
	\frac{1}{1 - \dltwbun} \, 
	\biggl( \dltwun \, \| \IFC \| + \frac{\dltwun}{2}
		+ \frac{\dltwu_{3} \, \crossgrad^{2}}{3} 
	\biggr) \, ;
\label{fukllkf9de93yujrdee}
\end{EQA}
cf. 
\eqref{36gfijh94ejdvtwekoised}.
%
%
%
%
In the next result, we discuss a local convergence of the AO procedure 
started at a point \( \targv_{0} \in \Targ \).

\begin{theorem}
\label{TAOconv}
Assume \nameref{LLxoT3ref} with \( \DFN^{2} \leq \IFT_{\targv\targv} \), 
\( \HFN^{2} \leq \IFT_{\nuiv\nuiv} \),
and let \( \dltwbun = \dltwun \, (\rrnu \vee \rrta) < 1 \) and with \( \crossgrad \) from \eqref{f8vuehery6gv65ftehwee1r}
\begin{EQA}[c]
	\rrta
	\geq 
	\crossgrad^{2} \, 
	\| \DFN (\targv_{0} - \targvs) \| 
	\, ,
	\qquad
	\rrnu
	\geq 
	\crossgrad \, 
	\| \DFN (\targv_{0} - \targvs) \| 
	\, ,
	\qquad
	\crossgrad \, \dltwu_{3} \, (\rrta \vee \rrnu)
	\leq 
	\frac{2}{3} 
	\, .
	\qquad
\label{6cuydjd5eg3jfggu8eyt4h}
\end{EQA}
Let with \( \dltwuns \) from \eqref{fukllkf9de93yujrdee},
the operator \( \IFC = \IFT_{\targv\targv}^{-1/2} \, \IFT_{\targv\nuiv} \, \IFT_{\nuiv\nuiv}^{-1/2} \) satisfy
\begin{EQA}[c]
	\| \IFC \, \IFC^{\T} \| + ( 1 + \crossgrad^{2}) \, \dltwuns \, \| \DFN (\targv_{0} - \targvs) \| 
	\leq 
	\rhos
	<
	1 
	\, .
\label{hddctegu4357fvyehex}
\end{EQA}
Then \( \IFT_{\targv\targv}^{1/2} (\targv_{t} - \targvs) \) converges to zero 
at linear rate \( \rhos \).
\end{theorem}

\begin{remark}
Condition \eqref{hddctegu4357fvyehex} ensures a geometric convergence \( \targv_{t} \to \targvs \).
It the current value \( \targv_{t} \) is sufficiently close to \( \targvs \), 
then the result of the theorem can be restated for the starting point \( \targv_{t} \),
which leads to the convergence rate \( \rhos \) very close to \( \| \IFC \, \IFC^{\T} \| \).
\end{remark}

%% file: alternate_proof.tex
\def\HFN{\HPN}


\paragraph{Proof of Theorem~\ref{TAOconv}}
For every \( t \geq 1 \), by \eqref{rhDGtuGmusGU0a2spt} of Theorem~\ref{PsemiAvex} with \( \| \cdot \|_{\nano} = \| \cdot \| \)
\begin{EQ}[rcl]
	\| \HFN (\nuiv_{t} - \nuivs) \|
	& \leq &
	\crossgrad \, \| \DFN (\targv_{t-1} - \targvs) \| 
	\, ,
	\\
	\| \DFN (\targv_{t} - \targvs) \|
	& \leq &
	\crossgrad \, \| \HFN (\nuiv_{t} - \nuivs) \| 
	\, . 
\label{rhDGtuGmusGU0a2sptAO}
\end{EQ}
In particular, for \( t=1 \),
\begin{EQA}[rcl]
	\| \HFN (\nuiv_{1} - \nuivs) \|
	& \leq &
	\crossgrad \, \| \DFN (\targv_{0} - \targvs) \| 
	\, .
\end{EQA}
Further, with \( \IFC = \IFT_{\targv\targv}^{-1/2} \, \IFT_{\targv\nuiv} \, \IFT_{\nuiv\nuiv}^{-1/2} \), denote
\begin{EQA}[rcccl]
	\epsv_{t}
	& \eqdef &
	\IFT_{\targv\targv}^{1/2} (\targv_{t} - \targvs)
	&+& \IFC \IFT_{\nuiv\nuiv}^{1/2} (\nuiv_{t} - \nuivs)
	\, ,
	\\
	\alpv_{t}
	& \eqdef &
	\IFT_{\nuiv\nuiv}^{1/2} (\nuiv_{t} - \nuivs)
	&+& \IFC^{\T} \IFT_{\targv\targv}^{1/2} (\targv_{t-1} - \targvs)
	\, .
\end{EQA}
Then 
\begin{EQA}[rcccl]
	\IFT_{\targv\targv}^{1/2} (\targv_{t} - \targvs)
	&-& \IFC \, \IFC^{\T} \IFT_{\targv\targv}^{1/2} (\targv_{t-1} - \targvs)
	&=&
	\epsv_{t} - \IFC \alpv_{t}
	\, ,
	\\
	\IFT_{\nuiv\nuiv}^{1/2} (\nuiv_{t} - \nuivs)
	&-& \IFC^{\T} \IFC \, \IFT_{\nuiv\nuiv}^{1/2} (\nuiv_{t-1} - \nuivs)
	&=&
	\alpv_{t} - \IFC^{\T} \epsv_{t-1}
	\, .
\label{dgvi3erd7buyjk32rr433}
\end{EQA}
Expansion \eqref{36gfijh94ejdvtwekoise} of Theorem~\ref{PsemiAvex} with \( \QP = \IFT_{\targv\targv}^{1/2} \) yields by \eqref{rhDGtuGmusGU0a2sptAO}
\begin{EQA}[rcl]
\label{36gfijh94ejdvtwekoiser}
	\| \epsv_{t} \|
	& \leq &
	\dltwuns \, \| \HFN (\nuiv_{t} - \nuivs) \|^{2}
	\leq 
	\dltwuns \, \crossgrad^{2} \, 
	\| \DFN (\targv_{t-1} - \targvs) \|^{2} 
	\, .
	\qquad
\label{36gfijh94ejdvtwekoisedr}
\end{EQA}
A similar result applies to \( \alpv_{t} \):
\begin{EQA}[rcl]
	\| \alpv_{t} \|
	& \leq &
	\dltwuns \, \| \DFN (\targv_{t-1} - \targvs) \|^{2}
	\leq 
	\dltwuns \, \crossgrad^{2} \, 
	\| \HFN (\nuiv_{t-1} - \nuivs) \|^{2} 
	\, .
	\qquad
\end{EQA}
As \( \| \IFC \| < 1 \) and \( \DFN^{2} \leq \IFT_{\targv\targv} \), 
we conclude that
\begin{EQA}[rcl]
	&& \nquad
	\bigl\| \IFT_{\targv\targv}^{1/2} (\targv_{t} - \targvs)
		- \IFC \, \IFC^{\T} \IFT_{\targv\targv}^{1/2} (\targv_{t-1} - \targvs) 
	\bigr\|
	=
	\| \epsv_{t} - \IFC \alpv_{t} \|
	\\
	& \leq &
	( 1 + \crossgrad^{2}) 
	\dltwuns \| \DFN (\targv_{t-1} - \targvs) \|^{2}
	\leq 
	( 1 + \crossgrad^{2} ) 
	\dltwuns \, \| \DFN (\targv_{t-1} - \targvs) \| \,
	\| \IFT_{\targv\targv}^{1/2} (\targv_{t-1} - \targvs) \|
	\, 
\label{dgvi3erd7buyjk32rr433AO}
\end{EQA}
yielding 
\begin{EQA}[c]
	\| \IFT_{\targv\targv}^{1/2} (\targv_{t} - \targvs) \|
	\leq 
	\Bigl\{ \| \IFC \IFC^{\T} \| + (1 + \crossgrad^{2}) \, \dltwuns 
		\| \DFN (\targv_{t-1} - \targvs) \| 
	\Bigr\} 
	 \, \| \IFT_{\targv\targv}^{1/2} (\targv_{t-1} - \targvs) \|
	 \, .
\end{EQA}
This and \eqref{hddctegu4357fvyehex} imply linear convergence of 
\( \| \IFT_{\targv\targv}^{1/2} (\targv_{t} - \targvs) \| \) to zero
at the rate \( \| \IFC \IFC^{\T} \| \).
Similarly, 
\begin{EQA}[rcl]
	&& \nquad
	\bigl\| \IFT_{\nuiv\nuiv}^{1/2} (\nuiv_{t+1} - \nuivs)
		- \IFC^{\T} \IFC \, \IFT_{\nuiv\nuiv}^{1/2} (\nuiv_{t} - \nuivs) 
	\bigr\|
	=
	\| \alpv_{t+1} - \IFC^{\T} \epsv_{t} \|
	\\
	& \leq &
	( 1 + \crossgrad^{2} ) \, \dltwuns \| \HFN (\nuiv_{t} - \nuivs) \|^{2}
	\leq 
	( 1 + \crossgrad^{2}) \,
	\dltwuns \| \HFN (\nuiv_{t} - \nuivs) \| 
	 \, \| \IFT_{\nuiv\nuiv}^{1/2} (\nuiv_{t} - \nuivs) \|
	\, .
\label{dgvi3erd7buyjk32rr433AOnu}
\end{EQA}
This ensures a linear convergence of 
\( \| \IFT_{\nuiv\nuiv}^{1/2} (\nuiv_{t} - \nuivs) \| \) to zero
at the same rate \( \| \IFC \IFC^{\T} \| \).

%% file: marg_opt_sup.tex

\Section{Sup-norm expansions in perturbed optimization}
\label{Ssupnormco}
For a convex smooth function \( \fs(\upsv) \) and a perturbation \( \langle \AAv,\upsv \rangle \),
consider
\begin{EQA}[c]
	\upsvs 
	= \arginf_{\upsv} \fs(\upsv)
	\, ,
	\qquad
	\upsvd = \arginf_{\upsv} \fd(\upsv)
	=
	\arginf_{\upsv} \bigl\{ \fs(\upsv) + \langle \AAv,\upsv \rangle \bigr\}
	\, .
\end{EQA} 
This section studies the componentwise difference \( \upsvd - \upsvs \).
The results from \cite{Sp2024PO} on perturbed optimization are coordinate-free.
A choice of a basis in the parameter space plays no role. 
This is, however, not the case in the sup-norm estimation. 
A different approach and another proof strategy are called for. 
The basic idea is to consider the sup-norm estimation 
as a kind of linewise optimization with one fixed component as a target
and the remaining components as nuisance. 


Define
\( \IFT = \nabla^{2} \fs(\upsvs) = (\IFT_{jm}) \) and
\( \DF^{2} = \diag(\IFT) \):
\begin{EQA}[c]
	\DF^{2} \eqdef \diag(\DPN_{1}^{2}, \ldots, \DPN_{\dimp}^{2})
	=
	\diag(\IFT_{11}, \ldots, \IFT_{\dimp\dimp})
	\, .
\label{7dycf8mwdy6ey43hf98yh}
\end{EQA}
Define \( \zzinf \eqdef \| \DF^{-1} \AAv \|_{\infty} \) and,
with \( \BBF \eqdef \DF^{-1} \IFT \, \DF^{-1} \), introduce the values \( \crosssup \) by
\begin{EQA}[rclcl]
	\crosssup
	& \eqdef &
	\sup_{\| \uv \|_{\infty} \leq 1} \| (\BBF - \Id_{\dimp}) \uv \|_{\infty}
	&=&
	\| \BBF - \Id_{\dimp} \|_{1}
	=
	\max_{j = 1,\ldots,\dimp} \frac{1}{\DPN_{j}}
	\sum_{m \neq j} \frac{|\IFT_{jm}|}{\DPN_{m}}
	\, .
	\qquad
\label{7tdsyf8iuwopkrtg4576}
\end{EQA}
The proposed approach requires \( \crosssup < 1 \).
Also, define \( \rrinf \eqdef 3 \zzinf/(1 - \crosssup) \) and introdice the local rectangle set
\begin{EQA}[c]
	\Upsd = \{ \upsv \colon \| \DF (\upsv - \upsvs) \|_{\infty} \leq \rrinf \} 
	=
	\Bigl\{ \upsv \colon \max_{j \leq \dimp} \bigl| \DPN_{j} (\ups_{j} - \upss_{j}) \bigr| 
	\leq \rrinf \Bigr\}.
\label{du7duydy7cy6fc6df3}
\end{EQA}
We will show that under mild smoothness conditions on \( \fs(\upsv) \), the point \( \upsvd \) belongs to \( \Upsd \).
The formulation uses
for each \( j \leq \dimp \) the representation \( \upsv = (\ups_{j},\nuiv_{j}) \), where 
\( \nuiv_{j} \in \R^{\dimp-1} \) collects all the remaining entries of \( \upsv \).
Similarly, denote \( \DF = (\DPN_{j},\HPN_{j}) \) with \( \HPN_{j} = \diag(\DPN_{m}) \) for \( m \neq j \).
%
Assume the following condition.

\begin{description}
    \item[\label{LLpsupref} \( \bb{(\mathcal{T}_{\infty}^{*})} \)]
      \emph{For each \( j \leq \dimp \), 
     the function \( \fs(\ups_{j},\nuiv_{j}) \) fulfills}
\begin{EQA}[rcl]
	\sup_{\nuiv_{j} \colon \| \HPN_{j} (\nuiv_{j} - \nuivs_{j}) \|_{\infty} \leq \rrinfi} \,\, 
	\sup_{\ups \colon \DPN_{j} |\ups - \upss_{j}| \leq \rrinfi} \,
	\frac{\bigl| \nabla_{\ups_{j}\ups_{j} \ups_{j}}^{(3)} \fs(\ups_{j},\nuiv_{j}) \bigr|}
		 {\DPN_{j}^{3}}
	& \leq &
	\dltwu_{3} \, ,
\label{d6f53ye5vry4fddfgsup}
\end{EQA}
\emph{and}
\begin{EQA}[lccl]
	\sup_{\upsv = (\upss_{j},\nuiv_{j}) \in \Upsd} \,
	\sup_{\zv_{j} \in \R^{\dimp-1}} \,
	&
	\frac{| \langle \nabla_{\ups_{j}\ups_{j}\nuiv_{j}}^{(3)} \fs(\upss_{j},\nuiv_{j}), \zv_{j} \rangle |}
		 {\DPN_{j}^{2} \, \| \HPN_{j} \zv_{j} \|_{\infty}}
	& \leq &
	\dltwunn 
	\, ,
\label{c6hyi8ietctwhcsup}
	\\
	\sup_{\upsv = (\upss_{j},\nuiv_{j}) \in \Upsd}  \,
	\sup_{\zv_{j} \in \R^{\dimp-1}} \,
	&
	\frac{| \langle \nabla_{\ups_{j}\nuiv_{j}\nuiv_{j}}^{3} \fs(\upss_{j},\nuiv_{j}), \zv_{j}^{\otimes 2} \rangle |}
		 {\DPN_{j} \, \| \HPN_{j} \zv_{j} \|_{\infty}^{2}}
	& \leq &
	\dltwun 
	\, .
	\qquad
\label{c6ceyecc5e5etcT2sup}
\end{EQA}
\end{description}

In a sense, this condition simply introduces the values \( \dltwun , \dltwunn , \dltwu_{3} \), bounding 
directional third-order derivatives of \( \fs \).
Before presenting the concentration result, define some technical values 
\( \dltwbun \eqdef \dltwunn \, \rrinf \) 
(later we assume \( \dltwbun \leq 1/4 \)) and 
\begin{EQA}[rcl]
	\dltwuns
	& \eqdef &
	\frac{1}{1 - \dltwbun} \, 
	\biggl( \crosssup \, \dltwunn + \frac{\dltwun}{2}
		+ \frac{3(\crosssup + \dltwbun/2)^{2} \, \dltwu_{3}}{4(1 - \dltwbun)^{2}} 
	\biggr) 
	\, ,
	\qquad
\label{36gfijh94ejdvtwekoisedi}
	\\
	\alp_{0}
	& \eqdef &
	(2 \dltwu_{3} + \dltwunn/2) \, \zzinf 
	\, ,
	\qquad
	\alp_{2} \eqdef (\dltwuns + \dltwunn) \, \zzinf
	\, ;
\label{kkoodfyf6e63tyhgds}
\end{EQA}
cf. \eqref{36gfijh94ejdvtwekoised}.
In typical situations with \( \dltwunn \leq \dltwu_{3} \), and \( \dltwu_{3} \, \rrinf \leq 1/4 \), it holds
\( \dltwuns \leq 4 \dltwu_{3} \).
Informally, the required smoothness condition means \( \alp_{0} \ll 1 \) and \( \alp_{2} \ll 1 \).
%

\begin{theorem}
\label{PsemibiassupC}
Let 
\par\noindent \( \bullet \)
\( \fs(\upsv) \) be strongly convex, \( \upsvs = \arginf_{\upsv} \fs(\upsv) \),
\( \IFT = \nabla^{2} \fs(\upsvs) \), and \( \DF^{2} = \diag(\IFT) \);

\par\noindent \( \bullet \)
\( \fd(\upsv) = \fs(\upsv) + \langle \AAv,\upsv \rangle \), where \( \AAv \) fulfills 
\( \| \DF^{-1} \AAv \|_{\infty} \leq \zzinf \);

\par\noindent \( \bullet \)
with \( \BBF \eqdef \DF^{-1} \IFT \, \DF^{-1} \), 
it holds \( \crosssup = \| \BBF - \Id_{\dimp} \|_{1} < 1 \); 
see \eqref{7tdsyf8iuwopkrtg4576}.

\par\noindent \( \bullet \)
\nameref{LLpsupref} hold on \( \Upsd \) from \eqref{du7duydy7cy6fc6df3} with \( \rrinf = \frac{3 \zzinf}{1 - \crosssup} \),
and \( \alp_{0} \) and \( \alp_{2} \) from \eqref{kkoodfyf6e63tyhgds} satisfy
\begin{EQA}[c]
	\alp_{0}
	\leq 
	1/2
	\, ,
	\qquad
	4 \alp_{2} (1 + \alp_{0})
	\leq 
	(1 - \crosssup)^{2}
	\, .
	\qquad
\label{dfujkkkhgfe32333}
\end{EQA}
Then \( \upsvd \eqdef \arginf_{\upsv} \fd(\upsv) \in \Upsd \) or, equivalently, 
\begin{EQA}[c]
	\| \DF(\upsvd - \upsvs) \|_{\infty} 
	\leq 
	\rrinf 
	= 
	3 \zzinf/(1 - \crosssup) .
\end{EQA}
\end{theorem}


\noindent\textbf{Sup-norm expansion for a linear perturbation.}
Here, we describe the sup-norm accuracy of approximation 
\( \upsvd - \upsvs \approx \IFT^{-1} \AAv \)
and \( \upsvd - \upsvs \approx \DF^{-2} \AAv \) with \( \DF^{2} = \diag(\IFT) \).
With \( \BBF = \DF^{-1} \IFT \DF^{-1} \), define the value \( \crossB \) by
\begin{EQA}[rcl]
	\crossB^{-1}
	& = &
	\sup_{\| \uv \|_{\infty} \leq 1} \| \BBF^{-1} \uv \|_{\infty} 
	=
	\| \BBF^{-1} \|_{1}
	\, .
\label{hdcukekei76erddfree}
\end{EQA}
By Lemma~\ref{LinvGersh} ahead, it holds \( \crossB \geq 1 - \crosssup \).
However, it is a rough lower bound.
Also, with \( \dltwuns \) from \eqref{36gfijh94ejdvtwekoisedi}, 
and \( \CONSTzz = 3/(1 - \crosssup) \), introduce
\begin{EQA}[c]
	\dltwu_{\CONSTzz}
	\eqdef
	2 \dltwu_{3} + \frac{\dltwunn}{2} + \CONSTzz^{2} (\dltwuns + \dltwunn)   
	\, .
	\qquad
\label{usdhyw6hikhurnetrspp2}
\end{EQA}

\begin{theorem}
\label{Psemibiassup}
Under the conditions of Theorem~\ref{PsemibiassupC}, it holds
with \( \dltwu_{\CONSTzz} \) from \eqref{usdhyw6hikhurnetrspp2}
\begin{EQ}[c]
\label{usdhyw6hikhurnetrspp}
	\bigl\| \DF^{-1} \bigl\{ \IFT (\upsvd - \upsvs) + \AAv \bigr\} \bigr\|_{\infty}
	\leq 
	\dltwu_{\CONSTzz} \, \| \DF^{-1} \AAv \|_{\infty}^{2}
	\, .
\end{EQ}
Furthermore, with \( \BBF = \DF^{-1} \IFT \DF^{-1} \) and \( \crossB \) from \eqref{hdcukekei76erddfree},
\begin{EQA}[c]
\label{usdhyw6hikhurnetrsppD}
	\| \DF (\upsvd - \upsvs) + \DF^{-1} \AAv \|_{\infty} 
	\leq 
	\crossB^{-1} \bigl( \| (\BBF - \Id_{\dimp}) \, \DF^{-1} \AAv \|_{\infty} + \dltwu_{\CONSTzz} \, \| \DF^{-1} \AAv \|_{\infty}^{2} \bigr)
	\, .
	\qquad
\end{EQA}
\end{theorem}


\begin{remark}
Expansion \eqref{usdhyw6hikhurnetrsppD} is meaningful only if the remainder in its right-hand side 
is small compared with the leading term \( \| \DF^{-1} \AAv \|_{\infty} \).
Let all the smoothness constants \( \dltwun, \dltwunn,\dltwuns \) are of order 
\( \dltwu_{3} \).
In addition to the conditions \( \crosssup < 1 \) and \( \dltwu_{3} \, \zzinf \ll 1 \) required for the componentwise consistency result of Theorem~\ref{PsemibiassupC}, we need that 
\( \| (\BBF - \Id_{\dimp}) \DF^{-1} \AAv \|_{\infty} \ll \| \DF^{-1} \AAv \|_{\infty} \) and 
\( \crossB \) is separated away from zero.
\end{remark}


\paragraph{A  separable perturbation}
A linear perturbation is a special case of a separable perturbation of the form 
\( \pert(\upsv) = \sum_{j} \pert_{j}(\ups_{j}) \).
A popular example is a Tikhonov regularization with \( \pert_{j}(\ups_{j}) = \lambda \ups_{j}^{2} \) or, 
more generally, a Sobolev-type penalty \( \pert_{j}(\ups_{j}) = \gp_{j}^{2} \ups_{j}^{2} \) for a given sequence 
\( \gp_{j}^{2} \).
Later we assume each component \( \pert_{j}(\ups_{j}) \) to be sufficiently smooth.
This excludes the cases of a sparse penalty \( \pert_{j}(\ups_{j}) = \lambda |\ups_{j}| \) or
complexity penalty \( \pert_{j}(\ups_{j}) = \lambda \Ind(\ups_{j} \neq 0) \).
Let \( \fs(\upsv) \) be convex, \( \upsvs = \arginf_{\upsv} \fs(\upsv) \).
Define a separable perturbation \( \fd(\upsv) = \fs(\upsv) + \pert(\upsv) \) of \( \fs(\upsv) \) by smooth functions \( \pert_{j}(\ups_{j}) \):
\begin{EQA}
	\fd(\upsv)
	&=&
	\fs(\upsv) + \pert(\upsv)
	=
	\fs(\upsv) + \sum_{j=1}^{\dimp} \pert_{j}(\ups_{j})
	\, , 
\label{c7c7ctefwct5enfvyewsd}
\end{EQA} 
and let \( \upsvn = \arginf_{\upsv} \fd(\upsv) \).
We intend to bound the corresponding change \( \upsvn - \upsvs \) in a sup-norm.
Separability of the penalty is very useful; the cross derivatives of \( \fd \) are the same as for \( \fs \).
We only update the tensor \( \DF(\upsv) \) assuming convexity of each \( \pert_{j}(\cdot) \):
\begin{EQA}
	\DF(\upsv)
	& \eqdef &
	\diag\{ \DPN_{1}(\upsv),\ldots,\DPN_{\dimp}(\upsv) \},
	\qquad
	\DPN_{j}^{2}(\upsv)
	\eqdef
	\IF_{jj}(\upsv) + \pert''_{j}(\ups_{j}) \, .
\label{7fyc5wgdf7vrewgdhfywh}
\end{EQA}

\begin{theorem}
\label{Psemibiassupp}
Let 
\par\noindent \( \bullet \)
the solution \( \upsvs = \arginf_{\upsv} \fs(\upsv) \) exist and unique;

\par\noindent \( \bullet \)
\( \fd(\upsv) \) from \eqref{c7c7ctefwct5enfvyewsd} be strongly convex;

\par\noindent \( \bullet \)
\( \upsvn = \arginf_{\upsv} \fd(\upsv) \),
\( \DF = \DF(\upsvn) \),
\( \AAv \eqdef \nabla \pert(\upsvs) = (\pert'_{j}(\upss_{j})) \),
\( \| \DF^{-1} \AAv \|_{\infty} \leq \zzinf \);

\par\noindent \( \bullet \) 
\( \fd(\upsv) \) satisfy \nameref{LLpsupref} with  \( \rrinf = \CONSTzz \, \zzinf \) and \( \CONSTzz = 3 / (1 - \crosssup) \),
and conditions \eqref{dfujkkkhgfe32333} hold.

Then with \( \dltwu_{\CONSTzz} \) from \eqref{usdhyw6hikhurnetrspp2},
it holds
\begin{EQA}[rcl]
	\| \DF(\upsvd - \upsvs) \|_{\infty} 
	& \leq &
	\rrinf
	=
	\CONSTzz \, \zzinf
	\, ,
	\\
	\bigl\| \DF^{-1} \bigl\{ \IFT (\upsvd - \upsvs) + \AAv \bigr\} \bigr\|_{\infty}
	& \leq &
	\dltwu_{\CONSTzz} \, \zzinf^{2}
	\, ,
	\\
	\| \DF (\upsvd - \upsvs) + \DF^{-1} \AAv \|_{\infty} 
	& \leq &
	\crossB^{-1} \bigl( \| (\BBF - \Id_{\dimp}) \, \DF^{-1} \AAv \|_{\infty} + \dltwu_{\CONSTzz} \, \zzinf^{2} \bigr)
	\, .
	\qquad
\end{EQA}
\end{theorem}

\begin{remark}
\label{RLLsup}
If each \( \pert_{j}(\ups_{j}) \) is quadratic, then condition \nameref{LLpsupref} for \( \fd(\upsv) \) follows from 
the same condition for \( \fs(\upsv) \) as the third derivatives of these two functions coincide.
For general smooth perturbations \( \pert_{j}(\cdot) \), \eqref{d6f53ye5vry4fddfgsup} should be checked
for \( \nabla^{3}_{\ups_{j}\ups_{j} \ups_{j}} \fd(\upsv) = \nabla^{3}_{\ups_{j}\ups_{j} \ups_{j}} \fs(\upsv) + \pert'''_{j}(\ups_{j}) \).
The other conditions in \nameref{LLpsupref} can be checked for \( \fs \).
\end{remark}


%% file: marg_opt_sup_proofs.tex

\def\nuivd{\nuiv^{\circ}}
\def\IFN{\IF}


\paragraph{Proof of Theorem~\ref{PsemibiassupC}}
With \( \| \DF^{-1} \AAv \|_{\infty} = \zzinf \), define \( \CONSTzz \) by the quadratic equation
\begin{EQA}[c]
	1 + \crosssup \, \CONSTzz + (\dltwuns + \dltwunn) \, \CONSTzz^{2} \zzinf
	+ (2 \dltwu_{3} + \dltwunn/2) \, \zzinf
	= 
	\CONSTzz
	\, ,
	\qquad
\label{jvikmmkkkide4iffde3}
\end{EQA}
or, with \( \alp_{0} \eqdef (2 \dltwu_{3} + \dltwunn/2) \, \zzinf \), 
\( \alp_{2} \eqdef (\dltwuns + \dltwunn) \, \zzinf \)
\begin{EQA}[c]
	\alp_{2} \, \CONSTzz^{2} - (1 - \crosssup) \CONSTzz + 1 + \alp_{0}
	=
	0
	\, .
\end{EQA}
In view of \eqref{dfujkkkhgfe32333},
\( (1 - \crosssup)^{2} \geq 4 \alp_{2} (1 + \alp_{0}) \), and the smallest root is given by
\begin{EQA}[c]
	\CONSTzz
	=
	\frac{1 - \crosssup - \sqrt{(1 - \crosssup)^{2} - 4 \alp_{2} (1 + \alp_{0})}}{2 \alp_{2}}
	=
	\frac{2 (1 + \alp_{0})}{1 - \crosssup + \sqrt{(1 - \crosssup)^{2} - 4 \alp_{2} (1 + \alp_{0})}}
\end{EQA}
and \( \alp_{0} \leq 1/2 \) implies \( \CONSTzz \leq 3/(1 - \crosssup) \).
Let the set \( \Upsd \) be fixed as in \eqref{du7duydy7cy6fc6df3} with \( \rrinf = \CONSTzz \, \zzinf \).
We aim at showing that \( \upsvd = \arginf_{\upsv} \fn(\upsv) \in \Upsd \).
Strong convexity of \( \fs(\upsv) \) and hence, of \( \fn(\upsv) \), ensures that the solution \( \upsvd \) is unique.
We apply a special alternating procedure with the starting point at \( \upsvs \):
at each step, only one component \( \ups_{j} \) is optimized while the remaining ones are kept fixed.
We show by induction that if the current iterate \( \upsv^{(k)} \) belongs to \( \Upsd \) then the same applies
to \( \upsv^{(k+1)} \).
In other words, the procedure never leaves \( \Upsd \) and, by strong convexity, it converges to the global solution \( \upsvd \in \Upsd \).
Obviously, it suffices to consider only one iterate.
Let us fix any \( \upsv \in \Upsd \), an index \( j \leq \dimp \), e.g. \( j=1 \), and consider the partial optimization problem 
w.r.t. \( \ups_{1} \) while the remaining entries of \( \upsv \) are kept fixed.
Represent \( \upsv = (\ups_{1},\nuiv_{1}) \), where \( \nuiv_{1} = (\ups_{2},\ldots,\ups_{\dimp})^{\T} \),
and \( \DF = (\DPN_{1},\HPN) \). 
For any \( \nuiv_{1} \) with \( \| \HPN (\nuiv_{1} - \nuivs_{1}) \|_{\infty} \leq \rrinf \), define
\begin{EQA}
	\upsd_{1}(\nuiv_{1})
	& \eqdef &
	\arginf_{\ups_{1}} \fn(\ups_{1},\nuiv_{1}) .
\label{fifwi9jhrio3jiofcwgyt7we}
\end{EQA}
We check that \( (\upsd_{1}(\nuiv_{1}),\nuiv_{1}) \in \Upsd \) or, equivalently,
\begin{EQA}[c]
	\bigl| \DPN_{1} \{ \upsd_{1}(\nuiv_{1}) - \upss_{1} \} \bigr|
	\leq 
	\rrinf
	=
	\CONSTzz \, \zzinf
	\, .
\end{EQA}
Bound \eqref{usdhyw6hikhurnetr256} of Theorem~\ref{Pconcsupp} with \( \QP = \DPN_{1} = \IFN_{11}^{1/2} \) implies
by \( \| \DPN_{1}^{-1} \AAv_{1} \|_{\infty} \leq \zzinf \) and \( \| \HPN^{-1} (\nuiv_{1} - \nuivs_{1}) \|_{\infty} \leq \CONSTzz \, \zzinf \)
\begin{EQA}[rcl]
	\bigl| \DPN_{1} \{ \upsd_{1}(\nuiv_{1}) - \upss_{1} + \IFN_{11}^{-1} \IFT_{\targv\nuiv} (\nuiv_{1} - \nuivs_{1}) + \IFN_{11}^{-1} \AAv_{1} \} \bigr|
	& \leq &
	(\dltwuns + \dltwunn) \, \CONSTzz^{2} \zzinf^{2}
		+ (2 \dltwu_{3} + \dltwunn/2) \, \zzinf^{2}
\end{EQA}
and by \( \| \DFN_{1}^{-1} \IFT_{\targv\nuiv} (\nuiv_{1} - \nuivs_{1}) \|_{\infty} \leq \crosssup \| \HPN (\nuiv_{1} - \nuivs_{1}) \|_{\infty} 
\leq \crosssup \, \CONSTzz \, \zzinf \) and \eqref{jvikmmkkkide4iffde3}
\begin{EQA}[rcl]
	\bigl| \DPN_{1} \{ \upsd_{1}(\nuiv_{1}) - \upss_{1} \} \bigr|
	& \leq &
	\zzinf \bigl\{ 1 + \crosssup \, \CONSTzz + (\dltwuns + \dltwunn) \, \CONSTzz^{2} \zzinf
	+ (2 \dltwu_{3} + \dltwunn/2) \, \zzinf \bigr\} 
	\leq 
	\CONSTzz \, \zzinf
	\, .
\end{EQA}
This yields the assertion of the theorem.

\paragraph{Proof of Theorem~\ref{Psemibiassup}}
As in the proof of Theorem~\ref{PsemibiassupC}, assume \( \| \DF^{-1} \AAv \|_{\infty} = \zzinf \).
Let us fix any \( \upsv \in \Upsd \), an index \( j \leq \dimp \), e.g. \( j=1 \), and consider the partial optimization problem 
w.r.t. \( \ups_{1} \) with the remaining entries of \( \upsv \) fixed.
Represent \( \upsv = (\ups_{1},\nuiv_{1}) \), where \( \nuiv_{1} = (\ups_{2},\ldots,\ups_{\dimp})^{\T} \). 
We apply Theorem~\ref{Pconcsupp} 
with \( \QP = \DPN_{1} \), \( \IFN = \IFT_{11} \), \( \HPN = \diag(\DPN_{2},\ldots,\DPN_{\dimp}) \),
and \( \| \HPN (\nuiv_{1} - \nuivs_{1}) \|_{\infty} \leq \rrinf \leq \CONSTzz \, \zzinf \).
Note first that \( \crossnorm \) from \eqref{f8vuehery6gv65ftehwee} coincides with \( \crosssup \) \eqref{7tdsyf8iuwopkrtg4576}.
Indeed, 
\begin{EQA}[c]
	\| \DPN_{1}^{-1} \IFT_{\targv\nuiv} \, \HPN^{-1} \|_{\dual}
	=
	\inf_{\| \zv \|_{\infty} \leq 1}
	\biggl| \sum_{m > 1} \frac{1}{\DPN_{1} \, \DPN_{m}} \, \IFT_{1m} z_{m} \biggr| 
	= 
	\crosssup 
	\, .
\end{EQA}
Bound \eqref{usdhyw6hikhurnetr256} yields by \( \DPN_{1} \IFT_{11}^{-1} = \DPN_{1}^{-1} \), 
\( \| \DF (\upsv - \upsvs) \|_{\infty} \leq \CONSTzz \, \zzinf \), and \( \| \DF^{-1} \AAv \|_{\infty} = \zzinf \)
\begin{EQA}
	&& \nquad
	\bigl| \DPN_{1} \{ \upsd_{1}(\nuiv_{1}) - \upss_{1} + \IFT_{11}^{-1} \AAv_{1} + \IFT_{11}^{-1} \IFT_{\ups_{1}\nuiv_{1}}(\nuiv_{1} - \nuivs_{1}) \} \bigr|
	\\
	& \leq &
	(\dltwuns + \dltwunn) \, \| \HPN (\nuiv_{1} - \nuivs_{1}) \|_{\infty}^{2}
		+ \bigl(2 \dltwu_{3} + \frac{\dltwunn}{2} \bigr) \, \DPN_{1}^{-2} \AAv_{1}^{2}
	\\
	& \leq &
	(\dltwuns + \dltwunn) \| \DF (\upsv - \upsvs) \|_{\infty}^{2} + (2 \dltwu_{3} + \frac{\dltwunn}{2} ) \, \| \DF^{-1} \AAv \|_{\infty}^{2} 
	\leq 
	\dltwu_{\CONSTzz} \, \zzinf^{2} 
	\, .
	\qquad
\label{usdhyw6hikhurnetr256sup}
\end{EQA}
The use of \( \nuiv_{1} = \nuivd_{1} \) yields \( \upsv = \upsvd \), \( \upsd_{1}(\nuivd_{1}) = \upsd_{1} \), and 
\begin{EQA}[rcl]
	&& \nquad
	\bigl| \DPN_{1} \{ \upsd_{1} - \upss_{1} + \IFT_{11}^{-1} \AAv_{1} + \IFT_{11}^{-1} \IFT_{\ups_{1}\nuiv_{1}}(\nuivd_{1} - \nuivs_{1}) \} \bigr|
	\\
	& \leq &
	\Bigl| \DPN_{1}^{-1} \bigl\{ \IFT_{11} (\upsd_{1} - \upss_{1}) + \IFT_{\ups_{1}\nuiv_{1}}(\nuivd_{1} - \nuivs_{1}) + \AAv_{1} \bigr\} \Bigr|
	\leq 
	\dltwu_{\CONSTzz} \, \zzinf^{2}
	\, .
\end{EQA}
Applying a similar bound to each component \( \upsvd_{j} \) for \( j \leq \dimp \) leads to
\begin{EQA}[c]
	\bigl\| \DF^{-1} \{ \IFT (\upsvd - \upsvs) + \AAv \} \bigr\|_{\infty}
	\leq 
	\dltwu_{\CONSTzz} \, \zzinf^{2}
	\, .
\label{usdhyw6hikhurnetrspp2p}
\end{EQA}
This completes the proof of the first statement.
%

Further, with \( \DeltaB \eqdef \DF^{-1} \IFT \, \DF^{-1} - \Id_{\dimp} \), 
bound \eqref{usdhyw6hikhurnetrspp2p} can be rewritten as
\begin{EQA}[c]
	\bigl\| \DF^{-1} \IFT (\upsvd - \upsvs + \DF^{-2} \AAv) - \DeltaB \, \DF^{-1} \AAv \bigr\|_{\infty}
	\leq 
	\dltwu_{\CONSTzz} \, \zzinf^{2}
	\, .
\end{EQA}
With \( \BBF \eqdef \DF^{-1} \IFT \, \DF^{-1} \) and 
\( \epsv \eqdef \DF (\upsvd - \upsvs + \DF^{-2} \AAv) \), this yields
\begin{EQA}[c]
	\bigl\| \BBF \, \epsv \bigr\|_{\infty}
	\leq 
	\| \DeltaB \, \DF^{-1} \AAv \|_{\infty}	+ \dltwu_{\CONSTzz} \, \zzinf^{2}
	\, .
\end{EQA}
As \( \| \epsv \|_{\infty} = \| \BBF^{-1} \BBF \epsv \|_{\infty} \leq 
\| \BBF^{-1} \|_{1} \; \| \BBF \epsv \|_{\infty} = \crossB^{-1} \| \BBF \epsv \|_{\infty} \),
it follows
\begin{EQA}[c]
	\crossB \| \epsv \|_{\infty} 
	\leq 
	\| \DeltaB \, \DF^{-1} \AAv \|_{\infty} + \dltwu_{\CONSTzz} \, \zzinf^{2}
	\, .
\end{EQA}
This proves \eqref{usdhyw6hikhurnetrsppD}.

\paragraph{Proof of Theorem~\ref{Psemibiassupp}}
The idea is to reduce the case of a separable perturbation to the case of a linear perturbation
of the function \( \fd(\upsv) = \fs(\upsv) + \sum_{j=1}^{\dimp} \pert_{j}(\ups_{j}) \).
Namely, \( \upsvs = \arginf_{\upsv} \fs(\upsv) \) yields \( \nabla \fs(\upsvs) = 0 \) and hence,
\( \nabla \fd(\upsvs) = \AAv = \sum_{j=1}^{\dimp} \pert'_{j}(\upss_{j}) \).
Consider the function \( g(\upsv) = \fd(\upsv) - \langle \AAv,\upsv \rangle \)
which is a linear perturbation of \( \fd(\upsv) \).
This function \( g(\upsv) \) is strongly convex and \( \nabla g(\upsvs) = 0 \).
Therefore, \( \upsvs = \arginf_{\upsv} g(\upsv) \) and the statements of the theorem follows
from Theorem~\ref{PsemibiassupC} and Theorem~\ref{Psemibiassup} applied to the function \( \fd(\upsv) \)
and a linear perturbation with the slope vector \( -\AAv \).

\paragraph{A bound \( \crossB = \| \DeltaB \|_{1} \geq 1 - \crosssup \)}
\begin{lemma}
\label{LinvGersh}
With \( \BBF = \DF^{-1} \IFT \DF^{-1} \), let \( \DeltaB = \BBF - \Id_{\dimp} \) satisfy 
\begin{EQA}[c]
	\| \DeltaB \|_{1} \leq \crosssup < 1 .
\label{8f98dej237yt656543n}
\end{EQA}
Then 
\begin{EQA}[c]
	\| \BBF^{-1} \|_{1} 
	=
	\sup_{\| \uv \|_{\infty} \leq 1} \| \BBF^{-1} \uv \|_{\infty}
	\leq 
	(1 - \crosssup)^{-1} 
	\, .
\end{EQA}
\end{lemma}

\begin{proof}
The use of \( \BBF^{-1} = (\Id_{\dimp} + \DeltaB)^{-1} 
= \Id_{\dimp} - \DeltaB + \DeltaB^{2} - \ldots \) yields for any \( \uv \in \R^{\dimp} \)
\begin{EQA}[c]
	\| \BBF^{-1} \uv \|_{\infty}
	\leq 
	\sum_{m=0}^{\infty} \| \DeltaB^{m} \uv \|_{\infty} \, .
\end{EQA}
Further, with \( \uv_{m} \eqdef \DeltaB^{m} \uv \), it holds
\begin{EQA}
	\| \DeltaB^{m+1} \uv \|_{\infty}
	&=&
	\| \DeltaB \uv_{m} \|_{\infty}
	\leq 
	\crosssup \| \uv_{m} \|_{\infty} \, .
\label{cd7dysrtd4ftgewhged4r}
\end{EQA}
By induction, this yields \( \| \DeltaB^{m} \uv \|_{\infty} \leq \crosssup^{m} \| \uv \|_{\infty} \) and thus,
\begin{EQA}[c]
	\| \BBF^{-1} \uv \|_{\infty}
	\leq 
	\sum_{m=0}^{\infty} \crosssup^{m} \| \uv \|_{\infty}
	= 
	\frac{1}{1 - \crosssup} \, \| \uv \|_{\infty}
\end{EQA}
as claimed. 
\end{proof}

%% file: BTL_short.tex

\def\nbin{n}
\Chapter{Estimation for Bradley-Terry-Luce model}
\label{ScoBTL}

Let \( \Graph = (\GrV, \GrE) \) stand for a {comparison graph}, where the
vertex set \( \GrV = \{ 1,2,\ldots,\dimGrV \} \) represents the \( \dimGrV \) items of interest. 
The items \( j \) and \( m \) are {compared} if and only if \( (j,m) \in \GrE \).
One observes independent paired comparisons \( Y^{(\ell)}_{jm} \), \( \ell = 1,\ldots,\nbin_{jm} \),
and \( Y^{(\ell)}_{jm} = 1 - Y^{(\ell)}_{mj} \).
BTL model from \cite{BT1952}, \cite{Lu1959} suggests that the chance of each item winning a paired comparison 
is determined by the {relative scores} 
\begin{EQA}
	\P\bigl( \text{item \( j \) is preferred over item } m \bigr)
	&=&
	\P\bigl( Y_{jm}^{(\ell)} = 1 \bigr)
	=
	\frac{\ex^{\upss_{j}}}{\ex^{\upss_{j}} + \ex^{\upss_{m}}} 
	=
	\frac{1}{1 + \ex^{\upss_{m} - \upss_{j}}} \, .
\label{Pjoieujeui}
\end{EQA}
Estimation and inference are reduced to recovering the {score vector} 
\( \upsvs = (\upss_{1},\ldots,\upss_{\dimGrV})^{\T} \).
We assume that \( \Graph \) is connected; otherwise, 
each connected component should be considered separately. 
%
For any pair \( j < m \),
denote \( S_{jm} = \sum_{\ell=1}^{\nbin_{jm}} Y_{jm}^{(\ell)} \) and \( S_{jm} = 0 \) if \( \nbin_{jm} = 0 \).
With \( \cdens(\ups) = \log(1 + \ex^{\ups}) \), the negative log-likelihood function reads as follows:
\begin{EQA}
	L(\upsv)
	&=&
	- \sum_{j=1}^{\dimGrV} \sum_{m = j+1}^{\dimGrV} 
	\bigl\{ (\ups_{j} - \ups_{m}) S_{jm} - \nbin_{jm} \cdens(\ups_{j} - \ups_{m}) \bigr\} 
	\, ,
	\qquad
\label{LusijEuiuj}
\end{EQA}
leading to the MLE \( \tilde{\upsv} = \arginf_{\upsv} L(\upsv) \).
The stochastic component \( \zeta(\upsv) = L(\upsv) - \E L(\upsv) \) is linear in \( \upsv \) and satisfies
\begin{EQA}[c]
	\nabla \zeta(\upsv) \equiv \AAv = (\AAv_{j}) \) with 
	\( \AAv_{j} = - \sum_{m \neq j} (S_{jm} - \E S_{jm}) 
	\, .
\label{jvjjweweseivgnjjjhr}
\end{EQA}
The function \( \cdens(\ups) = \log(1 + \ex^{\ups}) \) is convex, hence, \( L(\upsv) \) is {convex} too.
%
%
The corresponding Fisher information matrix \( \IF(\upsv) = \nabla^{2} \E L(\upsv) \)
has a special stochastically dominating structure.
%
Namely, with \( \cdens''(\ups) = \frac{\ex^{\ups}}{(1 + \ex^{\ups})^{2}} \), 
its entries \( \IF_{jm}(\upsv) \) satisfy
\begin{EQ}[rcl]
\label{sfhce8fy8efhidveuuedf8}
	\IF_{jm}(\upsv)
	&=&
	- \nbin_{jm} \cdens''(\ups_{j} - \ups_{m}) ,
	\quad
	j \neq m 
	\, ,
	\\
	\IF_{jj}(\upsv)
	&=&
	\sum_{m \neq j} \nbin_{jm} \cdens''(\ups_{j} - \ups_{m})
	=
	- \sum_{m \neq j} \IF_{jm}(\upsv)
	\, .
\end{EQ}
This property is important for proving a concentration of the MLE 
\( \tilde{\upsv} = \arginf_{\upsv} L(\upsv) \).
%
However, representation \eqref{LusijEuiuj} reveals \emph{lack-of-identifiability} problem: 
\( \tilde{\upsv} \) is {not unique}, any shift 
\( \upsv \to \upsv + a \ev \)
does not change \( L(\upsv) \), \( \ev = \dimGrV^{-1/2} (1,\ldots,1)^{\T} \in \R^{\dimGrV} \).
Therefore, the Fisher information matrix \( \IF(\upsv) = \nabla^{2} L(\upsv) \) is not positive definite and
\( L(\upsv) \) is not strongly convex.
For a connected graph \( \Graph \), assumed later on,
this issue can be resolved by fixing one component of \( \upsv \), e.g. \( \ups_{1} = 0 \), or by the condition 
\( \sum_{j} \ups_{j} = 0 \).
In general, we need {one condition} per {connected component} of the graph \( \Graph \).
Alternatively, one can use a penalized MLE 
with a {quadratic} penalty \( \| \GP \upsv \|^{2}/2 \).
The constraint \( \sum_{j} \ups_{j} = 0 \) can be replaced by the penalty 
\( \| \GP \upsv \|^{2} = \gp^{2} \langle \upsv,\ev \rangle^{2} \)
with \( \GP^{2} = \gp^{2} \ev \, \ev^{\T} \), 
because it vanishes for any \( \upsv \) with \( \sum \ups_{j} = 0 \).
If the true skill vector \( \upsvs \) satisfies \( \sum \upss_{j} = 0 \) then 
\( \upsvs_{\GP} = \arginf_{\upsv} \E \LGP(\upsv) = \upsvs \).
Hence, the penalty \( \| \GP \upsv \|^{2} = \gp^{2} \langle \upsv,\ev \rangle^{2} \) does not yield any bias of estimation
even for \( \gp^{2} \) large.
Another option is to use Tikhonov's regularization \( \GP^{2} = \gp^{2} \Id_{\dimGrV} \); cf. \cite{CFMW2019}.
Even for a disconnected graph \( \Graph \), it ensures the desired stochastically dominant structure 
of the Fisher information matrix: each diagonal element is larger than the sum of the off-diagonal elements for this row. 
This penalty leads to some bias, namely, to a shrinking effect: the estimated value \( \upsvs_{\GP} \) is shrunk towards zero.
Later we fix \( \GP^{2} = \gp^{2} \Id_{\dimGrV} \) and consider the regularized MLE \( \tilde{\upsv} \) 
\begin{EQA}[c]
		\tilde{\upsv}_{\GP}
	=
	\arginf_{\upsv} \LGP(\upsv) 
	=
	\arginf_{\upsv} \bigl( L(\upsv) + \gp^{2} \| \upsv \|^{2}/2 \bigr) .
\end{EQA}
Define \( \DF^{2} = \diag(\DPN_{1}^{2},\ldots,\DPN_{\dimGrV}^{2}) \) with \( \DPN_{j}^{2} = \IF_{jj} + \gp^{2} \) for 
\( \IF = \IF(\upsvs_{\GP}) \), and
\begin{EQA}[c]
	\IFmax
	\eqdef
	\max_{j \leq \dimp} \IF_{jj}
	\, ,
	\qquad
	\IFmin
	\eqdef
	\min_{j \leq \dimp} \IF_{jj}
	\, ,
	\qquad
	\IFrange
	\eqdef
	{\IFmin}/{\IFmax}
	\, .
\end{EQA}
Also, define \( \BBF \eqdef \DF^{-1} \IF_{\GP} \, \DF^{-1} \),
\begin{EQ}[rclrcl]
	\crosssup
	&=&
	\max_{j = 1,\ldots,\dimGrV} \biggl( \frac{1}{\DPN_{j}} \sum_{m \neq j} \frac{|\IF_{jm}|}{\DPN_{m}} \biggr)
	\, ,
	\quad
	&
	\crosssupp
	& \eqdef &
	\max_{j = 1,\ldots,\dimGrV} \biggl( \sum_{m \neq j} \frac{|\IF_{jm}|}{\DPN_{m}^{2}} \biggr)
	\, ,
	\\
	\crossD
	& \eqdef &
	\max_{j = 1,\ldots,\dimGrV} \biggl( \frac{1}{\DPN_{j}^{2}} \sum_{m \neq j} \frac{\IF_{jm}^{2}}{\DPN_{m}^{2}} \biggr)^{1/2}
	\, ,
	\quad
	&
	\crossB^{-1} 
	&=&
	\sup_{\| \uv \|_{\infty} \leq 1} \| \BBF^{-1} \uv \|_{\infty}
	=
	\| \BBF^{-1} \|_{1}
	\, .
\label{7tdsyf8iuwopkrtgBTL}
\end{EQ}
Lemma~\ref{Lcrosssup} provides upper bounds on these values in a balanced case when
all the \( \DPN_{j} \)'s are of the same order and the values \( \nbin_{jm} \) 
(numbers of games between \( j \) and \( m \)) satisfy
\( \nbin_{jm} \leq \nbin_{\max} \).
Namely, 
\( \crosssup < 1 \), \( \crosssupp \) is of order 1, \( \crossB \geq 1 - \crossnorm \), while
\( \crossD \asymp (\neff/\nbin_{\max})^{-1/2} \) for \( \neff \eqdef \IFmin + \gp^{2} \).

\begin{theorem}
\label{TBTLall}
Let \( Y_{jm}^{(\ell)} \) follow the BTL model with the true vector \( \upsvs \).
Let \( \IFmin \eqdef \min_{j \leq \dimGrV} \IF_{jj} \geq \CONSTi_{2} \{ \log(\dimGrV) + \xx \} \) 
for some fixed constant \( \CONSTi_{2} \).
Fix \( \GP^{2} = \gp^{2} \Id_{\dimp} \) with \( \gp^{2} = \IFmax \). 
Then on a random set \( \Omega(\xx) \) with \( \P\bigl( \Omega(\xx) \bigr) \geq 1 - 4 \ex^{-\xx} \), it holds
with \( \AAv = \nabla \zeta \) from \eqref{jvjjweweseivgnjjjhr} and \( \DeltaB = \DF^{-1} \IF_{\GP} \, \DF^{-1} - \Id_{\dimGrV} \)
\begin{EQA}[rcccl]
	\| \DF^{-1} \AAv \|_{\infty} 
	& \leq &
	\zzinf
	& \eqdef &
	\sqrt{\log(\dimGrV) + \xx}
	\, ,
	\\
	\| \DeltaB \, \DF^{-1} \AAv \|_{\infty} 
	& \leq &
	\crossD \, \zzinf
	&=&
	\crossD \sqrt{\log(\dimGrV) + \xx} 	
	\, ,
	\qquad
\end{EQA}
Also, 
on \( \Omega(\xx) \), it holds 
with \( \dltwu_{\CONSTzz} \) from \eqref{usdhyw6hikhurnetrspp2},  
\begin{EQ}[rcl]
	\| \DF (\tilde{\upsv}_{\GP} - \upsvs_{\GP}) \|_{\infty}
	& \leq &
	\frac{3 \zzinf}{1 - \crosssup}
	\, ,
	\\
	\bigl\| \DF^{-1} \bigl\{ \IFN_{\GP} (\tilde{\upsv}_{\GP} - \upsvs_{\GP}) + \AAv \bigr\} 
	\bigr\|_{\infty}
	& \leq &
	\dltwu_{\CONSTzz} \, \| \DF^{-1} \AAv \|_{\infty}^{2}
	\leq 
	\dltwu_{\CONSTzz} \, \zzinf^{2}
	\, ,
	\\
	\| \DF (\tilde{\upsv}_{\GP} - \upsvs_{\GP}) + \DF^{-1} \AAv \|_{\infty} 
	& \leq &
	\crossB^{-1} \bigl( 
		\| \DeltaB \, \DF^{-1} \AAv \|_{\infty} + \dltwu_{\CONSTzz} \| \DF^{-1} \AAv \|_{\infty}^{2} 
	\bigr)
	\\
	& \leq &
	\crossB^{-1} \bigl( \crossD \, \zzinf + \dltwu_{\CONSTzz} \, \zzinf^{2} \bigr)
	\, .
	\qquad
\label{usdhyw6hikhurnetrspBTL}
\end{EQ}
\end{theorem}

\paragraph{Critical sparseness}
\cite{GSZ2023} highlighted importance of considering a sparse graph regime with \( \alpG \asymp \dimGrV^{-1} \bigl( \log \dimGrV \bigr)^{a} \) and stated the results for \( a = 3/2 \).
Result \eqref{usdhyw6hikhurnetrspBTL} can be applied to the extremely sparse regime with 
\( \alpG \asymp \dimGrV^{-1} \log \dimGrV  \).
The only requirement is 
that its right-hand side is smaller in order than the leading term
\( \| \DF^{-1} \AAv \|_{\infty} \approx \zzinf \) of the left-hand side.
This leads to \( \crossD \ll 1 \) and \( \dltwu_{\CONSTzz} \, \zzinf \ll 1 \).
In view of \( \dltwu_{\CONSTzz} \asymp \IFmin^{-1/2} \), these condition meet as soon as \( \IFmin \gg \log \dimGrV \),
which improves \cite{GSZ2023}.

\paragraph{Bias issue}
Due to \eqref{usdhyw6hikhurnetrspBTL}, the penalized MLE \( \tilde{\upsv}_{\GP} \) estimates
\( \upsvs_{\GP} \) rather then \( \upsvs \).
The bias \( \upsvs_{\GP} - \upsvs \) induced by penalization \( \gp^{2} \| \upsv \|^{2}/2 \) is not negligible. 
Theorem~\ref{Psemibiassupp} can be applied to evaluate the leading term of this bias:
\( \upsvs_{\GP} - \upsvs \approx - \gp^{2} \IF_{\GP}^{-1} \upsvs = - \gp^{2} (\IF + \gp^{2} \Id_{\dimp})^{-1} \upsvs \);
cf. e.g. \cite{CFMW2019}.
One can treat \( \upsvs_{\GP} \) as a shrinkage of \( \upsvs \), and the ranking results based on 
\( \upsvs_{\GP} \) are very similar to the true ranking based on \( \upsvs \).

\begin{figure}[t]
\centering
\includegraphics[width=0.6\linewidth,height=0.25\textheight]{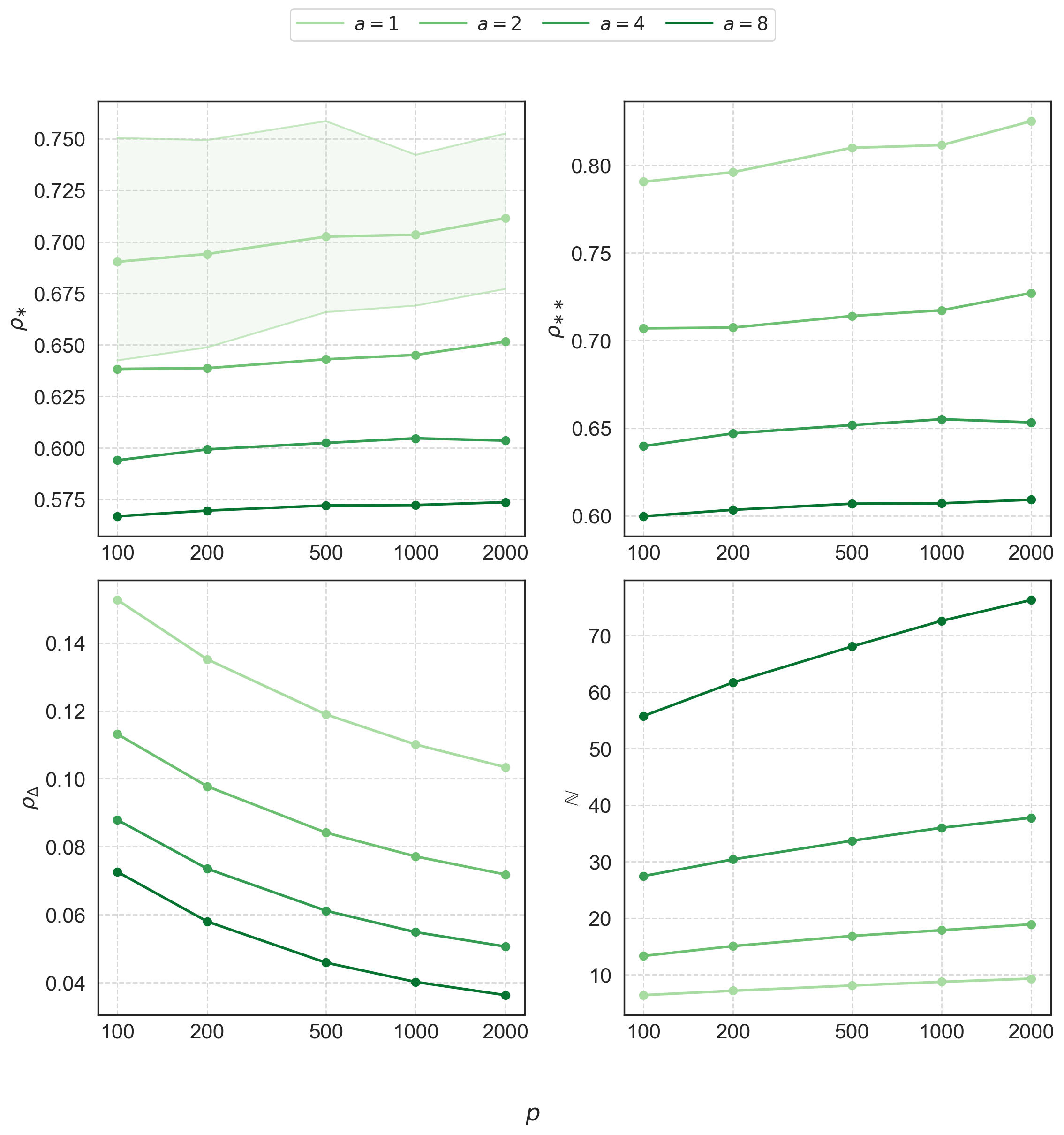}
\caption{\( \crosssup \) (top left), \( \crosssupp \) (top right), \( \crossD \) (bottom left), and \( \neff \) (bottom right) for different choices of \( \dimGrV \) and \( a \). }
\label{btl_quantities}
\end{figure}

\paragraph{Numerical results}
We use the following experimental setup, which extends \cite{GSZ2023}.
In particular, it allows an unbounded range of the \( \upss_{j} \)'s, more realistic design 
of the graph comparison, and, most important, a much weaker condition of graph sparsity.
The values \( \upss_{i} \) are i.i.d  samples from  \( \text{Logistic} (0, 1/4 ) \), 
the comparison graph is random with \( \dimGrV \) vertices and probability of the edge \( (jm) \) given by
\begin{EQA}[c]
	\alpG_{jm} = \frac{8 a \log \dimGrV \, \cdens''(\upss_{j} - \upss_{m})}{\dimGrV \, (\Phi_{j} + \Phi_{m}) }
	\, , 
	\qquad
	\Phi_{j}
	\eqdef
	\sum_{m \neq j} \cdens''(\upss_{j} - \upss_{m})
	\, ,
\end{EQA}
where \( a \in \{ 1,2,4,8 \} \).
This rule benefits balanced pairs with nearly the same skills at the cost of unbalanced games between players with very different skills.
The number of comparisons between pairs of items is \( L \equiv 1 \).
The regularization parameter is \( \gp^{2} = a \, \log \dimp \).
For experiments, we used 
\( \dimGrV \in \{100, 200, 500, 1000, 2000\} \). 
The study presents the values 
{ 
\( \crosssup \) (top left), \( \crosssupp \), \( \crossD \) from \eqref{7tdsyf8iuwopkrtgBTL} (top),
and \( \neff = \IFmin + \gp^{2} \)
}; see Figure~\ref{btl_quantities}.
Even in the sparse regime, the required conditions on these quantities are satisfied, i.e. 
\( \crosssup \) is bounded from 1, \( \crosssupp \) is of order 1, \( \crossD \ll 1 \), and \( \neff \) grows with \( \dimGrV \).

\begin{figure}[t]
\centering
\includegraphics[width=0.7\linewidth,height=0.2\textheight]{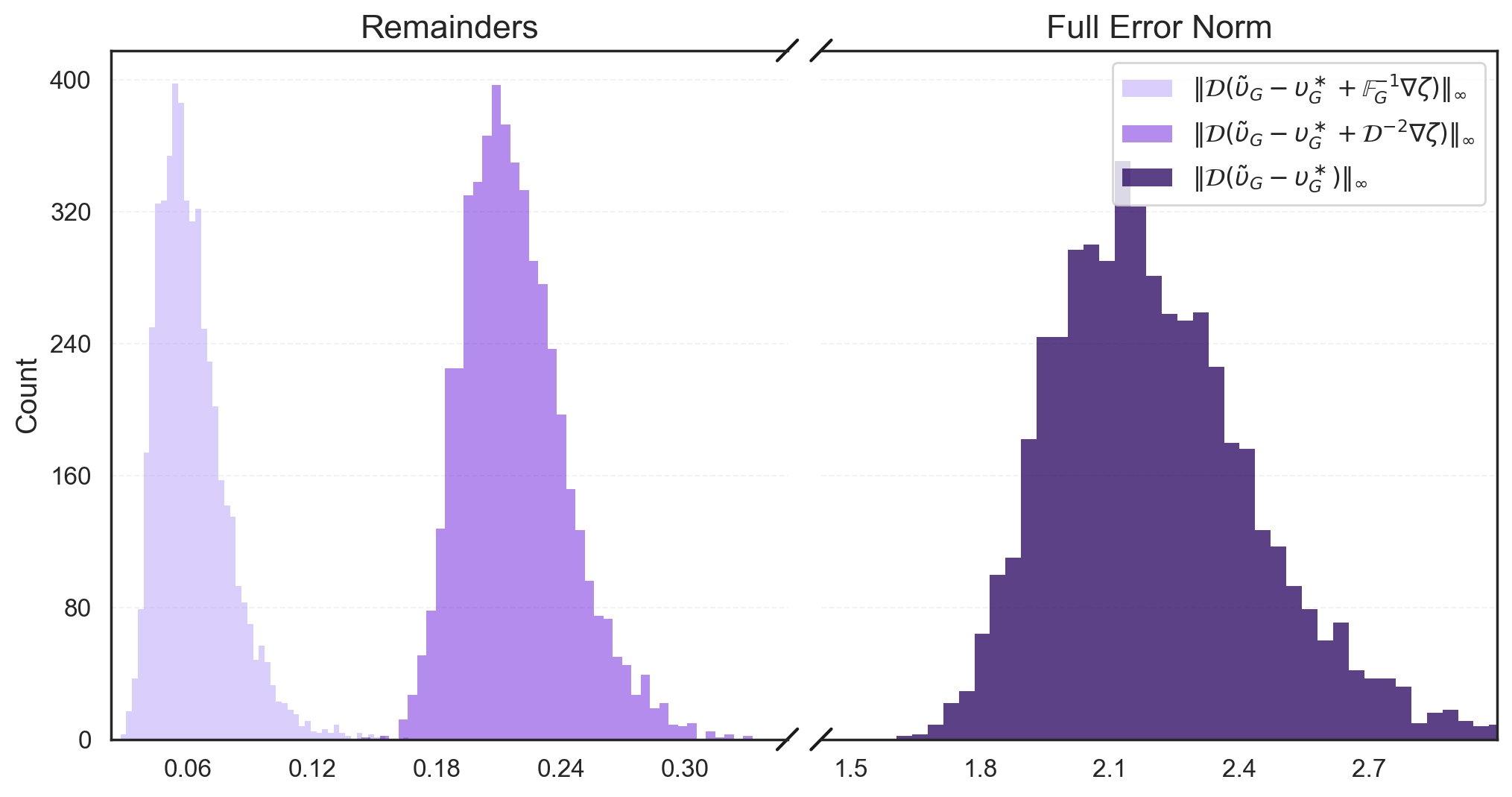}
\caption{The results  for \( \dimGrV = 500 \), \( a=2 \).
Left: approximation errors \( \| \DF (\tilde{\upsv}_{\GP} - \upsvs_{\GP} + \IF_{\GP}^{-1} \nabla \zeta) \|_{\infty} \).
Middle: approximation errors \( \| \DF( \tilde{\upsv}_{\GP} - \upsvs_{\GP} + \DF^{-2} \nabla \zeta) \|_{\infty} \). 
Right: the norm of the estimation error \( \| \DF (\tilde{\upsv}_{\GP} - \upsvs_{\GP}) \|_{\infty} \).
}
\label{fig:BTLapprox_p_500}
\end{figure}

\begin{figure}[t]
\centering
\includegraphics[width=0.6\linewidth,height=0.2\textheight]{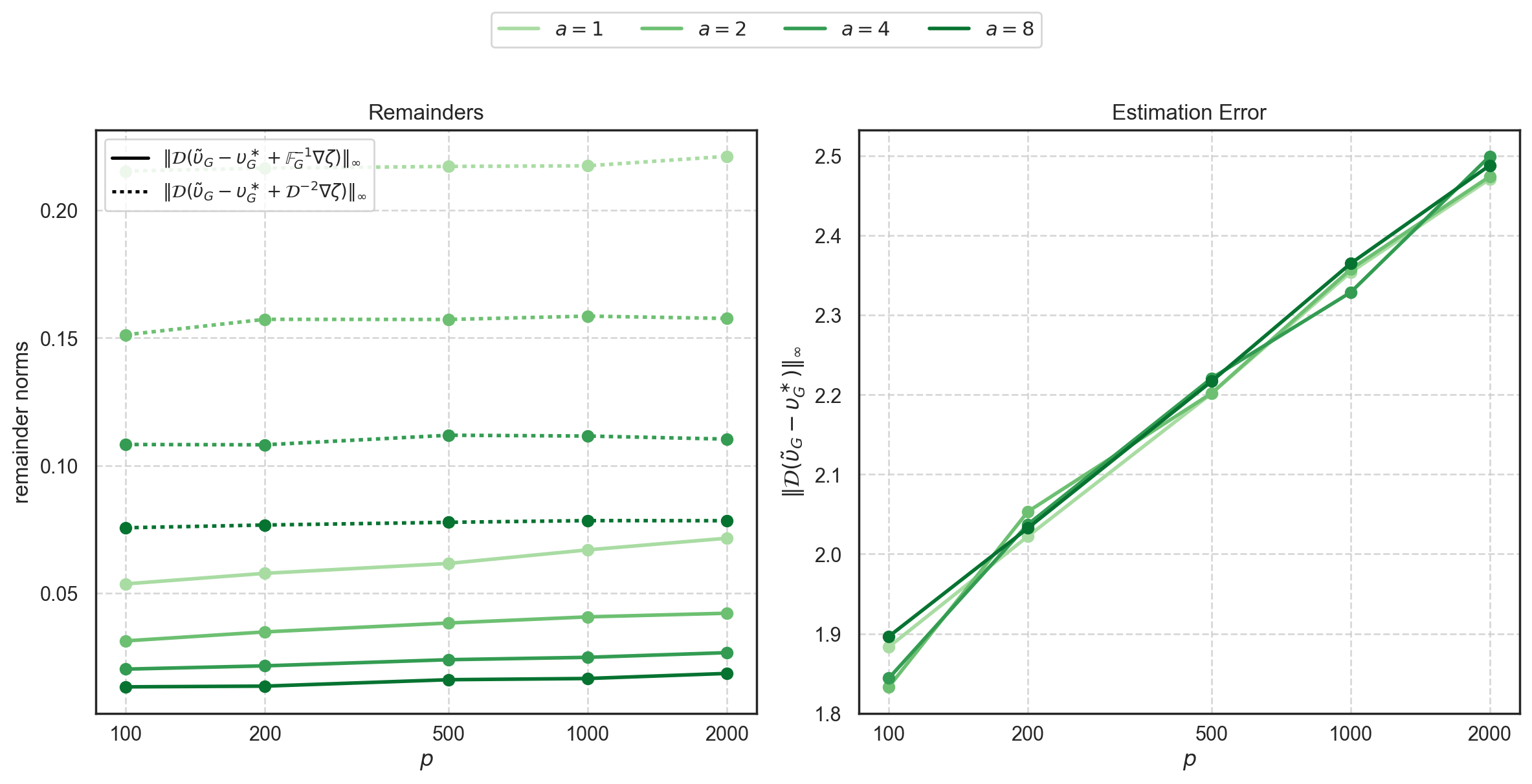}
\caption{Same as Figure~\ref{fig:BTLapprox_p_500}, but estimated for different choices of \( \dimGrV \) and \( a \). 
}
\label{fig:BTLapprox_all_p}
\end{figure}

Figure~\ref{fig:BTLapprox_p_500} demonstrates the quality of 
expansion \eqref{usdhyw6hikhurnetrspBTL} for \( \dimGrV = 500 \) and \( a=2 \).
We can see that the remainders (mean = 0.07 for the Fisher correction \( \IF_{\GP}^{-1} \nabla \zeta \), or mean = 0.2
for the coordinatewise correction \( \DF^{-2} \nabla \zeta \)) are significantly smaller than the full error (mean = 2.21), confirming the quality of the expansion.
Figure~\ref{fig:BTLapprox_all_p} summarizes the results for different choices of \( \dimGrV \) and \( a \).

%% file: BTL_short_proofs.tex
\paragraph{Proof of Theorem~\ref{TBTLall}}
Here we sketch the proof of the theorem.
The detailed derivations and more results can be found in the forthcoming paper \cite{PankSp2025}.
We are going to apply Theorem~\ref{Psemibiassup} to the penalized MLE 
\begin{EQA}
	\tilde{\upsv}_{\GP}
	&=&
	\arginf_{\upsv} \LGP(\upsv) 
	=
	\arginf_{\upsv} \bigl( L(\upsv) + \| \GP \upsv \|^{2}/2 \bigr) 
	\, ,
	\qquad
	\| \GP \upsv \|^{2} = \gp^{2} \| \upsv \|^{2}
	\, ,
\label{ftbc43e3e3rfdf6ddfd}
\end{EQA}
with \( L(\upsv) \) from \eqref{LusijEuiuj}, 
and its population counterpart 
\begin{EQA}
	\upsvs_{\GP}
	&=&
	\arginf_{\upsv} \E \LGP(\upsv) 
	=
	\arginf_{\upsv} \bigl( \E L(\upsv) + \| \GP \upsv \|^{2}/2 \bigr) 
	\, .
\label{ftbc43e3e3rfdf6ddfd}
\end{EQA}
Fix \( \IF = \IF(\upsvs_{\GP}) \) and \( \DF^{2} = \diag(\DPN_{1}^{2},\ldots,\DPN_{\dimGrV}^{2}) \) with \( \DPN_{j}^{2} = \IF_{jj} + \gp^{2} \).
Under the deterministic design \( (\Psiv_{i}) \), the stochastic component \( \zeta(\upsv) \eqdef \LL(\upsv) - \E \LL(\upsv) \)
is linear in \( \upsv \) and it holds \( \nabla \zeta = (\nabla \zeta_{j}) \) with
\begin{EQA}[c]
	\nabla \zeta_{j} = - \sum_{m \neq j} (S_{jm} - \E S_{jm})
	\, ,
	\qquad
	S_{jm} = \sum_{\ell=1}^{\nbin_{jm}} Y_{jm}^{(\ell)}
	\, .
\end{EQA}
Conditions of Theorem~\ref{Psemibiassup} require bounding in probability the sup-norm of the vectors \( \DF^{-1} \nabla \zeta \) 
and \( \DeltaB \, \DF^{-1} \nabla \zeta \), and
evaluating \( \crosssup \) from \eqref{7tdsyf8iuwopkrtg4576} and checking
\nameref{LLpsupref}.

First, we provide deviation bounds for the norms \( \| \DF^{-1} \nabla \zeta \|_{\infty} \) and 
\( \| \DeltaB \, \DF^{-1} \nabla \zeta \|_{\infty} \).

\begin{lemma}
\label{LVarBTL}
It holds with \( \nabla \zeta_{j} = - \sum_{m \neq j} (S_{jm} - \E S_{jm}) \)
\begin{EQA}[c]
	\Var (\nabla \zeta_{j})
	=
	\IF_{jj}
	\leq 
	\DPN_{j}^{2}
	\, ,
	\qquad
	\P\bigl( \| \DF^{-1} \nabla \zeta \|_{\infty} > 2 \sqrt{\log(\dimGrV) + \xx} \bigr)
	\leq 
	2 \ex^{-\xx}
	\, .
\label{cvjuvcu843jrtgytrud}
\end{EQA}
Similarly, with \( \DeltaB = \DF^{-1} \IFT \DF^{-1} - \Id_{\dimGrV} \), it holds
\( (\DeltaB \, \DF^{-1} \nabla \zeta)_{j} = \sum_{m \neq j} \DeltaB_{jm} \DPN_{mm}^{-1} \nabla \zeta_{m} \) and
\begin{EQA}[c]
	\Var\bigl\{ (\DeltaB \, \DF^{-1} \nabla \zeta)_{j} \bigr\}
	=
	\sum_{m \neq j} \Var(\DeltaB_{jm} \DPN_{mm}^{-1} \nabla \zeta_{m})
	\leq 
	\sum_{m \neq j} \DeltaB_{jm}^{2}
	=
	\sum_{m \neq j} \frac{\IF_{jm}^{2}}{\DPN_{jj}^{2} \, \DPN_{mm}^{2}}
	\leq 
	\crossD
	\, ,
	\qquad
\label{fdjihui4rboe73urtgg}
\end{EQA}
and
\begin{EQA}[c]
	\P\bigl( \| \DeltaB \, \DF^{-1} \nabla \zeta \|_{\infty} > 2 \crossD \sqrt{\log(\dimGrV) + \xx} \bigr)
	\leq 
	2 \ex^{-\xx}	
	\, .
\label{ifgivgkhjivtwjifg8gh}
\end{EQA}
\end{lemma}

\begin{proof}
For every \( j \leq \dimGrV \), applying the deviation bound from \cite{Sp2023d}, Proposition 3.1, for a weighted sum of Bernoulli r.v.'s yields
\begin{EQA}[c]
	\P\bigl( \bigl| \DPN_{j}^{-1} \nabla \zeta_{j} \bigr| > 2 \sqrt{\log(\dimGrV) + \xx} \bigr)
	\leq 
	2 \ex^{- \log(\dimGrV) - \xx}
	\, .
\end{EQA}
Therefore, 
\begin{EQA}[c]
	\P\bigl( \| \DF^{-1} \nabla \zeta \|_{\infty} > 2 \sqrt{\log(\dimGrV) + \xx} \bigr)
	\leq 
	\sum_{j=1}^{\dimGrV} 2 \ex^{- \log(\dimGrV) - \xx}
	=
	\ex^{-\xx}
	\, .
\end{EQA}
This proves \eqref{cvjuvcu843jrtgytrud}. 
Bound \eqref{ifgivgkhjivtwjifg8gh} can be checked similarly; it suffices to verify \eqref{fdjihui4rboe73urtgg}.
\end{proof}

The next important step is to check that \( \crosssup < 1 \) and evaluate \( \crosssupp \) and \( \crossD \); see \eqref{7tdsyf8iuwopkrtgBTL}.

\begin{lemma}
\label{Lcrosssup}
Denote \( \IFmax \eqdef \max_{j} \IFN_{jj} \) and \( \IFmin \eqdef \min_{j} \IFN_{jj} \).
Then
\begin{EQA}[rcccl]
\label{diciiiivfy6f76rf76rf7}
	\frac{\IFmax}{\IFmax + \gp^{2}}
	& \leq &
	\crosssup
	& \leq &
	\frac{\IFmax}{\sqrt{(\IFmin + \gp^{2}) (\IFmax + \gp^{2})}}
	\, ,
	\\
	\frac{\IFmax}{\IFmax + \gp^{2}}
	& \leq &
	\crosssupp
	& \leq &
	\frac{\IFmax}{\IFmin + \gp^{2}}
	\, .
\end{EQA}
With \( \IFrange \eqdef \IFmin/\IFmax \) and \( \gp^{2} = \IFmax \), it holds
\begin{EQA}[c]
	\frac{1}{2}
	\leq
	\crosssup
	\leq 
	\frac{1}{\sqrt{2(\IFrange + 1)}}
	\, ,
	\qquad
	\frac{1}{2}
	\leq
	\crosssupp
	\leq 
	\frac{1}{{\IFrange + 1}}
	\, .
\end{EQA}
Also,
\begin{EQA}[c]
	\crossD^{2}
	\leq 
	\frac{\crosssupt}{\IFmin + \gp^{2}}
	\, ,
	\qquad
	\crosssupt
	\eqdef
	\max_{j \leq  \dimp} \frac{1}{\DPN_{j}^{2}} \sum_{m \neq j} \IF_{jm}^{2}
	\leq 
	\max_{m \neq j} \IF_{jm}
	\, .
\end{EQA}
\end{lemma}

\begin{proof}
As \( \DPN_{m}^{2} = \IF_{mm} + \gp^{2} \leq \IFmax + \gp^{2} \) and \( \sum_{m \neq j} |\IF_{jm}| = \IF_{jj} \), it holds for \( j \leq \dimGrV \),
\begin{EQA}[c]
	\frac{1}{\DPN_{j}} \sum_{m \neq j} \frac{|\IF_{jm}|}{\DPN_{m}}
	\geq 
	\frac{1}{\DPN_{j} \sqrt{\IFmax + \gp^{2}}} \sum_{m \neq j} |\IF_{jm}|
	=
	\frac{\IF_{jj}}{\sqrt{(\IF_{jj} + \gp^{2}) (\IFmax + \gp^{2})}} 
\end{EQA}
yielding
\begin{EQA}[c]
	\crosssup
	=
	\max_{j \leq \dimp} \frac{1}{\DPN_{j}} \sum_{m \neq j} \frac{|\IF_{jm}|}{\DPN_{m}}
	\geq 
	\max_{j \leq \dimp} \frac{\IF_{jj}}{\sqrt{(\IF_{jj} + \gp^{2}) (\IFmax + \gp^{2})}}
	=
	\frac{\IFmax}{\IFmax + \gp^{2}}
	\, .
\end{EQA}
Similarly,
\begin{EQA}[c]
	\frac{1}{\DPN_{j}} \sum_{m \neq j} \frac{|\IF_{jm}|}{\DPN_{m}}
	\leq 
	\frac{\IF_{jj}}{\sqrt{(\IF_{jj} + \gp^{2}) (\IFmin + \gp^{2})}} 
	\leq 
	\frac{\IFmax}{\sqrt{(\IFmin + \gp^{2}) (\IFmax + \gp^{2})}} 
\end{EQA}
and \eqref{diciiiivfy6f76rf76rf7} follows.
Further, for every \( j \leq \dimp \),
\begin{EQA}[c]
	\sum_{m \neq j} \frac{|\IF_{jm}|}{\DPN_{m}^{2}}
	\leq 
	\frac{1}{\IFmin + \gp^{2}} \sum_{m \neq j} |\IF_{jm}|
	\leq 
	\frac{\IF_{jj}}{\IFmin + \gp^{2}}
	\leq 
	\frac{\IFmax}{\IFmin + \gp^{2}}
	\, ,
\end{EQA}
and similarly
\begin{EQA}[c]
	\max_{j \leq \dimp} \sum_{m \neq j} \frac{|\IF_{jm}|}{\DPN_{m}^{2}}
	\geq 
	\max_{j \leq \dimp} \frac{1}{\IFmax + \gp^{2}} \sum_{m \neq j} |\IF_{jm}|
	\geq 
	\max_{j \leq \dimp} \frac{\IF_{jj}}{\IFmax + \gp^{2}} 
	= 
	\frac{\IFmax}{\IFmax + \gp^{2}}
	\, .
\end{EQA}
Finally, for every \( j \leq \dimp \), by \( \IF \)
\begin{EQA}[c]
	\frac{1}{\DPN_{j}^{2}} \sum_{m \neq j} \frac{\IF_{jm}^{2}}{\DPN_{m}^{2}}
	\leq 
	\frac{1}{(\IFmin + \gp^{2}) \DPN_{j}^{2}} \sum_{m \neq j} \IF_{jm}^{2}
	\leq 
	\frac{\crosssupt}{\IFmin + \gp^{2}}
	\, .
\end{EQA}
This completes the proof.
\end{proof}


\begin{lemma}
For the values \( \crosssup \), \( \crosssupp \) from \eqref{7tdsyf8iuwopkrtgBTL},
and \( \neff \eqdef \IFmin + \gp^{2} \),
condition \nameref{LLpsupref} is satisfied with \( \dltwu_{3} = \ex \neff^{-1/2} \) and 
\( \dltwun = \sqrt{\ex} \, \crosssup \, \neff^{-1/2} \),
\( \dltwunn = \sqrt{\ex} \, \crosssupp \, \neff^{-1/2} \).
\end{lemma}

\begin{proof}
Fix \( \upsv = (\ups_{j}) \) such that \( |\ups_{j} - \upss_{j}| \leq 1/2 \)
for  \( j \leq \dimGrV \).
It follows by \( |\cdens'''(\ups)| \leq \cdens''(\ups) \) for any \( j,m \leq \dimGrV \)
\begin{EQA}[c]
	\bigl| \cdens'''(\ups_{j} - \ups_{m}) \bigr|
	\leq 
	\cdens''(\ups_{j} - \ups_{m})
	\leq 
	\ex \, \cdens''(\upss_{j} - \upss_{m}) 
	\, .
\end{EQA}
Fix \( j \leq \dimGrV \), represent \( \upsv = (\ups_{j},\nupv_{j}) \),
and consider a vector \( \zv = (z_{m}) \in \R^{\dimGrV-1} \) with 
\( \| \HPN \zv \|_{\infty} = \max_{m} |\DPN_{mm} \, z_{m}| = 1 \).
Then 
\begin{EQA}[rcl]
	\bigl| \nabla_{\ups_{j}\ups_{j} \ups_{j}}^{(3)} L(\ups_{j},\nuiv_{j}) \bigr|
	& \leq &
	\sum_{m \neq j} \nbin_{jm} \bigl| \cdens'''(\ups_{j} - \ups_{m}) \bigr| 
	\leq 
	\ex \sum_{m \neq j} \nbin_{jm} \cdens''(\upss_{j} - \upss_{m})
	=
	\ex \, \IFN_{jj}
\end{EQA}
and
\begin{EQA}[rcl]
	&&\nquad
	\bigl| \langle \nabla_{\ups_{j}\ups_{j}\nuiv_{j}}^{(3)} L(\upss_{j},\nuiv_{j}), \zv_{j} \rangle \bigr|
	\\
	& \leq &
	\sum_{m \neq j} \nbin_{jm} \bigl| \cdens'''(\upss_{j} - \ups_{m}) z_{m} \bigr| 
	\leq 
	\sqrt{\ex} \, \sum_{m \neq j} \frac{\nbin_{jm} \cdens''(\upss_{j} - \upss_{m})}{\DPN_{m}}
	\leq 
	\sqrt{\ex} \, \crosssup \, \DPN_{j} 
	\, .
\end{EQA}
Similarly,
\begin{EQA}[rcl]
	\bigl| \langle \nabla_{\ups_{j}\nuiv_{j}\nuiv_{j}}^{(3)} L(\upss_{j},\nuiv_{j}), \zv_{j}^{\otimes 2} \rangle \bigr|
	& \leq &
	\sqrt{\ex} \, \sum_{m \neq j} \frac{\nbin_{jm} \cdens''(\upss_{j} - \upss_{m})}{\DPN_{m}^{2}}
	\leq 
	\sqrt{\ex} \, \crosssupp 
	\, .
\end{EQA}
This implies \nameref{LLpsupref} with \( \dltwu_{3} = \ex \neff^{-1/2} \) and 
\( \dltwun = \sqrt{\ex} \, \crosssup \, \neff^{-1/2} \),
\( \dltwunn = \sqrt{\ex} \, \crosssupp \, \neff^{-1/2} \).
\end{proof}

We are now prepared to finalize the proof by applying Theorem~\ref{Psemibiassup}.
Conditions \( \dltwunn \, \rrinf \leq 1/4 \) and
\( \dltwun \, \rrinf \leq 1/4 \) are fulfilled if \( \neff \gtrsim \rrinf \gtrsim \sqrt{\log \dimGrV} \).


%% file: localbounds_short.tex

\def\AFN{\mathbbmsl{U}}
\def\Avm{\bb{M}}
\def\upsvp{\upsv^{\circ}}
\def\fp{\fs^{\circ}}
\def\fn{\fp}
\def\upsvn{\upsvp}
\def\IFN{\IF}
\def\IFP{\IF_{\circ}}

\Chapter{Perturbed optimization}
\label{Squadnquad}
This section provides an overview of some results from \cite{Sp2024PO}.
Let \( \fs(\upsv) \) be a smooth strongly convex function, 
\begin{EQA}
	\upsvs
	&=&
	\arginf_{\upsv} \fs(\upsv),
\label{fg5hg3gf98tkj3dciryt}
\end{EQA}
and \( \IFN = \nabla^{2} \fs(\upsvs) \).
We study the question of how the point of minimum and the value of minimum of \( \fs \) change 
under some perturbation of \( \fs \).
More precisely, with \( \fp(\upsv) \) being a perturbed version of \( \fs(\upsv) \),
we consider in parallel two optimization problems 
\begin{EQA}
	\upsvs 
	&=& 
	\arginf_{\upsv} \fs(\upsv),
	\quad
	\upsvp = \arginf_{\upsv} \fp(\upsv) 
	\, .
\label{8cvkfc9fujf6jnmcer4gen}
\end{EQA}
The goal is to describe the discrepancies \( \upsvp - \upsvs \) and \( \fp(\upsvp) - \fp(\upsvs) \) induced by this perturbation.
We also present the results for two specific cases: a linear and a quadratic perturbation.


\Section{Smoothness conditions}
\label{Ssmoothshort}

Assume \( \fs(\upsv) \) to be convex and three or sometimes even four times Gateaux differentiable in \( \upsv \in \Ups \).
For a fixed \( \upsv \in \Ups \), we need such a condition uniformly in 
a local vicinity of a fixed point \( \upsv \in \Ups \).
The notion of locality is given in terms of a metric tensor \( \DFN(\upsv) \in \Matr_{\dimp} \),
where \( \Matr_{\dimp} \) is the set of symmetric positive \( \dimp \)-matrices.
In most of the results later on, one can use \( \DFN(\upsv) = \IFN^{1/2}(\upsv) \), where \( \IFN(\upsv) = \nabla^{2} \fs(\upsv) \).
Introduce the following conditions.

\begin{description}
    \item[\label{LLsT3ref} \( \bb{(\mathcal{T}_{3}^{*})} \)]
    \emph{\( \fs(\upsv + \uv) \) is three times differentiable in \( \uv \) and with \( \DFN \) possibly depending on \( \upsv \)
	}
\begin{EQA}
    \sup_{\uv \colon \| \DFN \uv \| \leq \rr} \,\, \sup_{\zv \in \R^{\dimp}} \,\, 
    \frac{\bigl| \langle \nabla^{3} \fs(\upsv + \uv), \zv^{\otimes 3} \rangle \bigr|}
		 {\| \DFN \zv \|^{3}} 
	& \leq &
	\dltwu_{3} \, .
\label{jcxhydtferyu9j3d6vhew}
\end{EQA}

    \item[\label{LLsT4ref} \( \bb{(\mathcal{T}_{4}^{*})} \)]
    \emph{\( \fs(\upsv + \uv) \) is four times differentiable in \( \uv \) and 
	}
\begin{EQA}
    \sup_{\uv \colon \| \DFN \uv \| \leq \rr} \,\, \sup_{\zv \in \R^{\dimp}} \,\, 
    \frac{\bigl| \langle \nabla^{4} \fs(\upsv + \uv), \zv^{\otimes 4} \rangle \bigr|}
		 {\| \DFN \zv \|^{4}} 
	& \leq &
	\dltwu_{4} \, .
\label{jcxhydtferyu9j3d6vhew4}
\end{EQA}

\end{description}

%
\noindent
By Banach's characterization \cite{Banach1938}, \nameref{LLsT3ref} implies
\begin{EQA}
	\bigl| \langle \nabla^{3} \fs(\upsv + \uv), \zv_{1} \otimes \zv_{2} \otimes \zv_{3} \rangle \bigr|
	& \leq &	 
	\dltwu_{3} \| \DFN \zv_{1} \| \, \| \DFN \zv_{2} \| \, \| \DFN \zv_{3} \| \, 
\label{jbuyfg773jgion94euyyfg}
\end{EQA}
for any \( \uv \) with \( \| \DFN \uv \| \leq \rr \) and all \( \zv_{1} , \zv_{2}, \zv_{3} \in \R^{\dimp} \).
\iffourG{
Similarly under \nameref{LLsT4ref}
\begin{EQA}
	\bigl| \langle \nabla^{4} \fs(\upsv + \uv), \zv_{1} \otimes \zv_{2} \otimes \zv_{3} \otimes \zv_{4} \rangle \bigr|
	& \leq &	 
	\dltwu_{4} \prod_{k=1}^{4} \| \DFN \zv_{k} \| \, ,
	\quad 
	\zv_{1} , \zv_{2}, \zv_{3}, \zv_{4} \in \R^{\dimp} \, .
	\qquad
\label{jbuyfg773jgion94euyyfg4}
\end{EQA}
}{}

\iffourG{The values \( \dltwu_{3} \) and \( \dltwu_{4} \) are usually very small.}
{The value \( \dltwu_{3} \) is usually very small.}
If the function \( \fs(\upsv) \) can be written in the form \( \fs(\upsv) = n \hL(\upsv) \) 
for a fixed smooth function \( h(\upsv) \) with the Hessian \( \nabla^{2} \hL(\upsv) \),
and a factor \( n \) meaning of the sample size, then \( \dltwu_{3} \asymp n^{-1/2} \).
%
%
\iffourG{\nameref{LLsT3ref} and \nameref{LLsT4ref} are local versions}
{\nameref{LLsT3ref} is a local version}
of the so-called self-concordance condition; see \cite{Ne1988} and \cite{OsBa2021}.

\Section{A linear perturbation }
For the objective function \( \fs \), 
let another function \( \fn(\upsv) \) satisfy for some vector \( \Av \)
\begin{EQA}
	\fn(\upsv) - \fn(\upsvs) 
	&=&
	\bigl\langle \upsv - \upsvs, \Av \bigr\rangle + \fs(\upsv) - \fs(\upsvs) .
\label{4hbh8njoelvt6jwgf09}
\end{EQA}
Define
\begin{EQA}
	\upsvn
	& \eqdef &
	\arginf_{\upsv} \fn(\upsv),
	\qquad
	\fn(\upsvn)
	=
	\inf_{\upsv} \fn(\upsv) .
\label{6yc63yhudf7fdy6edgehy} 
\end{EQA}
The analysis aims to evaluate the quantities \( \upsvn - \upsvs \) and \( \fn(\upsvn) - \fn(\upsvs) \).

\begin{theorem}
\label{PFiWigeneric2}
Let \( \fs(\upsv) \) be a strongly convex function with \( \fs(\upsvs) = \min_{\upsv} \fs(\upsv) \)  
and \( \IFN = \nabla^{2} \fs(\upsvs) \).
Assume \nameref{LLsT3ref} at \( \upsvs \) with \( \DFN^{2} \), \( \rrn \), and \( \dltwu_{3} \) such that 
\begin{EQA}[c]
	\DFN^{2} \leq \dmax^{2} \, \IFN ,
	\quad
	\rrn \geq \frac{3}{2} \| \DFN \, \IFN^{-1} \Av \| \, ,
	\quad
	\dmax^{2} \dltwu_{3} \| \DFN \, \IFN^{-1} \Av \| < \frac{4}{9} \, .
\label{8difiyfc54wrboer7bjfr}
\end{EQA}
Then 
\begin{EQA}[rcl]
    \| \DFN (\upsvn - \upsvs) \|
    & \leq &
    \frac{3}{2} \| \DFN \, \IFN^{-1} \Av \| 
    \, ,
    \\
    \| \DFN^{-1} \IFN (\upsvn - \upsvs + \IFN^{-1} \Av) \|
    & \leq &
    \frac{3\dltwu_{3}}{4} \| \DFN \, \IFN^{-1} \Av \|^{2} 
	\, .
\label{DGttGtsGDGm13rGa2}
\end{EQA}
Moreover, 
\begin{EQA}
    \Bigl| 2 \fn(\upsvn) - 2 \fn(\upsvs) + \| \IFN^{-1/2} \Av \|^{2} \Bigr|
    & \leq &
    \frac{\dltwu_{3}}{2} \, \| \DFN \, \IFN^{-1} \Av \|^{3} \, .
    \qquad
\label{3d3Af12DGttGa2}
\end{EQA}
\end{theorem}

\begin{remark}
\label{Rbiasgeneric}
The roles of \( \fs \) and \( \fn \) can be exchanged.
In particular, \eqref{DGttGtsGDGm13rGa2} applies with \( \IFN = \IFN(\upsvn) \) provided that 
\nameref{LLsT3ref} is also fulfilled at \( \upsvn \).
\end{remark}

\Section{Quadratic penalization}
\label{Slinquadr}
Here we discuss the case when \( \fn(\upsv) - \fs(\upsv) \) is quadratic.
The general case can be reduced to the situation with \( \fn(\upsv) = \fs(\upsv) - \| \GP \upsv \|^{2}/2 \).
To make the dependence of \( \GP \) more explicit, denote 
\( \fG(\upsv) \eqdef \fs(\upsv) - \| \GP \upsv \|^{2}/2 \),
\begin{EQA}
	\upsvs 
	&=& 
	\arginf_{\upsv} \fs(\upsv),
	\quad
	\upsvs_{\GP} = \arginf_{\upsv} \fG(\upsv) 
	=
	\arginf_{\upsv} \bigl\{ \fs(\upsv) - \| \GP \upsv \|^{2}/2 \bigr\}.
	\qquad
\label{8cvkfc9fujf6jnmcer4cd}
\end{EQA}
We study the bias \( \upsvs_{\GP} - \upsvs \) induced by this penalization.

\begin{proposition}
\label{Pbiasgeneric} 
Let \( \fG(\upsv) = \fs(\upsv) - \| \GP \upsv \|^{2}/2 \) be convex
and follow \nameref{LLsT3ref} with some \( \DFN^{2} \), \( \dltwu_{3} \), and \( \rrn \) satisfying for \( \dmax > 0 \)
\begin{EQA}[c]
	\DFN^{2} \leq \dmax^{2} \, \IFN_{\GP} \, ,
	\qquad 
	\rrn \geq 3 \bias_{\GP}/2 \, ,
	\qquad
	\dmax^{2} \dltwu_{3} \, \bias_{\GP} < 4/9 ,
\label{7fjgjgvuvy44erd52f}
\end{EQA}
where \( \bias_{\GP} \eqdef \| \DFN \, \IFN_{\GP}^{-1} \GP^{2} \upsvs \| \).
Then 
\begin{EQA}[rcl]
	\| \DFN (\upsvs_{\GP} - \upsvs) \| 
	& \leq &
	\frac{3 }{2} \, \bias_{\GP}
	\, ,
\label{odf6fdyr6e4deuewjug}
	\\
	\bigl\| \DFN^{-1} \IFN_{\GP} (\upsvs_{\GP} - \upsvs + \IFN_{\GP}^{-1} \GP^{2} \upsvs) \bigr\|
	& \leq &
	\frac{3\dltwu_{3}}{4} \, \bias_{\GP}^{2}
	\, ,
\label{11ma3eaelDebbgen}
	\\
	\Bigl| 2 \fG(\upsvs_{\GP}) - 2 \fG(\upsvs) - \frac{1}{2} \| \IFN_{\GP}^{-1/2} \GP^{2} \upsvs \|^{2} \Bigr|
	& \leq &
	\frac{\dltwu_{3}}{2} \, \bias_{\GP}^{3}
	\, .
\label{7ywsjhd7wjhjdbiui84kje}
\end{EQA}
\end{proposition}

\Section{A generic perturbation}
\label{Sgensmooth}
This section presents a generic perturbation lemma for \eqref{6yc63yhudf7fdy6edgehy}
without any assumption on the structure of the perturbation \( \fp(\upsv) - \fs(\upsv) \).
Define
\begin{EQA}
	\IFP
	\eqdef
	\nabla^{2} \fp(\upsvs),
	\quad
	\biasp
	& \eqdef & 
	\| \DFN \, \IFP^{-1} \nabla \fp(\upsvs)  \| 
	\, .
\label{fd9dfhy4ye6fuydfrerfGgen}
\end{EQA} 


\begin{theorem}
\label{Pbiasgen} 
Let a function \( \fp(\upsv) \) be strongly convex
and follow \nameref{LLsT3ref} at \( \upsvs \)
with some \( \DFN^{2} \), \( \dltwu_{3} \), and \( \rrn \) satisfying for \( \dmax > 0 \)
\begin{EQA}[c]
	\DFN^{2} \leq \dmax^{2} \hspm \IFP \, ,
	\qquad 
	\rrn \geq 3 \biasp/2 \, ,
	\qquad
	\dmax^{2} \hspm \dltwu_{3} \, \biasp < 4/9 
	\, .
\label{7fjgjgvuvy44erd52fgen}
\end{EQA}
Then 
\begin{EQA}[rcl]
	\| \DFN (\upsvp - \upsvs) \|
	& \leq &
	\frac{3}{2} \, \biasp
	\, ,
\label{odf6fdyr6e4deuewjugGgen}
	\\
	\bigl\| \DFN^{-1} \IFP (\upsvp - \upsvs + \IFP^{-1} \nabla \fp(\upsvs) ) \bigr\|
	& \leq &
	\frac{3\dltwu_{3}}{4} \, \biasp^{2}
	\, ,
\label{11ma3eaelDebbgen}
	\\
	\Bigl| 2 \fp(\upsvp) - 2 \fp(\upsvs) + \| \IFP^{-1/2} \nabla \fp(\upsvs)  \|^{2} \Bigr|
	& \leq &
	\frac{\dltwu_{3}}{2} \, \biasp^{3}
	\, .
\label{7ywsjhd7wjhjdbiui84gen}
\end{EQA}
\iffourG{
If, in addition, \( \fp(\upsv) \) satisfies \nameref{LLsT4ref} and \( \dmax^{2} \hspm \dltwu_{4} \, \biasp^{2} < \frac{1}{3} \), 
then with \( \Tensp(\uv) = \frac{1}{6} \langle \nabla^{3} \fp(\upsvs), \uv^{\otimes 3} \rangle \),
\( \nabla \Tensp = \frac{1}{2} \langle \nabla^{3} \fp(\upsvs), \uv^{\otimes 2} \rangle \), and
\begin{EQA}
	\bvp
	&=&
	- \IFP^{-1} \{ \nabla \fp(\upsvs)  + \nabla \Tensp(\IFP^{-1} \nabla \fp(\upsvs) ) \} \, ,
\label{8vfjvr43223efryfuweefG}
\end{EQA}
it holds \( \,\, \| \DFN^{-1} \IFP (\upsvp - \upsvs - \bvp) \| \leq 
    \mfrac{1}{2} \,(\dltwu_{4} + 2 \dmax^{2} \hspm \dltwu_{3}^{2}) \biasp^{3} 
 \),
\begin{EQA}
    && \nquad
    \Bigl| \fp(\upsvp) - \fp(\upsvs) 
    + \frac{1}{2} \| \IFP^{-1/2} \nabla \fp(\upsvs)  \|^{2} + \Tensp(\IFP^{-1} \nabla \fp(\upsvs) ) \Bigr|
    \\
    & \leq &
    \frac{\dltwu_{4} + 4 \dmax^{2} \hspm \dltwu_{3}^{2}}{8} \biasp^{4} 
    + \frac{\dmax^{2} \, (\dltwu_{4} + 2 \dmax^{2} \hspm \dltwu_{3}^{2})^{2} }{4} \, \biasp^{6} \, .
    \qquad
\label{3d3Af12DGttGa4bgen}
\end{EQA}
}{}
\end{theorem}